\documentclass[12pt]{article}
\usepackage{graphicx} 

\usepackage{amsmath}
\usepackage[colorlinks=true, allcolors=blue]{hyperref}
\usepackage{subfig}
\usepackage{amssymb}
\usepackage{bm}
\usepackage{float} 
     
\usepackage{algorithm}
\usepackage{algorithmic}
\usepackage[margin=1in]{geometry}

\usepackage{makecell}

\providecommand{\keywords}[1]
{
  \small	
  \textbf{\textit{Keywords---}} #1
}

\newtheorem{theorem}{Theorem}[section] 
\newtheorem{remark}[theorem]{Remark} 
\newtheorem{definition}{Definition}[section]

\usepackage{multirow} 
\usepackage{booktabs}

\title{ Post-Processing Reduced-Order Models for Transport-Dominated Problems by Gegenbauer Reconstruction}
\author{Lei Yan
    \thanks{School of Mathematical Sciences, University of Science and Technology of China, Hefei, Anhui 230026, China.
    E-mail: {\tt lei\_yan@mail.ustc.edu.cn}.}
\and Yan Jiang
    \thanks{School of Mathematical Sciences, University of Science and Technology of China, Hefei, Anhui 230026, China.
    Email: {\tt jiangy@ustc.edu.cn}.}
\and Chi-Wang Shu
    \thanks{Division of Applied Mathematics, Brown University, Providence, RI 02912, Email: {\tt chi-wang\_shu@brown.edu}. }
    }
\date{}

\begin{document}

\maketitle

\begin{abstract}
In this paper, we develop numerical techniques for post-processing of data-driven reduced-order models (ROMs) for transport-dominated problems. 
For transport-dominated problems, we encounter not only the challenge of slow decay in the Kolmogorov n-width, but also, due to the global nature of the reduced basis, the on-line computation of data-driven ROMs can result in unphysical oscillations when solving for solutions with shocks or sharp gradients. This is closely related to the Gibbs phenomenon, which is a well-known phenomenon in spectral approximations. This paper aims at addressing this issue using post-processing techniques from spectral methods by Gegenbauer polynomial reconstruction.

The idea of Gegenbauer reconstruction is to re-project the solution in each interval of analyticity by Gegenbauer polynomial basis. It has been shown to offer spectral accuracy for spectral basis including Fourier and orthogonal polynomial expansions. We adopt this approach to ROM data.
In particular, we consider three types of commonly used linear and nonlinear ROMs: Proper Orthogonal Decomposition (POD) based Galerkin ROM, Operator Inference (OpInf) and nonlinear manifold ROM based on convolutional auto-encoders (CAE). We show that post-processing by Gegenbauer reconstruction is effective for all three ROMs in the sense that it eliminates unphysical oscillations for solutions with discontinuities. The Gegenbauer reconstruction can be easily applied  in 1D once a discontinuity detector is used. We further propose a detailed reconstruction procedure for 2D problems combining line-by-line reconstruction with respect to each dimension. We demonstrate that Gegenbauer reconstruction can  reduce  the numerical error by one or two orders of magnitude for inviscid problems and the results significantly outperform standard method by total variation regularization in numerical accuracy and sharp resolution of the discontinuity. 
\end{abstract}


\section{Introduction}

This paper introduces a post-processing technique for computing accurate solutions from data-driven reduced-order models (ROMs) for transport-dominated problems. It is well known that numerical solutions to conservation laws may develop shock and other complex structures, motivating decades of research on non-oscillatory numerical methods \cite{hesthaven2017numerical}. Those structures present great challenges for ROMs  \cite{peherstorfer2022breaking}. The goal of this paper is to design a numerical procedure to eliminate the spurious oscillations of solutions from ROMs for shocks. 

The main insight of our work comes from the fact that most commonly used data-driven ROMs are \emph{global} methods, i.e. the basis functions or the representation of the solutions are defined on the whole  domain. 
This characteristic establishes a close and fundamental connection between reduced-order models and spectral methods such as Chebyshev or Fourier expansions. 
On function approximation level, the spectral expansion of discontinuous functions gives rise to oscillatory behavior of the finite sum truncation, called the Gibbs phenomenon \cite{gottlieb1997gibbs, tadmor2007filters}.  
Analogously, in reduced-order models, the combination of global basis approximation and severe modal truncation results in similar unphysical oscillations\cite{yu2022model}, a challenge widely observed in both proper orthogonal decomposition (POD)-based methods \cite{berkooz1993proper} and various nonlinear model reduction techniques \cite{cagniart2017model, lee2020model, maulik2021reduced} when dealing with strong shocks, insufficient resolution, or stochastic effects.
Building on this shared mathematical context, this paper investigates the applicability of an oscillation-eliminating post-processing technique originally developed and proven highly effective for spectral methods, and to determine whether it can perform robustly on the diverse family of data-driven ROM solutions.

However, unlike in spectral methods where orthogonal basis functions are known, the data-driven basis functions used in reduced-order model are typically problem-dependent, which makes the manifestation of Gibbs-type oscillations more complex and their removal more challenging. The oscillations originate not only from spatial discretization error but also from the intrinsic limitations of representing transport-dominated, high-dimensional dynamics with a low-dimensional global subspace. Consequently, post-processing techniques designed in spectral frameworks require careful adaptation for reduced-order model data.

To address spurious oscillations, numerous strategies have been developed to enhance the stability and accuracy of reduced-order models.
A direct approach is to apply filtering techniques as a post-processing step, as commonly used in both spectral methods and reduced-order models for numerical stabilization \cite{sirisup2004spectral,hesthaven2007spectral, bergmann2009enablers, wells2017evolve}. However, such methods typically degrade numerical accuracy near discontinuities and tend to overly smear sharp features such as shocks.
A more fundamental line of research aims to modify the ROM formulation itself to better accommodate the translational characteristics of transport-dominated problems. Representative examples include transforming the governing equations into a moving reference frame using a Lagrangian framework, with extensions based on hodograph transformations to improve shock resolution \cite{lu2020lagrangian,lu2021dynamic}.
Other approaches preprocess the snapshot data to align coherent structures prior to model reduction, such as in the shifted Proper Orthogonal Decomposition (sPOD) \cite{reiss2018shifted} and shifted Operator Inference (sOpInf) \cite{issan2023predicting}.
A distinct class of methods abandons the use of a fixed global basis altogether, instead employing local or time-dependent bases that dynamically adapt to evolving flow features during online simulation \cite{peherstorfer2020model}.

Although these advanced strategies have demonstrated strong performance for transport-dominated problems, they typically modify the underlying formulation or training procedure of reduced-order models. As a result, they may introduce additional computational overhead and tend to be problem-dependent, requiring problem-specific tuning or reformulation. In contrast, the proposed approach in our work is a general, non-intrusive post-processing technique that can be applied to a wide range of reduced-order models without altering the underlying model structure or solution procedure.
Specifically, we mitigate the spurious oscillations from reduced-order model by employing Gegenbauer reconstruction technique \cite{gottlieb1997gibbs}. 
The Gegenbauer reconstruction removes Gibbs phenomenon completely for spectral methods, i.e., it can obtain exponential accuracy in the maximum norm in any interval of analyticity.
The main idea is that one can show the knowledge of the expansion coefficients is sufficient for obtaining the point values of a piecewise smooth function, with the same order of accuracy as in the smooth case. This is done by using the finite expansion series to construct a different, rapidly convergent approximation. The application of Gegenbauer reconstruction is simple and computationally effective. It is a post-processing step done only once by  performing an edge detection, then re-projecting the solution separately on each subdomains of analyticity.
A series of works as mentioned in \cite{gottlieb1997gibbs,tadmor2007filters} showed great performance of Gegenbauer reconstruction for spectral methods for hyperbolic PDEs and image processing tasks \cite{archibald2002method,gelb2004parameter}. Recent work \cite{kawai2022gegenbauer} demonstrated that the applicability of Gegenbauer reconstruction towards multi-dimensional uncertainty propagation. 
Building on this Gegenbauer polynomial framework, the Bayesian Spectral Reprojection (BSR) and Generalized BSR (GBSR) methods \cite{li2025bayesian} further enabled robust reconstruction of spectral data corrupted by random noise.

In this paper, we test the Gegenbauer reconstruction methods on various reduced-order models, including POD Galerkin (G-ROM), Operator Inference (OpInf) \cite{peherstorfer2016data}, and  
nonlinear ROMs based on convolutional auto-encoders (CAE) \cite{kashima2016nonlinear, lee2020model, maulik2021reduced}. We test several benchmark problems in 1D and 2D, and find that the Gegenbauer reconstruction is successful in completely removing numerical oscillations. When the parameters in the reconstruction are chosen appropriately, it can reduce the numerical errors by one to two orders of magnitude. When comparing the results with standard techniques such as total variation regularization or ridge regression, Gegenbauer reconstruction shows clear advantages in  numerical accuracy  and sharp resolution of the discontinuity.

The rest of the paper is organized as follows. Section~\ref{Reduced order models} reviews the reduced-order models we considered in this paper, including POD-Galerkin ROM, OpInf and ROMs based on CAE.  Section~\ref{Enhancing ROMs} introduce Gegenbauer reconstruction for ROMs in 1D and 2D. In Section~\ref{Numerical Results}, we perform numerical tests on 1D and 2D transport-dominated problems. Finally, in Section~\ref{Conclusions and Future Work}, we present conclusions and future research directions.

\section{Review of Reduced-Order Models}
\label{Reduced order models}
This section first presents the general form of the problem, followed by an introduction and a brief review to three types of reduced-order models, including POD-based Galerkin ROM, OpInf and ROMs based on CAE. We note that data-driven ROMs \cite{brunton2022data} constitute a vast research area. The goal is not to review the whole discipline, but rather to provide a short introduction to those ROMs that we tested in this paper.

ROMs are designed to reduce the computational burden of parametric dynamical systems. 
We consider ordinary differential equations (ODEs) defined on the time domain $[t_i,t_f]$, with $t_i$ representing the initial time and $t_f$ the final time,
\begin{equation}
    \frac{d \bm q}{d t}(t;\mu) = \bm f(\bm q,t;\mu)\qquad\bm q(t_i;\mu) = \bm q_{0},
    \label{problem_setup}
\end{equation}
where $\bm q(t;\mu) \in \mathbb R^{n}$ is the $n$-dimensional vector of state variables at time $t$, $\mu\in \mathcal{D}\subset \mathbb{R}$ is a parameter, 
$\bm q_{0}$ is a specified initial condition, $f: \mathbb R^{n}\times [t_i,t_f]\times \mathcal{D} \to \mathbb R^n$ is a linear or nonlinear function defined as the time evolution of the system state.
Equation \eqref{problem_setup} encompasses a wide range of models, including those for partial differential equations (PDEs) after spatial discretizations have been employed. In this paper, we focus exclusively on equations that are derived from transport-dominated problems, in particular for which the solutions contain shocks. In the following, we discuss the setup of ROMs for \eqref{problem_setup}. All of the ROMs here adopt the offline-online decomposition, 
in which the offline stage constructs the reduced model by generating training data and computing reduced bases,  while the online stage evaluates the resulting reduced system at dramatically lower cost.

\subsection{\label{G-ROM} POD-Galerkin ROM}
Technically, POD \cite{berkooz1993proper} is primarily a mode-extraction technique; when combined with Galerkin projection, it provides a well-known approach for model reduction by identifying the most important modes and constructing a low-dimensional description of the system’s dynamics.
The method begins with the construction of the snapshot matrix. In the offline stage, we obtain the training data that is accessible over the time interval $[t_i, t_{tr}]$, 
where $t_i\leq t_1< \cdots< t_{n_t}\leq t_{tr}$, with $t_{tr} < t_f$ denoting the end of the training time horizon. Furthermore, 
a set of $n_\mu$ parameter samples $\{\mu_1, \cdots, \mu_{n_\mu}\}$ is drawn from the parameter domain $\mathcal{D}$.
For each selected parameter $\mu$, during the training phase, we obtain the high-dimensional state solution in $n_t$ time instances by solving the high-fidelity problem~\eqref{problem_setup}, which serves as training snapshots. Let $\bm q_{i,j} = \bm q(t_i; \mu_j)$ denote the snapshot corresponding to the $i$-th time instance and the $j$-th parameter, the snapshot matrix can then be expressed as:
\begin{equation}
    \bm{\mathcal{Q}} = [\bm q_{1,1},\cdots,\bm q_{n_t,1},\cdots,\bm q_{1,n_\mu},\cdots,\bm q_{n_t,n_\mu}]\in \mathbb R^{n\times (n_t\, n_\mu)}.
\end{equation}

The POD basis functions can be obtained by finding the eigenvectors to the matrix $\bm{\mathcal{Q}}\bm{\mathcal{Q}}^T$ that correspond to the largest $r$ eigenvalues. We denote those bases by \( \bm{\varphi}_1, \cdots, \bm{\varphi}_r  \in \mathbb R^n \), which is computed during the offline phase.

The Galerkin-ROM (G-ROM) approximates the state variable by
\begin{equation}
    \bm q_r(t;\mu) \approx \sum_{i=1}^{r} \hat{q}_i(t;\mu) \bm \varphi_i,
\end{equation}
where $\{\hat{q}_i(t)\}_{i=1}^r $ are the desired time-dependent coefficients representing the trajectories. By a Galerkin formulation, we obtain
\begin{equation}
    \Big( \frac{d\bm q_r(t;\mu)}{d t} , \bm \varphi_i \Big) 
    = \Big(\bm f (\bm q_r,t;\mu) , \bm \varphi_i \Big), 
    \qquad \forall \ i = 1, \cdots, r,
    \label{G_rom}
\end{equation}
where $(\cdot, \cdot)$ is the inner product of vectors.
Since the POD basis functions are mutually orthogonal, equation~\eqref{G_rom} can be transformed into a system of ordinary differential equations with respect to $\{\hat{q}_j(t)\}_{j=1}^r$, 
\begin{equation}
    \frac{d \hat{q}_i(t;\mu)}{dt} = \Big(\bm f (\bm q_r,t;\mu) ,\bm \varphi_i\Big), \qquad  \ i = 1, \cdots, r,
    \label{eq_GROM}
\end{equation}
with initial conditions
\begin{equation}
    \hat{q}_i(0) = (\bm q_0,\bm \varphi_i), \qquad  \ i =1,\cdots,r.
\end{equation}

Equation~\eqref{eq_GROM} is evolved \emph{online}. Since the G-ROM is evolving only $r$ unknowns, it yields substantial computational savings compared to the full-order methods (FOMs).
If the term $\bm f(\bm q(t), t)$ contains nonlinearities, we note that a hyper-reduction technique, e.g. Discrete Empirical Interpolation Method (DEIM) \cite{chaturantabut2010nonlinear}, will be employed to further enhance computational efficiency. In addition, filtering and spectral viscosity \cite{sirisup2004spectral,  bergmann2009enablers, wells2017evolve} have been developed to stabilize the ROMs, which is particularly important for nonlinear equations.

\subsection{Operator inference}

G-ROM belongs to the class of projection-based reduced-order models. In this paper, we also consider purely data-driven ROMs, which constitute an active area of research. Among various types of data-driven reduced-order models, in this work, we focus on the approach of OpInf as a representative example. 
OpInf \cite{peherstorfer2016data} is a technique for learning reduced-order models with polynomial structure from data. It approximates and analyzes dynamic systems by inferring operators from observed data that can capture the underlying system dynamics.

Following the same notations as in the previous subsection, we assume
\begin{equation}
    \bm q_{opinf}(t;\mu) \approx \sum_{i=1}^{r} \hat{q}_i(t;\mu) \bm \varphi_i,
\end{equation}
and $\hat{\bm q}=(\hat{q}_1, \ldots \hat{q}_r)^T.$
For the time-continuous or semi-discrete OpInf, the reduced operators of the reduced ODE system are learned by solving a linear least-squares problem. We assume the ROM to have  a quadratic form
\begin{equation}
    \frac{d \hat{\bm q}}{d t} = \hat{\bm A}_{sd} \hat{\bm q} + \hat{\bm H}_{sd} \bigl(\hat{\bm q} \otimes \hat{\bm q}\bigr),
    \label{eq_opinf_c}
\end{equation}
where $\hat{\bm A}_{sd}\in \mathbb R^{r\times r}$ and $\hat{\bm H}_{sd}\in \mathbb R^{r\times r^2}$ are to be trained based on the data $\hat{\bm{\mathcal Q}}= \bm U_r^T \bm{\mathcal{Q}} \in \mathbb{R}^{r\times  n_t},$ and \(\bm U_r = [\bm \varphi_1, \cdots, \bm \varphi_r]\in \mathbb{R}^{n\times r}\). 
It infers the reduced operators that optimally align with the projected snapshot data by minimizing the residuals through the solution of a linear least squares problem
\begin{equation}
    \mathop{\mathrm{arg\,min}}_{\hat{\bm A}_{sd},\hat{\bm H}_{sd}} \Vert \hat{\bm{\mathcal Q}}^T \hat{\bm A}_{sd}^T +(\hat{\bm{\mathcal Q}}\otimes \hat{\bm{\mathcal Q}})^T\hat{\bm H}_{sd}^T-\frac{d}{dt} \hat{\bm{\mathcal Q}}^T\Vert_F^2+\lambda_1 \Vert\hat{\bm A}_{sd}\Vert_F^2 +\lambda_2 \Vert \hat{\bm H}_{sd}\Vert_F^2,
    \label{sd_opinf}
\end{equation}
where $\Vert \cdot \Vert_F$ denotes the Frobenius norm. In order to prevent overfitting and to address potential model misspecification as well as numerical noise in the computed time derivatives, Tikhonov regularization is applied to~\eqref{sd_opinf}, with $\lambda_1$ and $\lambda_2$ serving as the scalar regularization hyperparameters.

For problems with large time steps, performing accurate calculations of the derivative terms of the projected snapshots, $\frac{d}{dt} \hat{\bm{\mathcal Q}}$, as described in~\eqref{sd_opinf}, can be challenging. In such cases, we can consider  time-discrete version of OpInf obtained by solving
\begin{equation}
    \hat{\bm q} [k+1] = \hat{\bm A}_d \hat{\bm q}[k] + \hat{\bm H}_d (\hat{\bm q}\otimes \hat{\bm q}),
    \label{eq_d_opinf}
\end{equation}
where $\hat{\bm q}[k]\in \mathbb R^r$ denotes the discrete reduced state at time $t_k$. The reduced operators in the discrete Operator Inference model can be obtained by solving
\begin{equation}
    \mathop{\mathrm{arg\,min}}_{\hat{\bm A}_{d},\hat{\bm H}_{d}} \Vert \hat{\bm {\mathcal Q}}_1^T \hat{\bm A}_d^T + (\hat{\bm{\mathcal Q}}_1\otimes \hat{\bm{\mathcal Q}}_1)^T\hat{\bm H}_{d}^T - \hat{\bm {\mathcal Q}}_2\Vert_F^2+\lambda_1 \Vert\hat{\bm A}_{d}\Vert_F^2 +\lambda_2 \Vert \hat{\bm H}_{d}\Vert_F^2,
    \label{d_opinf}
\end{equation}
where
\begin{equation}
    \hat{\bm {\mathcal Q}}_1 = [\hat{\bm q}_1,\cdots,\hat{\bm q}_{n_t-1}]\in \mathbb{R}^{r\times (n_t-1)}\qquad \hat{\bm {\mathcal Q}}_2= [\hat{\bm q}_2,\cdots,\hat{\bm q}_{n_t}]\in \mathbb{R}^{r\times (n_t-1)}.
\end{equation}

For more details on OpInf and the calculation of the optimization problem with the Tikhonov regularization mentioned above, one can refer to \cite{peherstorfer2016data, mcquarrie2021data, qian2022reduced}.

\subsection{Nonlinear ROMs via convolutional auto-encoders}
The G-ROM and OpInf project the governing equations onto a linear subspace of the original state space. Such linear ROMs need many modes when the problem exhibits a slowly decaying Kolmogorov \(n\)-width \cite{ohlberger2015reduced}. In recent years, nonlinear manifold ROMs have gained attention in the ROM community. In this work, we consider nonlinear ROMs based on CAE, namely CAE-Galerkin \cite{kashima2016nonlinear, lee2020model} and CAE-LSTM (Long Short term memory) \cite{maulik2021reduced}. The former is a projection-based ROM, and the later is a purely data-driven ROM.

\subsubsection{CAE-Galerkin}

The auto-encoder, $ \bm h: \bm q \mapsto \tilde{\bm q} $, is a feedforward neural network designed to learn identity mapping, which consists of an encoder and a decoder. The auto-encoder is a well-established approach for nonlinear dimension reduction \cite{demers1992non}. The encoder, $\bm h_{enc}: \bm q \mapsto \hat{\bm q}$ with $\bm h_{enc}: \mathbb{R}^n \to \mathbb{R}^r$, maps the high-dimensional vector $\bm q$ to a low-dimensional code $\hat{\bm q}$. The decoder, $\bm h_{dec}: \hat{\bm q} \mapsto \tilde{\bm q}$ with $\bm h_{dec}: \mathbb{R}^r \to \mathbb{R}^n$, reconstructs the high-dimensional vector $\tilde{\bm q}$ from the low-dimensional code $\hat{\bm q}$. Thus, the resulting auto-encoder takes the form
\begin{equation}
    \bm h : \bm q \mapsto \bm h_{dec} \circ \bm h_{enc} (\bm q).
\end{equation}
The goal is for $\bm h(\bm q)$ to approximate $\bm q$ as closely as possible, so the loss function is defined as
\begin{equation}
    Loss(\bm{\mathcal{Q}};\bm \theta) = \sum_{\bm q  \in \bm{\mathcal{Q}}} \Vert \bm h(\bm q) - \bm q\Vert_2^2,
\end{equation}
where $\bm{\mathcal{Q}}$ is the snapshot matrix.

In a feedforward network, each layer typically corresponds to a vector or tensor, whose values are computed by applying an affine transformation followed by a nonlinear activation function to the previous layer. An encoder with \( n_L \) layers takes the following form:
\begin{equation}
    \bm h_{enc}: (\bm q;\bm \Theta_{enc}) \mapsto \bm h_{n_L} (\cdot;\bm \Theta_{n_L})\circ \bm h_{n_L-1} (\cdot,\bm \Theta_{n_L-1})\circ \cdots \circ \bm h_1 (\bm q;\bm \Theta_1),
\end{equation}
where $\bm h_i (\cdot;\bm \Theta_{i}): \mathbb{R}^{p_{i-1}}\to \mathbb{R}^{p_i}, i = 1,\cdots, n_L$ and $p_i$ denotes the dimensionality of the output at layer $i$. The input dimension is $p_0 = n$ , and the final layer generates a dimension of $p_{n_L} = r$. The nonlinear activation function is utilized on the weights, biases, and outputs of the preceding layer, progressively from the first layer to the \( n_L\)-th layer, i.e.,
\begin{equation}
    \bm h_i: (\bm q;\bm \Theta_i)\mapsto \phi_i(h_i (\bm q;\bm \Theta_i)),
\end{equation}
where $\phi_i$ is an element-wise nonlinear activation function. Since we are constructing a convolutional auto-encoder, $h_i$ represents the corresponding convolution operator and $\bm \Theta_i$ provides the weights of the convolutional filter.

The decoder with $\overline{n}_L$ layers also corresponds to a feedforward network, which takes the form
\begin{equation}
    \bm h_{dec}: (\hat{\bm q};\bm \Theta_{dec}) \mapsto \overline{\bm h} _{\overline n_L} (\cdot;\overline{\bm \Theta}_{\overline n_L})\circ \overline{\bm h}_{\overline n_L-1} (\cdot,\overline{\bm \Theta}_{\overline n_L-1})\circ \cdots \circ \overline{\bm h}_1 (\hat{\bm q};\overline{\bm \Theta}_1),
\end{equation}
where $\overline{\bm h}_i (\cdot;\overline{\bm \Theta}_{i}): \mathbb{R}^{\overline p_{i-1}}\to \mathbb{R}^{\overline p_i}, i = 1,\cdots, \overline n_L$ and $\overline p_i$ denotes the dimensionality of the output at layer $i$. Similarly, the nonlinear activation is applied to the output of the previous layer, such that
\begin{equation}
    \overline{\bm h}_i: (\bm q;\overline{\bm \Theta}_i)\mapsto \phi_i (\overline h_i (\bm q;\overline{\bm \Theta}_i)),
\end{equation}
where $\overline{h}_i$ corresponds to a convolution operator, and $\overline{\bm \Theta}_i$ provides the transposed convolutional filter weights.

Assuming that the CAE model has been trained, the approximation of the true solution can be expressed as
\begin{equation}
    \tilde{\bm q}(t;\mu) = \bm q_{ref} (\mu) + \bm h_{dec} (\hat{\bm q}(t;\mu)),
\end{equation}
where $\tilde{\bm q} \in \mathbb{R}^n$ and $\hat{\bm q} \in \mathbb{R}^r$.
The reference solution $\bm q_{ref}\in\mathbb{R}^n $ can be computed from the initial conditions as:
\begin{equation}
    \bm q_{ref} = \bm q_0 - \bm h (\bm q_0).
\end{equation}

From application of the chain rule, we can obtain
\begin{equation}
    \frac{d\tilde{\bm q}}{dt} = J(\hat{\bm q}(t;\mu)) \frac{d\hat{\bm q}}{dt}, 
\end{equation}
where $J:\hat{\bm q}\mapsto \frac{d\bm h_{dec}}{d\hat{\bm q}} $ with $J:\mathbb{R}^r \to \mathbb{R}^{n\times r}$ denoting the Jacobian of the decoder.
The residual function for Problem \eqref{problem_setup} is defined as:
\begin{equation}
    \bm r: (\bm q, \bm \xi,t;\mu)\mapsto \bm q - \bm f(\bm \xi,t;\mu).
\end{equation}

Our goal is to achieve $\tilde{\bm q} \to \bm q$ and $\frac{d\tilde{\bm q}}{dt} \to \frac{d\bm q}{dt}$. By substituting $\tilde{\bm q}$ into problem~\eqref{problem_setup} and minimizing the \( L^2 \)-norm of the residual function, the loss function can be expressed as:
\begin{equation}
    \frac{d \hat{\bm q}}{dt}(t;\mu) = \arg\min \limits_{\hspace{-24pt}\hat{\bm v} \in \mathbb{R}^r} \Vert \bm r(J(\hat{\bm q}(t;\mu))\hat{\bm v},\bm q_{ref}+\bm h_{dec}(\hat{\bm q}(t;\mu)),t;\mu)\Vert_{L^2}^2,
\end{equation}
with initial condition $\hat{\bm q}_0 = \bm h_{enc}(\bm q_0)$. If the Jacobian matrix $J(\hat{\bm q}(t;\mu))$ has full column rank, then the above residual minimization problem has a unique solution:
\begin{equation}
    \frac{d \hat{\bm q}}{dt}(t;\mu) = J(\hat{\bm q}(t;\mu))^+ f(\bm q_{ref}(\mu)+\bm h_{dec}(\hat{\bm q}(t;\mu))),
    \label{CAE_Galerkin}
\end{equation}
where the superscript $+$ denotes the Moore–Penrose pseudoinverse, and Equation~\eqref{CAE_Galerkin} represents the CAE-Galerkin ROM.

\begin{remark}
    The above description primarily targets time-continuous problems. For time-discrete cases, article \cite{lee2020model} presents an alternative solution approach using the least-squares Petrov–Galerkin (LSPG) projection.
\end{remark}

\subsubsection{CAE-LSTM}
In the CAE-LSTM method, the training of the convolutional auto-encoder is identical to that in CAE-Galerkin and therefore will not be repeated here. Unlike the method used for processing low-dimensional data in CAE-Galerkin, the CAE-LSTM network employs the  LSTM  network \cite{hochreiter1997long}, which is a specialized type of recurrent neural network (RNN) for modeling sequential data \cite{al2024rnn}, to learn the relationships between the training data and process the low-dimensional data. Therefore, it is a data-driven ROM. 
The CAE-LSTM model is designed to capture the temporal dynamics of the reduced-order states in the latent space obtained through the CAE. The training process involves constructing a time series dataset from the encoded snapshots and input of parameters. A sliding time window approach is used to prepare the input and output sequences, ensuring that the model effectively learns the temporal dependencies.  The implementation detail and the network structure we used in this paper follows from  \cite{maulik2021reduced} and we skip the details for brevity.

\section{Post-processing ROM solutions by Gegenbauer   reconstruction}
\label{Enhancing ROMs}

The goal of this work is to address the issues of spurious oscillations of ROMs for convection dominated problems. The ROMs, while effective in cost reduction, come at a cost of employing global basis. Resolution of Gibbs phenomenon is studied in spectral methods \cite{gottlieb1997gibbs,tadmor2007filters}. Techniques including filters and Gegenbauer reconstruction techniques \cite{gottlieb1997gibbs}, the Fourier-Pad\'e approximation \cite{min2007fourier}, and the singular Pad\'e approximation \cite{driscoll2001pade} are designed for this purpose among others. In this paper, we use the Gegenbauer reconstruction technique. Below, we will review the Gegenbauer reconstruction for spectral methods, and discuss its adaptation to ROM solutions in 1D and 2D.

\subsection{Gegenbauer polynomial  reconstruction for spectral methods}
In the following, we will review Gegenbauer polynomial reconstruction for spectral methods, which removes the Gibbs phenomenon completely, by re-expanding a slowly converging global expansion using a different set of basis functions.

To introduce the concept, we assume that $f(x)$ is analytic on the subinterval $[a,b]$, and denote its finite continuous expansion by
\begin{equation}
\label{eq:spect}
    f_N(x) = \sum_{k=-N}^{N} (f,\Psi_k)\Psi_k(x),
\end{equation}
where $\Psi_k(x)$ are the spectral (Fourier/orthogonal polynomial) bases. We can re-expand $f_N(x)$, by an orthonormal family of functions $\{\Phi_k^\lambda\}$ defined in an interval of analyticity $[a,b]$,
\begin{equation}
    f_N^\lambda(x) = \sum_{k \ge 0} <f_N,\Phi_k^\lambda>_\lambda \Phi_k^\lambda(\xi(x)),\qquad \xi(x) = -1 + 2\frac{x-a}{b-a}.
\end{equation}
Here $<\cdot, \cdot>_\lambda$ refers to the weighted inner product on $[-1, 1]$ as
$$<f,g>_{\lambda} = \int_{-1}^1 (1-\xi^2)^{\lambda-\frac{1}{2}} f(\xi) g(\xi) d\xi $$
with $\lambda\geq0$. 
The bases $\Phi_k^\lambda$ are chosen as the Gibbs complementary to the original bases as defined below.
 
\begin{definition}[Gibbs complementary \cite{hesthaven2007spectral}]
    The new basis functions $\{\Phi_k^\lambda (\xi)\}$ are \textbf{Gibbs complementary} to the original basis functions $\{\Psi_k(x)\}$ if the following three conditions are satisfied:
    \begin{enumerate}
        \item Orthonormality: For any fixed $\lambda$,
        \begin{equation}
            <\Phi_k^\lambda (\xi),\Phi_l^\lambda(\xi)>_\lambda = \delta_{kl}.
        \end{equation}
        \item Spectral convergence: The expansion of a function $g(\xi)$ which is analytic in $\xi \in [-1,1]$, in the basis $\{\Phi_k^\lambda (\xi)\}$, converges exponentially fast with $\lambda$.
        \item The Gibbs condition: There exists a number $\alpha<1$ and $\beta < 1$ such that if $\lambda = \beta N$, then
        \begin{equation}
            |<\Phi_l^\lambda (\xi),\Psi_k(x(\xi))>_\lambda| \cdot \max_{-1\le \xi\le 1} |\Phi_l^\lambda(\xi)|\le (\frac{\alpha N}{k})^\lambda,
        \end{equation}
        for any $k>N, l<\lambda$.
    \end{enumerate}
\end{definition}

The new expansion $f_N^\lambda$ converges exponentially fast in the interval \([a,b]\) provided that the family $\{\Phi_k^\lambda(\xi)\}$ is Gibbs complementary. A detailed proof of this result is presented in \cite{hesthaven2007spectral}, where it is further demonstrated that the Gegenbauer polynomials $\{C_k^\lambda\}$ form a Gibbs complementary set with respect to the Fourier, Legendre and general Gegenbauer bases. 

We now review the definition of Gegenbauer polynomials.
\begin{definition}
    The Gegenbauer polynomial \( C_n^\lambda(x) \), for \( \lambda \ge 0 \), is the polynomial of degree \( n \) that satisfies
    \begin{equation*}
        \int_{-1}^{1} (1-x^2)^{\lambda - \frac{1}{2}} C_k^\lambda(x) C_n^\lambda(x) \, dx = h_n^\lambda\, \delta_{kn}, \quad \forall \ k.
    \end{equation*}
    The normalization constant \( h_n^\lambda \) associated with the Gegenbauer polynomial is defined by
    \begin{equation}
    \label{eq:normalc}
        h_n^\lambda = \int_{-1}^1 (1-x^2)^{\lambda - \frac{1}{2}} (C_n^\lambda(x))^2 \, dx = \pi^{\frac{1}{2}} C_n^\lambda(1) \frac{\Gamma(\lambda + \frac{1}{2})}{\Gamma(\lambda) (n+\lambda)},
    \end{equation}
        where
        \begin{equation}
        C_n^\lambda(1) = \frac{\Gamma(n+2\lambda)}{n! \Gamma(2\lambda)}.
    \end{equation}
\end{definition}

We can now  calculate the first $m+1$ Gegenbauer coefficients, 
\begin{equation}
    \hat{g}^\lambda_\epsilon(k) = \frac{1}{h_k^\lambda} \int_{-1}^1 (1-\xi^2)^{\lambda-\frac{1}{2}}f_N(x(\xi)) C_k^\lambda(\xi)d\xi, \quad k=0,\cdots, m,
\end{equation}
and finally, the Gegenbauer series can be defined as:
\begin{equation}
    f_N^{\lambda,m}(x) =  \sum_{k=0}^m \hat{g}_\epsilon^\lambda(k) C_k^\lambda(\xi).
\end{equation}

The analysis in \cite{gottlieb1997gibbs} demonstrated if $\lambda, m$ is chosen to be proportionate to $N,$ then $f_N^{\lambda,m}$ converges spectrally to $f$ on the interval $[a,b]$ uniformly.
When applied to spectral methods for PDEs, Gegenbauer polynomial reconstructions show excellent performance in eliminating oscillations and recovering spectral accuracy. The main intuition is the spectral coefficient in \eqref{eq:spect} contains enough (spectrally accurate) information about the solution, and the Gegenbauer reprojection can extract this information by a simple post-processing step.

\subsection{Gegenbauer polynomial  reconstruction for ROMs}

In this subsection, we present the details of applications of the Gegenbauer polynomial reconstruction to ROMs, addressing the 1D and 2D cases separately.

\subsubsection{1D case}
Suppose we have obtained the ROM solution at the final time as $u_{ROM}(x)$. The reconstruction procedure can be structured as follows:

\paragraph{Step 1: Discontinuity detection and domain partitioning.}
The first step is to use a discontinuity detector on $u_{ROM}(x)$ to partition the domain into subdomains, on each of which, the solution is smooth. 
In this work, we use Sobel filter \cite{sobel19683x3} for discontinuity detection. We note that for periodic problems, we treat the left-end and right-end points as the same. Therefore, we can merge the left-most and right-most subdomains into one interval using periodic boundary condition.

\paragraph{Step 2: Local variable definition.}
Suppose the subdomain is defined as $[a,b],$ we define the local variable $\xi \in [-1,1]$ as:
\begin{equation}
    \xi = \frac{x-\delta}{\epsilon},\qquad \delta=\frac{a+b}{2},\quad \epsilon = \frac{b-a}{2}.
\end{equation}

\paragraph{Step 3: Computation of Gegenbauer coefficients.}
The Gegenbauer coefficients are computed from the ROM solution as
\begin{equation}
    \hat{g}^\lambda_\epsilon(k) = \frac{1}{h_k^\lambda} \int_{-1}^1 (1-\xi^2)^{\lambda-\frac{1}{2}} u_{ROM}(x(\xi)) C_k^\lambda(\xi)d\xi,
    \label{eq_ROM_gegen}
\end{equation}
where $h_k^\lambda$ is the normalization constant defined in \eqref{eq:normalc}.
Equation \eqref{eq_ROM_gegen} is computed by a numerical quadrature, e.g.  Simpson's rule.

\paragraph{Step 4: Reconstruction of the Gegenbauer series.}
Finally, we construct the Gegenbauer series as:
\begin{equation}
    u_{ROM}^G(x) = \sum_{k=0}^m \hat{g}^\lambda_\epsilon(k) C_k^{\lambda}(\xi).
\end{equation}

The computational cost of the 1D Gegenbauer reconstruction is minor compared to the offline and online stages of the ROM procedure. 
It can be summarized as follows. 
The discontinuity detection using the Sobel filter involves a single local operation over the grid and thus costs $\mathcal{O}(N)$. 
For each smooth subdomain, the computation of $m$ Gegenbauer coefficients through numerical quadrature involves $\mathcal{O}(m Q_i)$ operations, where $Q_i$ is the number of quadrature (or reconstruction) points in the $i$-th subdomain. 
Reconstructing the Gegenbauer series is performed independently on each subdomain, also costing $\mathcal{O}(m Q_i)$ per subdomain. 
Summing over all subdomains and taking $\sum_i Q_i \approx N$, the overall complexity of the 1D reconstruction is
$T_{1\mathrm{D}} = \mathcal{O}\big(m N\big)$,
which scales linearly with the number of grid points and remains negligible compared to the total ROM computational cost.

For ROMs, unlike spectral methods, we do not expect the Gibbs condition to hold between data-driven basis and Gegenbauer polynomials. 
Therefore, we do not expect to recover spectral accuracy. Another central question is how to choose the parameters $\lambda, m.$ For spectral methods, the theory indicates choosing $\lambda, m$ to be proportionate to the number $N.$ Later studies performed in the literature showed various strategies of parameter choices, including using a fixed constant, choice according to the length of the subinterval and smoothness of the function in each subintervals \cite{shizgal2003towards,gelb2004parameter,jackiewicz2004determination}.
 Through extensive tests, we show that $\lambda$ and $m$ play a crucial role in the quality of reconstruction for ROM solutions. Since $m$ represents the degree of the reconstruction parameter, it should not be chosen too small or too large (to avoid under- and over-fitting). The choice of $\lambda$ is also important as demonstrated in detailed study in Section \ref{Numerical Results}. We defer the discussion and detailed comparisons to Section \ref{Numerical Results}.

\subsubsection{2D case}
 
Gegenbauer reconstructions can be extended to multivariate functions. For a regular rectangular domain, high-dimensional Gegenbauer reconstruction can be regarded as the tensor product of multiple one-dimensional Gegenbauer reconstructions. Thereby, the simplest method will be using a line-by-line approach. For each fixed $x$ (or $y$) value, one can perform a reconstruction along that dimension. On the other hand, it is also possible to do a direct reconstruction on a 2D rectangular domain by a tensor product of Gegenbauer polynomials in 1D. For non-rectangular domains, \cite{kawai2022gegenbauer} proposed utilizing the Rosenblatt transformation to map local variables from the irregular subdomains onto a rectangular domain, then a truly 2D reconstruction can be obtained. 

In this paper, we take the simple line-by-line approach. We perform two reconstruction along fixed $x$ and fixed $y$ direction, respectively. The main challenge arises when the detected subdomain is of a small length, which is known to cause large numerical inaccuracy. In those cases, we switch between the fixed $x$ and the fixed $y$ reconstruction to avoid the appearance of large numerical error.
 
Suppose that we have obtained the ROM solution at the final time as $u_{ROM}(x,y),$ which is represented by a solution matrix $X$.
Our approach involves performing Gegenbauer polynomial reconstruction by fixing either $x$ or $y$, followed by a post-processing step. For instance, by fixing the $x-$coordinate, the ROM solution   is reduced to a sequence of 1D functions in terms of $y$. For this function, discontinuities can be identified to partition the domain into multiple analytic subregions. The Gegenbauer polynomial reconstruction algorithm for the 1D case is then applied to each subregion. The procedures for fixed $x$ and fixed $y$ are summarized as Algorithm~\ref{alg_1} and Algorithm~\ref{alg_2}, respectively, which are straightforward extensions of the 1D methods. The output of the algorithms are the reconstructed values together with a metric denoting the size of the smallest analytic sub-interval. 

In Algorithm~\ref{alg_3}, we combine the two reconstructions. The threshold value, $\textit{threshold}$, is a predefined value to account for very small intervals. The idea is that we initialize with one of the reconstruction method, and if the interval is too small, we replace that line by the results from the alternative reconstruction from the other direction. The selection of this threshold is empirical. We note that it is also possible to initialize with a fixed $y$ reconstruction and correct with a fixed $x$ reconstruction. As we demonstrated in the numerical experiments, Algorithm~\ref{alg_3} can significantly improve upon the reconstruction results from Algorithm~\ref{alg_1} and Algorithm~\ref{alg_2}. However, there still exist situations that the analytic subintervals in both $x-$ and  $y-$ subintervals are small. Those will pose numerical challenges for the reconstruction.

The overall computational cost of the 2D line-by-line Gegenbauer reconstruction is inexpensive relative to the offline and online stages of the ROM procedure.
Assume that each row or column contains $N$ grid points. 
The first step consists of performing the Gegenbauer reconstruction along each row and each column, which is analogous to the 1D case. 
For each row or column, the reconstruction cost is $\mathcal{O}(m N)$, where $m$ is the Gegenbauer truncation degree. 
Since there are $N$ rows and $N$ columns, the total cost for these line-by-line reconstructions is $\mathcal{O} (m N^2)$. 
In addition, the threshold-based combination step, which selects between the fixed-$x$ and fixed-$y$ reconstructions based on edge distances, requires traversing all points and performing comparisons, leading to a cost of $\mathcal{O}(N^2)$. 
Consequently, the overall computational complexity of the 2D reconstruction is 
$T_{2\mathrm{D}} = \mathcal{O}(m N^2)$, which scales quadratically with the number of grid points along each dimension and remains negligible compared to the cost of generating the full 2D ROM solution.

\begin{algorithm}
\caption{Gegenbauer reconstruction for fixed $x$}
\label{alg_1}
\begin{algorithmic}[1]
\STATE \textbf{Input:} Solution matrix $X \in \mathbb{R}^{N_x \times N_y}$, Gegenbauer parameter $\lambda$ and degree $m$ for each analytical region. 

\STATE \textbf{Output:} $\text{distance}\_x \in \mathbb{R}^{N_x}$ and reconstructed matrix $\text{Result}\_x \in \mathbb{R}^{N_x \times N_y}$.

\FOR{$n = 1:N_x$}
    \STATE Perform discontinuity detection on \( X[n, :] \), and separate the domain into analytical regions based on the detected discontinuity points.
    
    \IF{no discontinuity is detected}
        \STATE Set $\text{distance}\_x[n] = N_y$;
    \ELSE
        \STATE Set $ \text{distance}\_x[n] = $ the number of sample points contained in the smallest analytical region.
    \ENDIF
    
    \STATE For each analytical region, perform Gegenbauer polynomial reconstruction with parameters  $\lambda$  and degree $m$, and store the result in $\text{Result}\_x[n, :]$.

\ENDFOR
\STATE \textbf{return} $\text{distance}\_x$, $\text{Result}\_x$.
\end{algorithmic}
\end{algorithm}

\begin{algorithm}
\caption{Gegenbauer reconstruction for fixed $y$}
\label{alg_2}
\begin{algorithmic}[1]
\STATE \textbf{Input:} Solution matrix $X \in \mathbb{R}^{N_x \times N_y}$, Gegenbauer parameter $\lambda$ and degree $m$ for each analytical region.

\STATE \textbf{Output:} $\text{distance}\_y \in \mathbb{R}^{N_y}$ and reconstructed matrix $\text{Result}\_y \in \mathbb{R}^{N_x \times N_y}$.

\FOR{$n = 1:N_y$}
    \STATE Perform discontinuity detection on $X[:, n]$, and separate the domain into analytical regions based on the detected discontinuity points.

    \IF{no discontinuity is detected}
        \STATE Set $\text{distance}\_y[n] = N_x$;
    \ELSE
        \STATE Set \( \text{distance}\_y[n] = \) the number of sample points contained in the smallest analytical region.
    \ENDIF
    
    \STATE For each analytical region, perform Gegenbauer polynomial reconstruction with parameters \( \lambda \) and degree \( m \), and store the result in $\text{Result}\_y[:, n]$.
\ENDFOR
\STATE \textbf{return} $\text{distance}\_y$, $\text{Result}\_y$.
\end{algorithmic}
\end{algorithm}

\begin{algorithm}
\caption{Reconstruction combining the fixed $x$ or $y$ reconstruction based on thresholded edge distances }
\label{alg_3}
\begin{algorithmic}[1]
\STATE \textbf{Input:} $\text{Result}\_x$, $\text{Result}\_y$, $\text{distance}\_x$, $\text{distance}\_y$, threshold values $\textit{threshold}$.
\STATE \textbf{Output:} Final reconstruction matrix 
 $\text{Result} \in \mathbb{R}^{N_x \times N_y}$.
 
\STATE $\text{Result} \gets \text{Result}\_x$
    \FOR{$n = 1:N_y$}
        \IF{$\text{distance}\_y[n] > \textit{threshold} \quad \textbf{and} \quad\text{distance}\_y[n] < N_x$}
            \STATE Set $\text{Result}[:, n] = \text{Result}\_y[:, n]$;
        \ENDIF
    \ENDFOR

\STATE \textbf{return} $\text{Result}$.
\end{algorithmic}
\end{algorithm}

\begin{remark}
Although our primary focus is on mitigating spurious oscillations in reduced-order models via Gegenbauer reconstruction, the proposed procedure is not intrinsically restricted to ROM solutions. 
From an algorithmic perspective, the method operates purely as a post-processing step applied to discrete solution vectors and aims at suppressing oscillatory behavior while enhancing resolution near discontinuities. 
Therefore, it is solver-independent and can be readily employed for outputs generated by other numerical or data-driven approaches. 
In particular, in the subsequent inviscid Burgers' example of subsection~\ref{eg2_inviscid_burgers_equation}, we demonstrate that the Gegenbauer reconstruction is equally effective for classical discretization methods and machine learning solvers, where it alleviates either oscillatory artifacts or excessive smoothing near shocks. 
This highlights the practical versatility and broad applicability of the approach.
\end{remark}

\section{Numerical Results}
\label{Numerical Results}
In this section, we perform benchmark numerical tests for 1D and 2D linear and nonlinear problems. For all examples, data snapshots are obtained using exact solutions.
And the FOM is defined by fifth order upwind finite difference scheme with the third-order Runge-Kutta method for time stepping. 
We remark that if we use a numerical method with damping to obtain the solution snapshots, the resulting basis functions are qualitatively similar to those obtained by exact solutions with the pre-processing filter. 
In the following, we use shorthand expressions such as \textit{G-ROM + Post}, where \textit{G-ROM} represents the ROM used, and 
\textit{Post} refers to post-processing by Gegenbauer polynomial reconstruction. 

We measure the numerical error using the relative error and the maximum error. The definitions for these metrics are provided below:
\begin{eqnarray*}
   && \text{Relative\ error} = 
    \frac{(\sum_{i} (u^{\text{reconstruction}}_{i} -  u^{\text{true}}_{i})^2 )^{1/2} }{(\sum_{i} u^{\text{true}}_{i} )^{1/2}}.\\
    &&    \text{Maximum\ error} = \max_{i} \vert u^{\text{reconstruction}}_{i} -  u^{\text{true}}_{i}\vert.
\end{eqnarray*}
Here, $ u^{\text{reconstruction}}_{i} $ represents the $i$-th element of the reconstruction result, and $ u^{\text{true}}_{i} $ denotes the $i$-th element of the corresponding true solution with the same mask applied. 
To take into account slight inaccuracies in discontinuity detection,  five grid points near the identified discontinuities in each spatial direction are excluded from the error computation.

The experiments are carried out on a MacOS 14.4.1 system (Build Version: 23E224) with Python 3.9.7. The environment is configured with TensorFlow 2.7.0, PyTorch 2.0.1, and other dependencies, including Pandas 2.0.3, NumPy 1.23.5, and SciPy 1.10.1. The system features an Apple M1 chip with 16 GB RAM, and no GPU acceleration was utilized during the experiments.

\subsection{Example 1: linear equation}
 
In this example, we consider a linear transport equation with discontinuous initial conditions. 
\begin{equation}
\label{eg1_hyperbolic_equation}
\begin{cases}
    \begin{aligned}
        &\frac{\partial u}{\partial t} = -2\pi \frac{\partial u}{\partial x}, \quad x\in[0,2\pi],\\
        &u(0,t) = u(2\pi ,t).
    \end{aligned}
\end{cases}
\end{equation}
The problem possesses periodic boundary conditions. Consequently, extending the initial condition periodically over the entire domain with period $2\pi$,
the exact solution of \eqref{eg1_hyperbolic_equation} can be expressed as $ u(x,t) = u(x - 2\pi t, 0) $.
We focus on postprocessing G-ROM solution in this example. When constructing the snapshot matrix with exact solution, the spatial domain is uniformly discretized as $0=x_1<\cdots<x_{256}=2\pi$, and the temporal domain is uniformly discretized as 
$0=t_0<t_1<\cdots<t_{800}=0.8$. 

\begin{figure}[!htbp]
    \centering
    \includegraphics[scale=0.36]{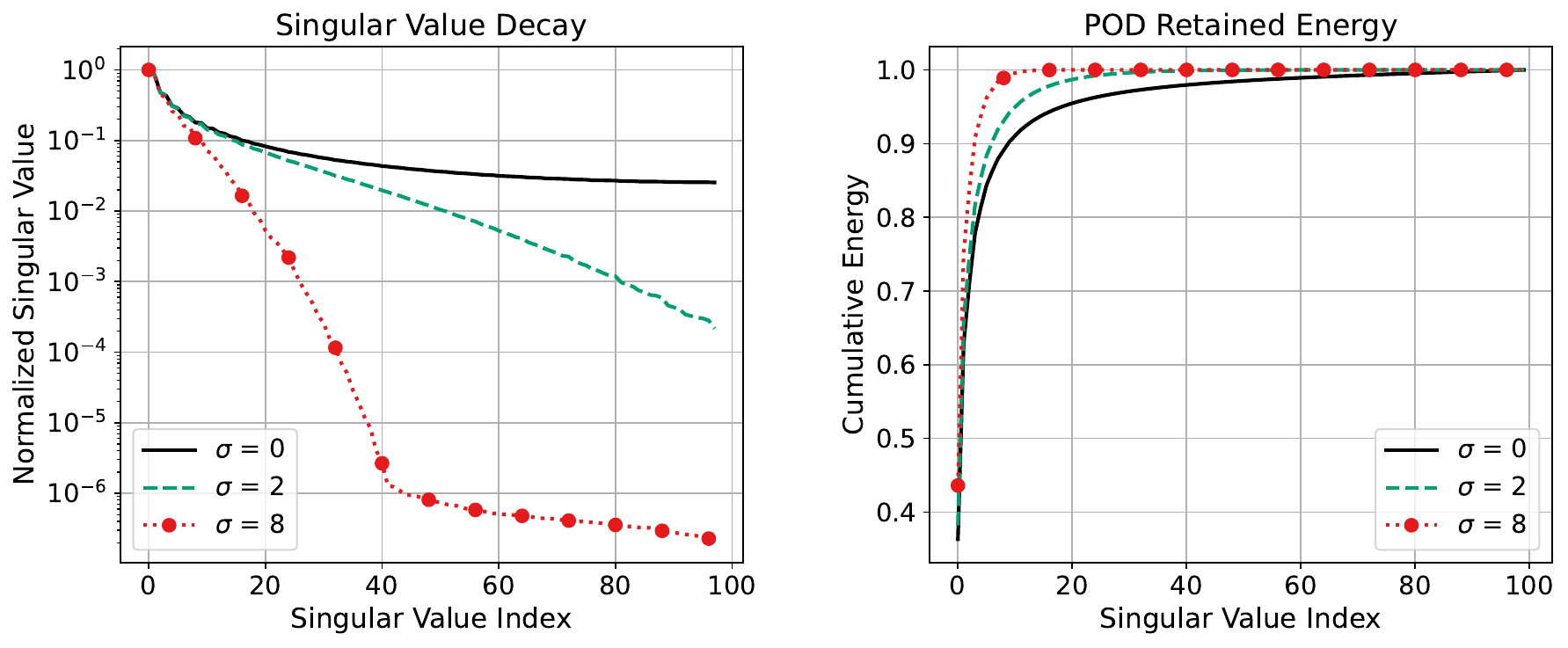}
    \caption{Example 1. Normalized singular values (left) and POD retained energy (right) for different filter widths $\sigma$ in the pre-processing.}     \label{eg1_svd_energy}
\end{figure}

First, we consider the following piecewise linear initial condition.
\begin{equation}  
    \label{eg1_linear_ic}
    u(x,0) =  
    \begin{cases}  
        x, & x \in [0, \pi], \\  
        x-2\pi, & x \in [\pi, 2\pi].  
    \end{cases}  
\end{equation}  
We pre-process the snapshots by Gaussian filters \cite{farcas2022filtering} (with parameter $\sigma$, and $\sigma=0$ means no pre-processing is used). Similar to what is observed in \cite{farcas2022filtering}, the effect of filtering on the decay of the singular values of the filtered snapshot matrix is obvious and is shown in Figure~\ref{eg1_svd_energy}.

\begin{figure}[!htbp]
    \centering
    \subfloat[$\sigma=0$]{\includegraphics[width=0.33\textwidth]{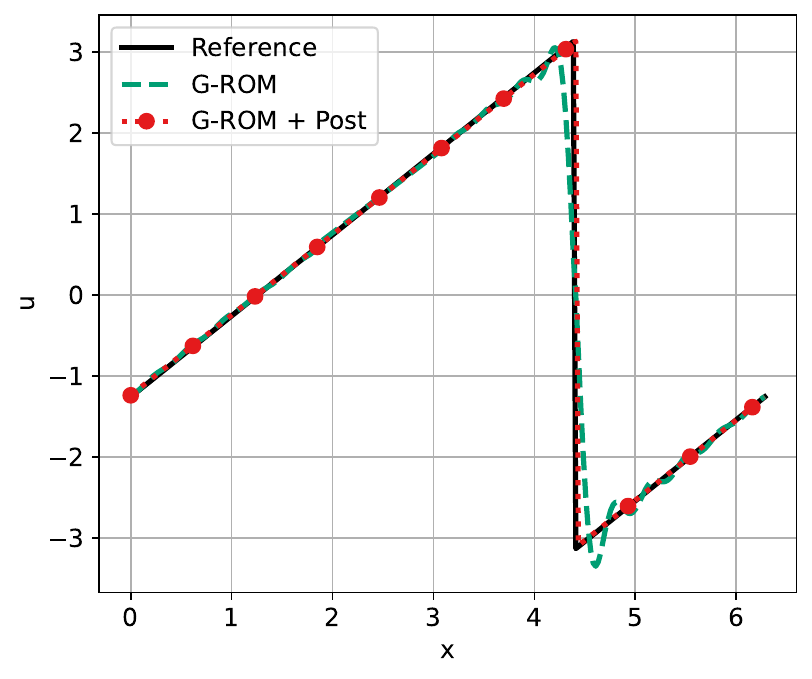}}
    \subfloat[$\sigma=2$]{\includegraphics[width=0.33\textwidth]{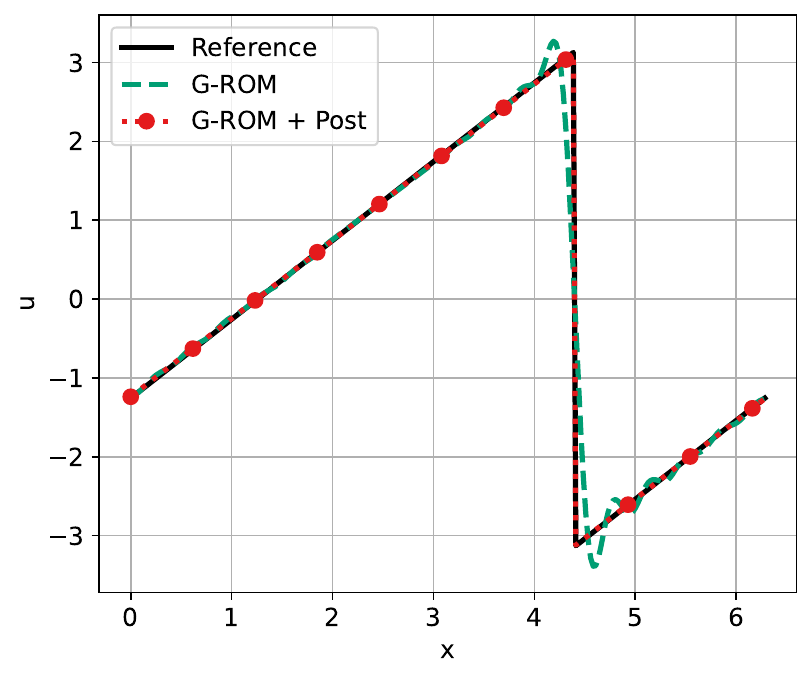}}
    \subfloat[$\sigma=8$]{\includegraphics[width=0.33\textwidth]{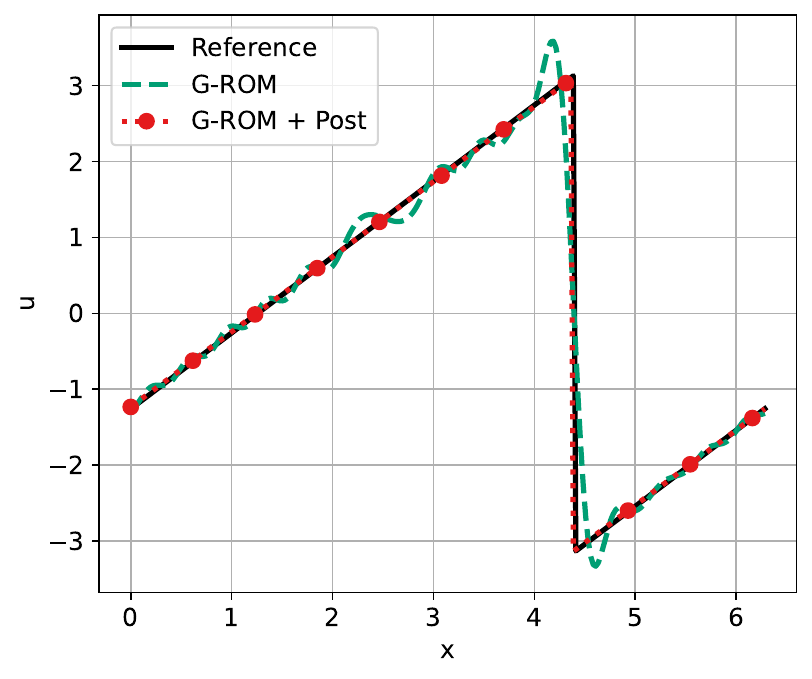}}
    \caption{Example 1. The G-ROM solution for initial condition \eqref{eg1_linear_ic} (with or without post-processing) at $t=0.2$ with different filter widths in the pre-processing.} 
    \label{eg1_t1}
\end{figure}

Figure~\ref{eg1_t1} shows the G-ROM solution with or without Gegenbauer post-processing at $t=0.2$ with reduced dimension $r=30.$ In Gegenbauer polynomial reconstruction, for now we use reconstruction parameters $\lambda = 2$ and $m = 1$ (motivated by the fact that the solution is linear).    As we can see, for all cases (different widths of the pre-processing filter), post-processing successfully removes spurious oscillations. Furthermore, Table \ref{tab:error1} collects the comparison of numerical errors with or without post-processing. It is clear that for all cases, the numerical error decreases by one to two orders of magnitude.

\begin{table}[!htbp]
    \centering
    \caption{\label{tab:error1}Example 1. Errors of the G-ROM solution for initial condition \eqref{eg1_linear_ic} at $ t = 0.2 $ with different filter widths in the pre-processing.}
    \begin{tabular}{c|l|c|c|c}
        \hline
        & & $\sigma = 0$ & $\sigma = 2$ & $\sigma = 8$ \\\hline
        \multirow{2}{*}{Relative error} 
        & G-ROM& 1.0078e-01 &  8.2851e-02 & 1.0172e-01 \\
        & G-ROM+Post & 5.8525e-03 & 4.8993e-03 & 7.3401e-03 \\\hline
        \multirow{2}{*}{Maximum error} 
        & G-ROM & 1.9711e+00 & 1.4634e+00& 1.3206e+00 \\
        & G-ROM+Post &2.0660e-02 & 1.7565e-02 & 2.4633e-02 \\\hline
    \end{tabular}   
\end{table}

We further compare Gegenbauer post-processing with the standard total variation (TV) regularization \cite{rudin1992nonlinear}.
We take the G-ROM solution at $ t = 0.2 $ with a filter width of $ \sigma = 2 $ for this comparison. We take the TV regularization parameter to be $ \lambda_{TV} $ is set to $ 0.2 $, $ 1.0 $, and $ 4.0. $ The results are reported in Table~\ref{tab_t1_tv}. It is clear that the TV regularization, while effective in controlling oscillations, does not work in reducing the numerical errors. 

\begin{table}[!htbp]
    \centering
    \caption{\label{tab_t1_tv}Example 1. Error of the G-ROM solutions  for initial condition \eqref{eg1_linear_ic} for $t = 0.2$ and $\sigma = 2$ with TV regularization (with regularization parameter $\lambda_{TV}$), in comparison with the Gegenbauer polynomial post-processing. }
    \begin{tabular}{c|c|c|c|c|c}\hline
        & G-ROM & $\lambda_{TV}=0.2$ & $\lambda_{TV} = 1.0$ & $\lambda_{TV} = 4.0$ & G-ROM + Post\\
        \hline
        Relative error& 8.2851e-02 &7.8780e-02&7.6022e-02 & 9.2020e-02 & 4.8993e-03  \\
        \hline
        Maximum error&1.4634e+00 &1.4644e+00 & 1.4806e+00 & 1.4660e+00 & 1.7565e-02\\
        \hline
    \end{tabular}
\end{table}

\begin{figure}[!htbp]
    \centering
    \subfloat[$m=1, \lambda=2$]{\includegraphics[width=0.33\textwidth]{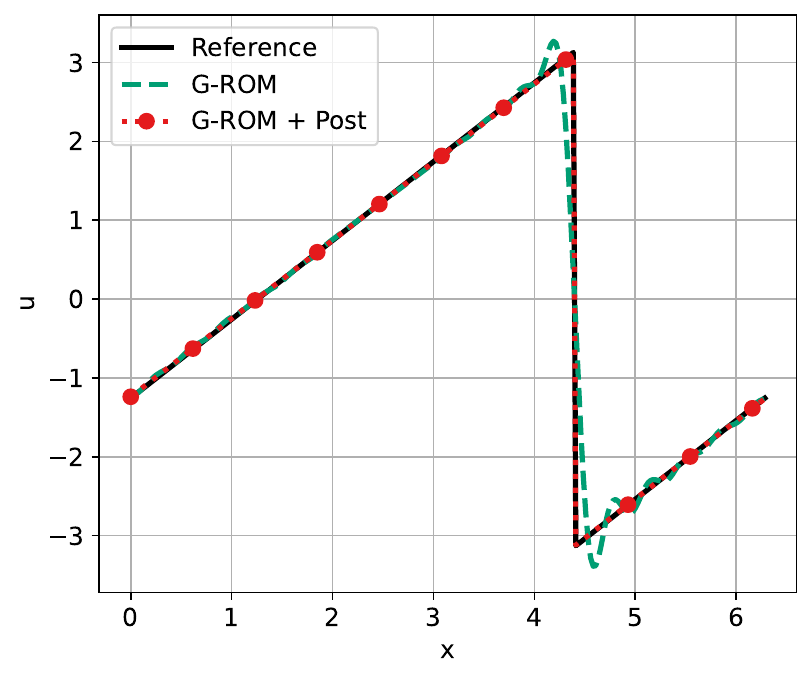}}
    \subfloat[$m=3, \lambda=3$]{\includegraphics[width=0.33\textwidth]{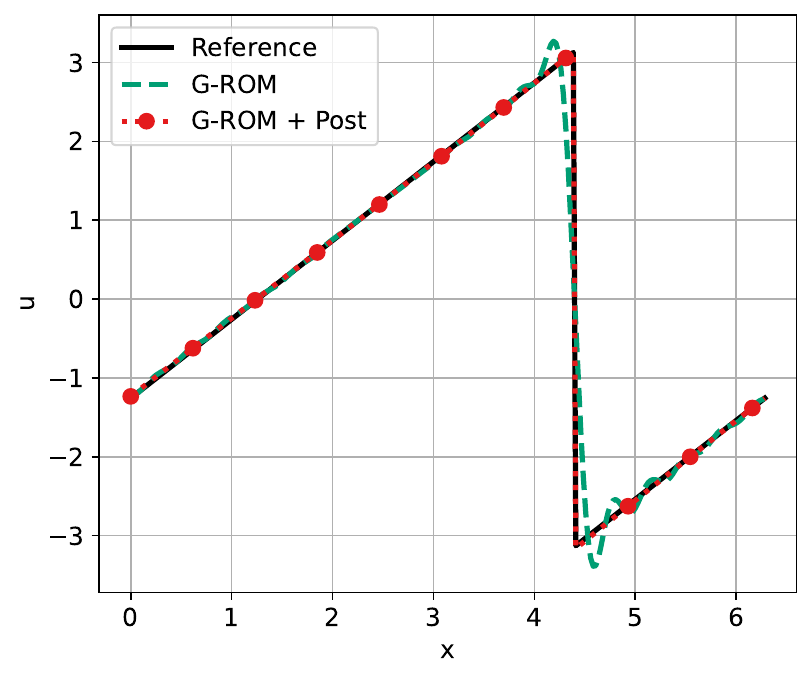}}
    \subfloat[$m=5, \lambda=3$]{\includegraphics[width=0.33\textwidth]{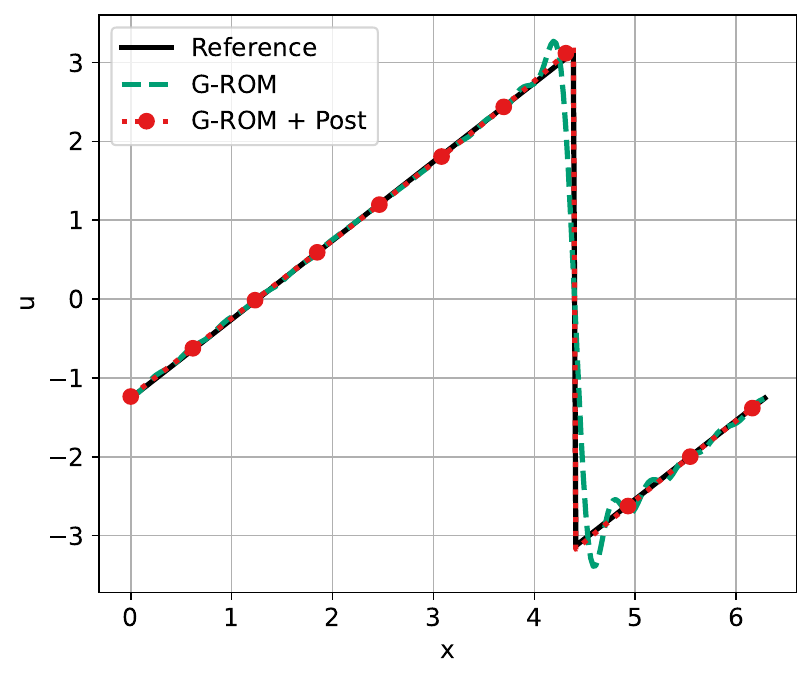}}\\
    \subfloat[$m=5, \lambda=10$]{\includegraphics[width=0.33\textwidth]{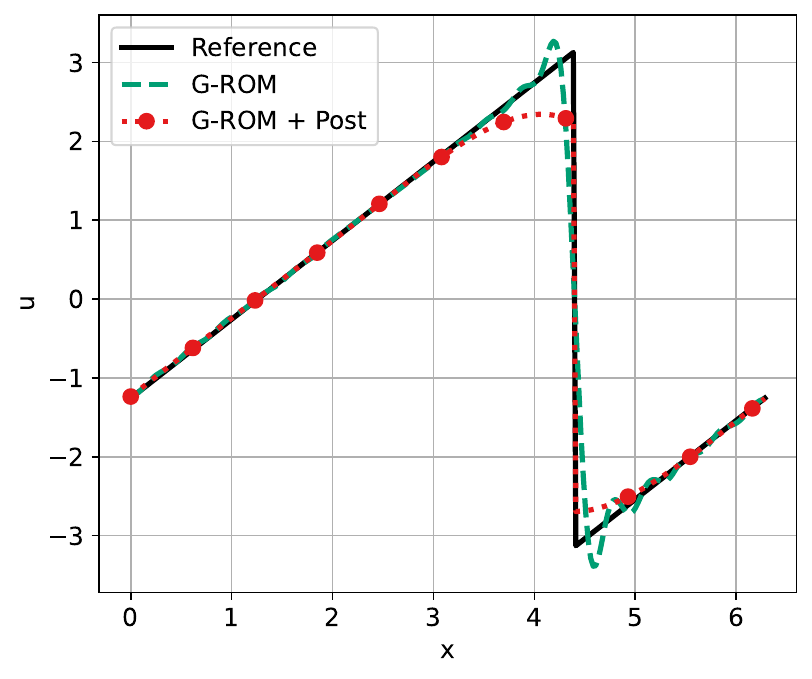}}
    \subfloat[$m=1, \lambda=0.5$]{\includegraphics[width=0.33\textwidth]{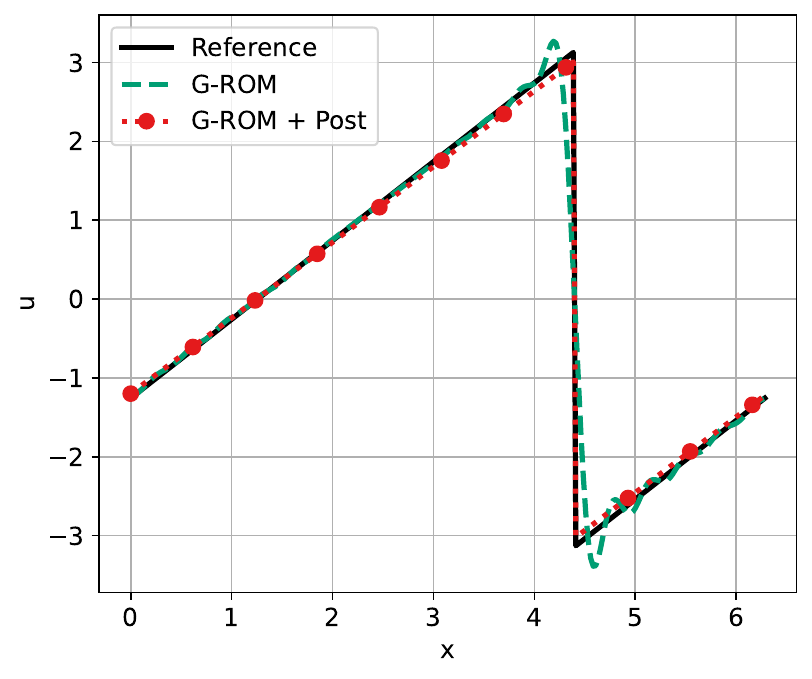}}
    \subfloat[$m=3, \lambda=0.5$]{\includegraphics[width=0.33\textwidth]{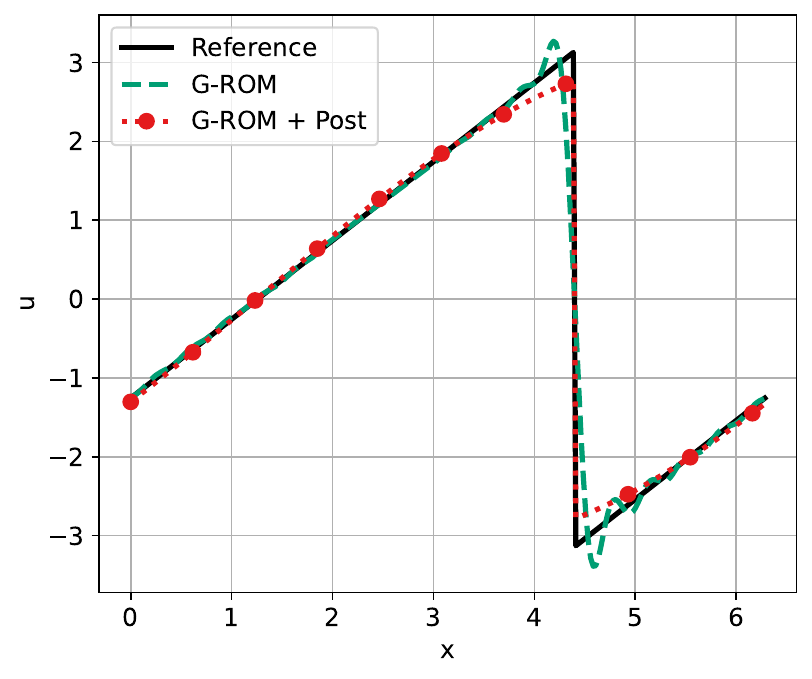}} 
    \caption{Example 1. The G-ROM solution  for initial condition \eqref{eg1_linear_ic} at $t=0.2$   with different reconstruction parameters $\lambda$ and $m$. 
    }
    \label{eg1_spectral_rom_linear_ic}
\end{figure}

\begin{figure}[!htbp]
    \centering
    \includegraphics[width=0.48\linewidth]{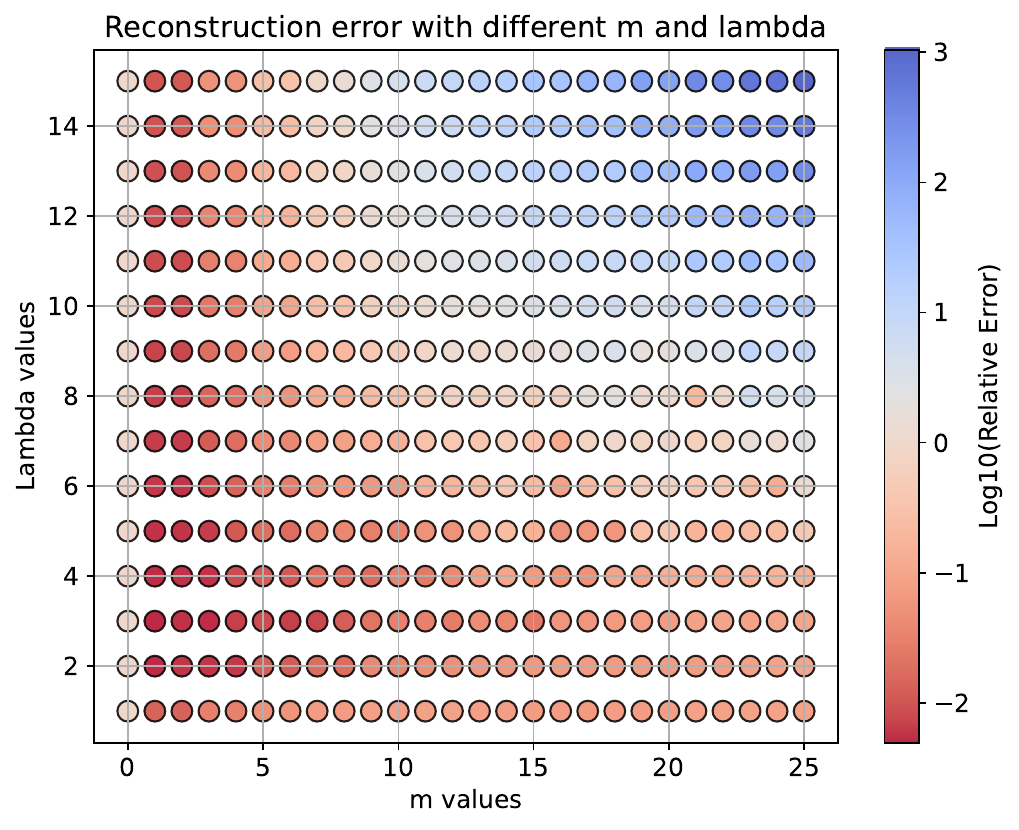}
    \caption{Example 1. Relative error of the G-ROM solution for initial condition \eqref{eg1_linear_ic} at $t=0.2$, with different reconstruction parameters $\lambda$ and $m$, in the case of $r=30$ and $\sigma = 2$.}
    \label{eg1_r_sigma_linear_ic}
\end{figure}

Next, we vary the selection of $\lambda$ and $m$ in the Gegenbauer reconstruction. In Figure \ref{eg1_spectral_rom_linear_ic}, we present the reconstruction results using a variety of $\lambda, m.$ We note that the case of $\lambda=0.5$ corresponds to Legendre polynomial, i.e. the reconstruction is the optimal $L^2$ projection of the G-ROM solution in each subinterval. From the figure, we can clearly see that choosing $\lambda=2$ for $m=1,$ and $\lambda=3$ for $m=3$ significantly improves upon the optimal $L^2$ projection obtained by $\lambda=0.5.$ However, if $\lambda, m$ are chosen to be too large, it negatively affects the quality of the reconstruction. Then we display the relative error of the reconstruction for values of $\lambda$ ranging from 1 to 15 and $m$ ranging from 0 to 25 in Figure \ref{eg1_r_sigma_linear_ic}.   The optimal parameters found through traversal are $m=1$ and $\lambda=2$, while the second-best parameters are $\lambda=3$ and $m=3.$ The overall quality of the reconstruction tends to be better when $m$ and $\lambda$ are chosen not to be exceedingly large. 

\begin{figure}[!htbp]
    \centering
    \subfloat[$m=2, \lambda=4$]{\includegraphics[width=0.33\textwidth]{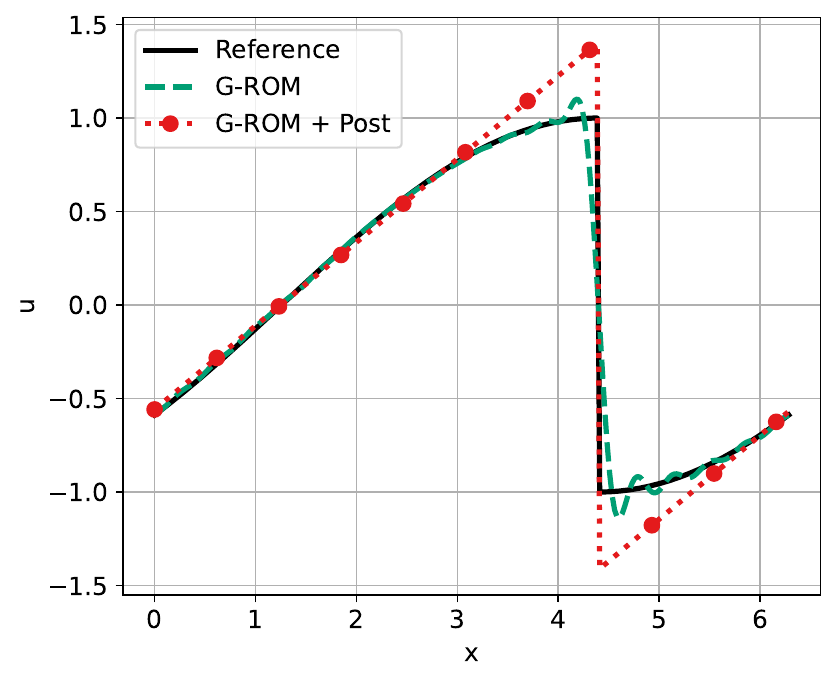}}
    \subfloat[$m=3, \lambda=6$]{\includegraphics[width=0.33\textwidth]{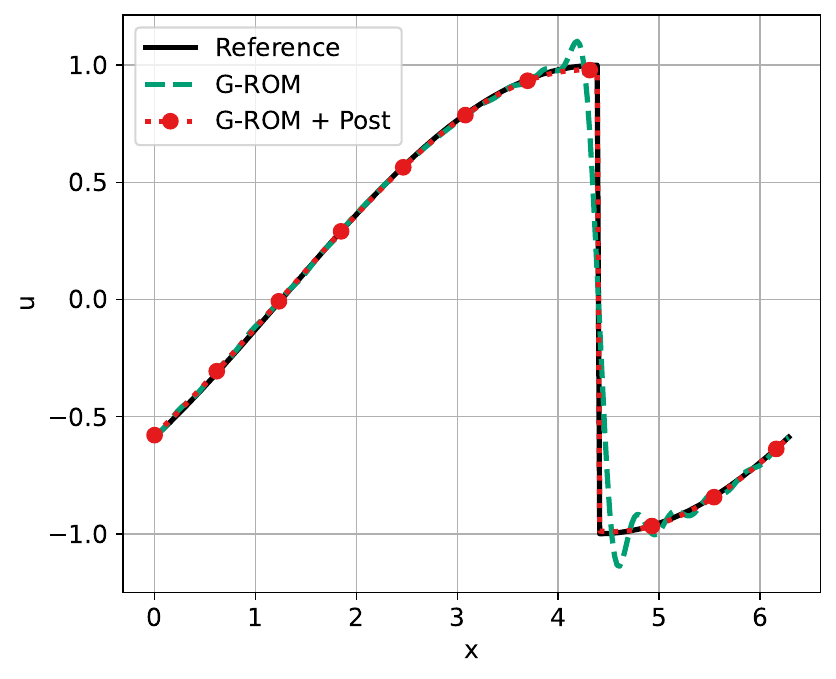}\label{eg1_0.2_new_b}}
    \subfloat[$m=8, \lambda=3$]{\includegraphics[width=0.33\textwidth]{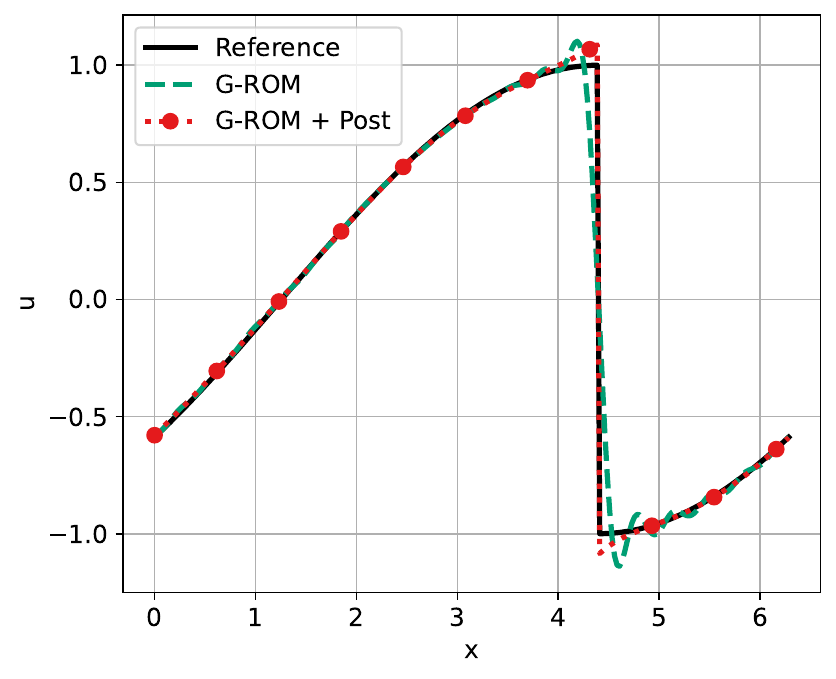}}\\
    \subfloat[$m=8, \lambda=16$]{\includegraphics[width=0.33\textwidth]{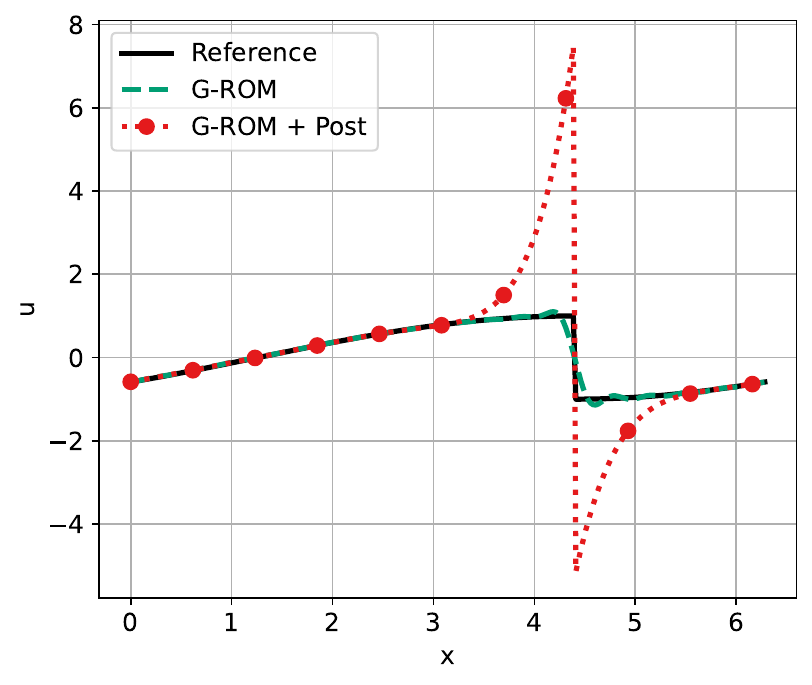}}
    \subfloat[$m=3, \lambda=0.5$]{\includegraphics[width=0.33\textwidth]{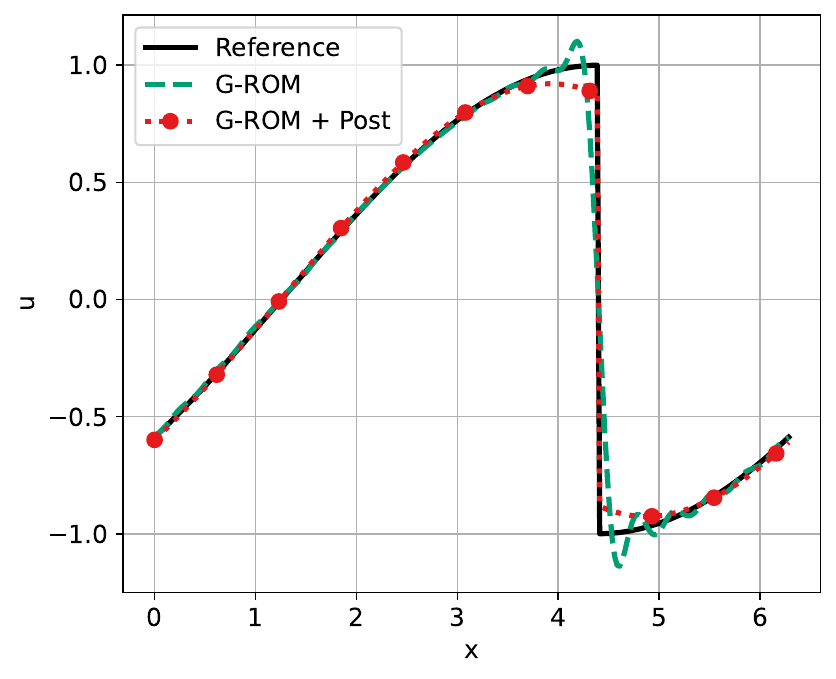}}
    \subfloat[$m=8, \lambda=0.5$]{\includegraphics[width=0.33\textwidth]{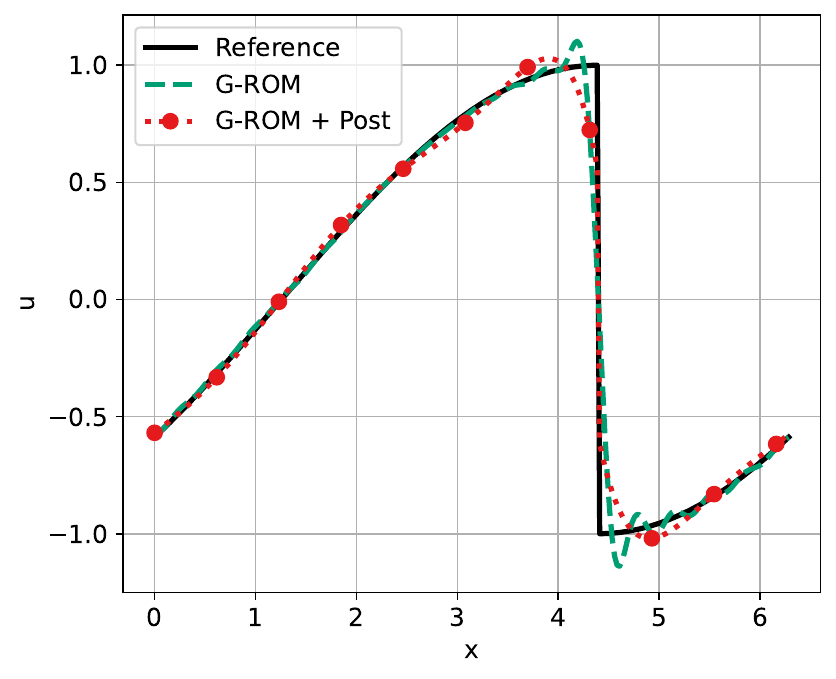}}
    \caption{Example 1. The G-ROM solution for initial condition \eqref{eg1_new_ic} at $t=0.2$ , with different reconstruction parameters $\lambda$ and $m$. 
    }
    \label{eg1_t1_new}
\end{figure}

In the next experiment, we switch the initial function to be
\begin{equation}  
    \label{eg1_new_ic}
    u(x,0) =  
    \begin{cases}  
        \sin(\frac{1}{2}x), & x \in [0, \pi], \\  
        \sin(\frac{1}{2}x+\pi), & x \in [\pi, 2\pi].  
    \end{cases}  
\end{equation}  
Under this initial condition, the solution at any time is discontinuous, however, it is not a polynomial. For $r=30$ and filter width at $\sigma=2,$ we test six different sets of reconstruction parameters, including the optimal ones, for Gegenbauer polynomial reconstruction, and the results are presented in Figure~\ref{eg1_t1_new}. We observe that when $m, \lambda$ are too small, underfitting occurs, whereas for larger $m, \lambda$,  overfitting effect emerges. The choice of $m=3, \lambda=6$ yields good results. On the other hand, the $L^2$ projection by setting $\lambda=0.5$ has significant errors.

\begin{figure}[!htbp]
    \centering
    \subfloat[$r=10$]{\includegraphics[width=0.4\textwidth]{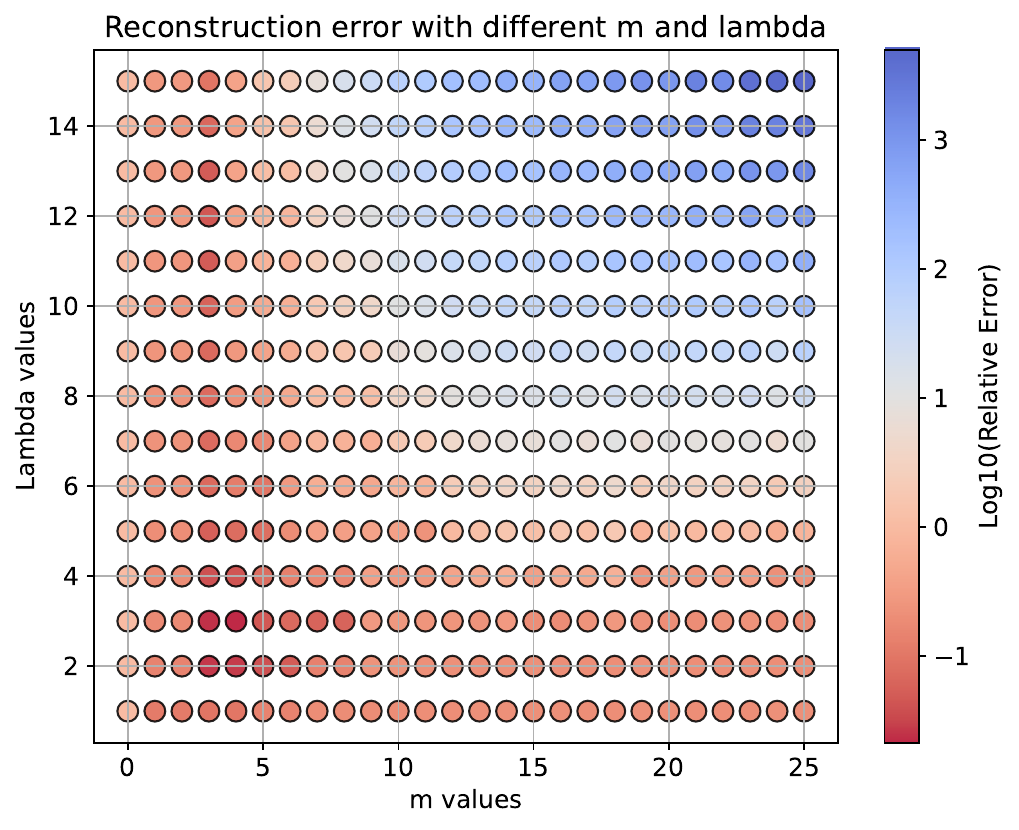}}
    \hspace{0.05\textwidth}
    \subfloat[$r=30$]{\includegraphics[width=0.4\textwidth]{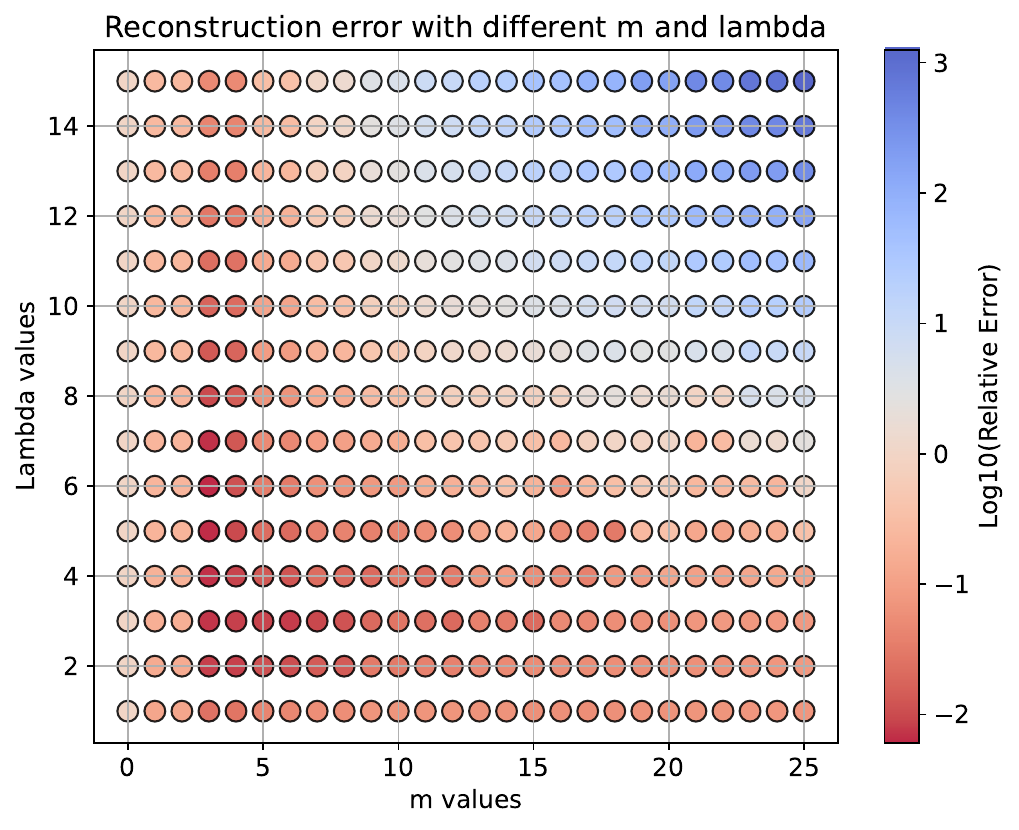}}
    \caption{Example 1. Relative error of the G-ROM solution at $t=0.2$ for initial condition \eqref{eg1_new_ic}, with  $\sigma=2$ and reduced dimensions $r=10,30$ in the pre-processing.}
    \label{eg1_r_sigma}
\end{figure}

To provide a more comprehensive understanding, Figure~\ref{eg1_r_sigma} presents the relative error of the reconstruction of the G-ROM solution for   reduced dimensions $r = 10, 30$ under various selections of $\lambda$ and $m$. It is clear that for bigger value of $r,$ larger values of $\lambda, m$ are admissible. The optimal value for both cases is about $m=3.$  
From this numerical experiment, we can clearly see the optimal choice of $\lambda, m$ depends on the solution. In general, they should be chosen neither too small nor too large. When chosen appropriately, Gegenbauer reconstruction clearly outperforms TV regularization and polynomial regression. 

\begin{figure}[!htbp]
    \centering
    \subfloat[Exact discontinuity location]{\includegraphics[width=0.33\textwidth]{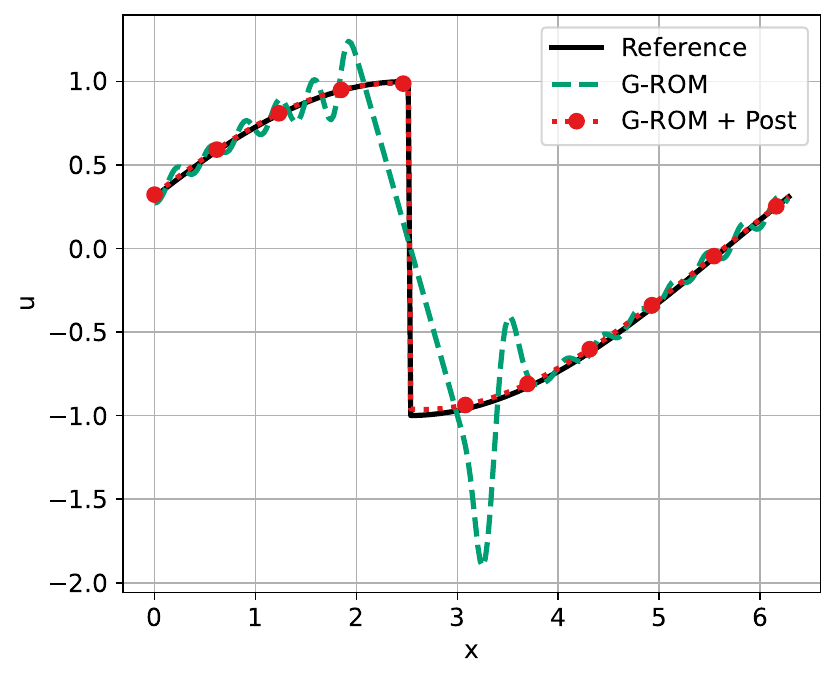} \label{fig:0.9a} }
    \subfloat[Incorrectly detected discontinuity]{\includegraphics[width=0.33\textwidth]{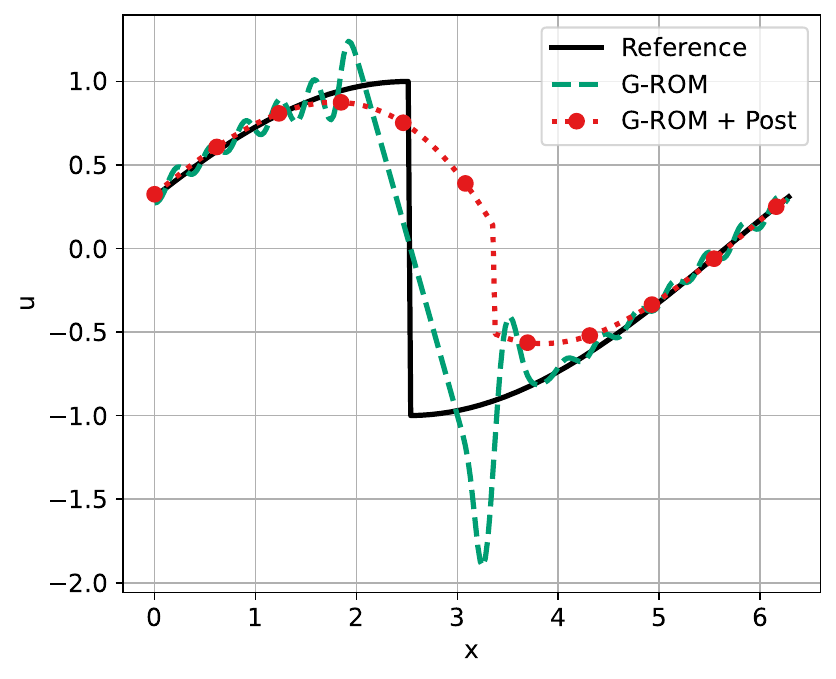} \label{fig:0.9b}} 
    \subfloat[Slightly misidentified (offset by 5 grid points)]{\includegraphics[width=0.33\textwidth]{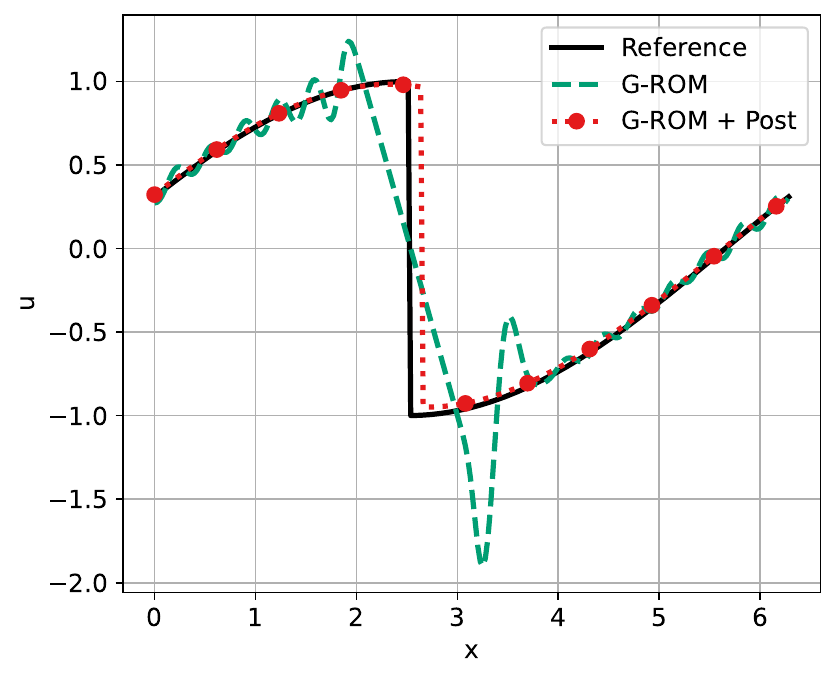} \label{fig:0.9c}} 
    \caption{Example 1: The G-ROM solutions for initial condition \eqref{eg1_new_ic} at $ t = 0.9 $ with reconstruction parameters $m = 3$ and $\lambda = 6$, corresponding to different assumptions on the discontinuity location: (a) Exact location, (b) Incorrect detection, and (c) Slight offset (5 grid points).}
    \label{eg1_2_t2}
\end{figure}

Finally,  we test the application of G-ROM outside the training window, which is known to be extremely challenging.
For $t = 0.9$, beyond the final training time, the G-ROM solution obtained with filter parameter $\sigma = 2$ and reduced dimension $r = 30$ is first computed and analyzed. 
Moreover, we assess the robustness of the proposed framework with respect to errors in discontinuity detection. 
Figure~\ref{eg1_2_t2} presents the reconstructed solutions corresponding to: 
(a) the exact discontinuity location, 
(b) an incorrectly detected discontinuity, 
and (c) a slightly misidentified discontinuity (offset by 5 grid points). 
The reconstruction parameters are fixed to $m = 3$ and $\lambda = 6$. 
We observe that the Gegenbauer post-processing fails to recover the correct solution when the discontinuity is detected incorrectly; in this case, the reconstructed profile deviates significantly from the reference solution. 
For the slightly misidentified case, the reconstruction can effectively capture the smooth components on each side of the jump.
This is because that the Gegenbauer reconstruction is inherently data-dependent. When the boundary region contains a limited number of prominent singularities or sharp variations, these can be interpreted as numerical oscillations. In such cases, effectively removing these oscillations allows the Gegenbauer reconstruction to yield accurate approximations of the smooth segments of the original function. 
However, when the boundary region involves a sufficiently large number of such irregular points or deformation features, their collective influence can no longer be disregarded, causing the reconstructed expression to be smeared and the discontinuity position to be shifted from the true location.
Therefore, a more accurate identification of the discontinuity provides better localized information and consequently yields a more accurate global reconstruction. 
These results confirm that the performance of the proposed framework is strongly dependent on the precision of the discontinuity detection.

\subsection{Example 2: inviscid Burgers' equation}
\label{eg2_inviscid_burgers_equation}

We now consider the inviscid Burgers' equation with moving discontinuity location:
\begin{equation}
\begin{cases}
    \begin{aligned}
        &u_t + uu_x = 0, \qquad x \in [0, 1],\\
        &u(x,0) = \frac{1}{2}+\sin(2\pi x) .
    \end{aligned}
\end{cases}
\end{equation}
The periodic boundary condition is used here. 
In this example, the spatial domain is uniformly discretized as $0=x_1<\cdots<x_{500}=1 $, and the temporal domain is uniformly discretized as $0=t_0<t_1<\cdots<t_{1000}=1$. 
When constructing the snapshot matrix, we use the method of characteristics to compute the solution at each time step.  
We consider ROMs obtained by G-ROM, OpInf, and CAE-LSTM to test at $ t = 0.5 $, when the solution exhibits discontinuities.

\begin{figure}[!htbp]
    \centering
    \subfloat[$\sigma=0$]{\includegraphics[width=0.32\textwidth]{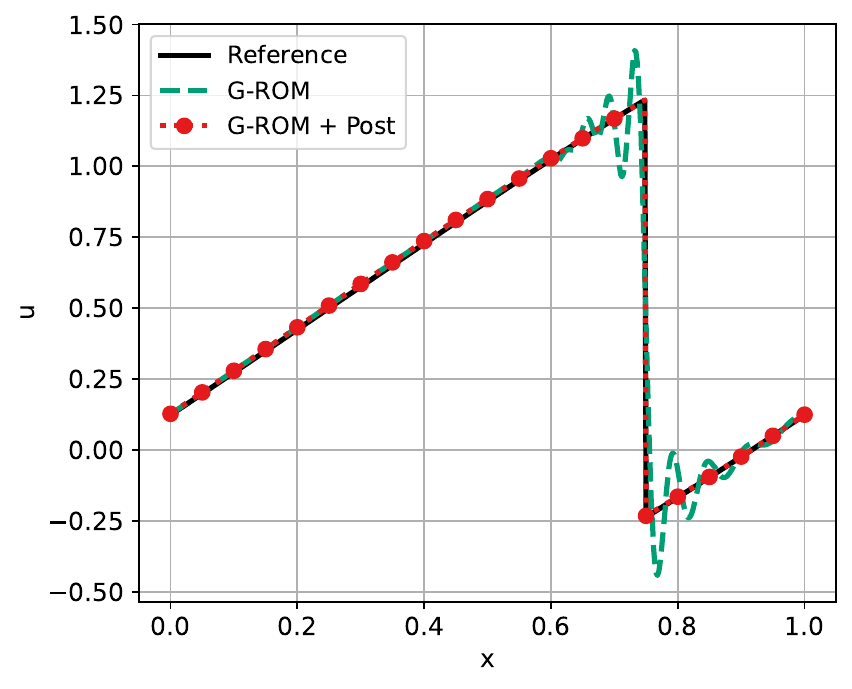}}
    \subfloat[$\sigma=2$]{\includegraphics[width=0.32\textwidth]{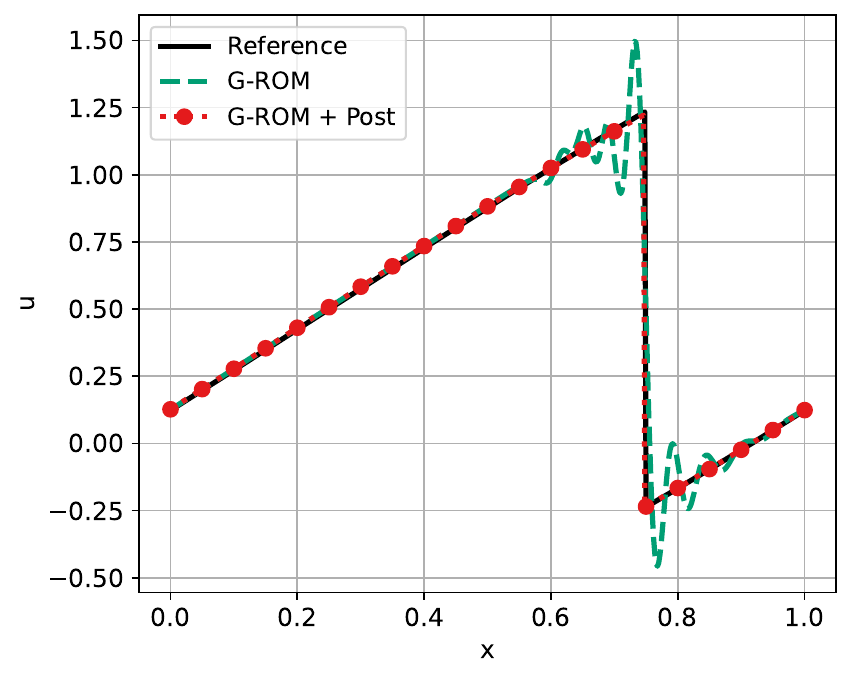}}
    \subfloat[$\sigma=10$]{\includegraphics[width=0.32\textwidth]{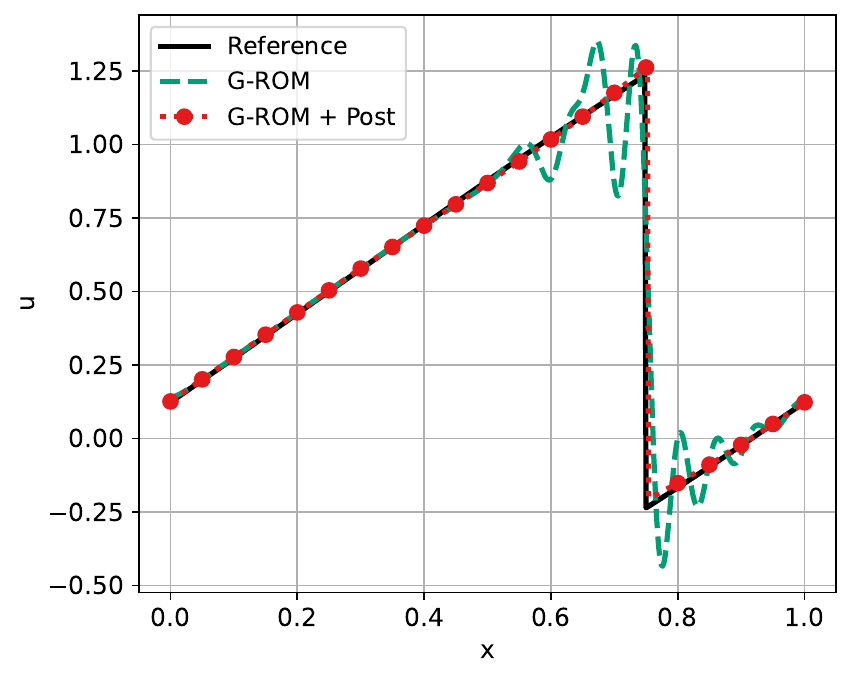}}
    \caption{Example 2. The G-ROM solutions at $ t=0.5 $ under different filter widths, with fixed reconstruction parameters $m = 4$ and $\lambda=3$. 
    }
    \label{eg3_t1}
\end{figure}

\begin{figure}[!htbp]
    \centering
    \subfloat[$\sigma=0$]{\includegraphics[width=0.33\textwidth]{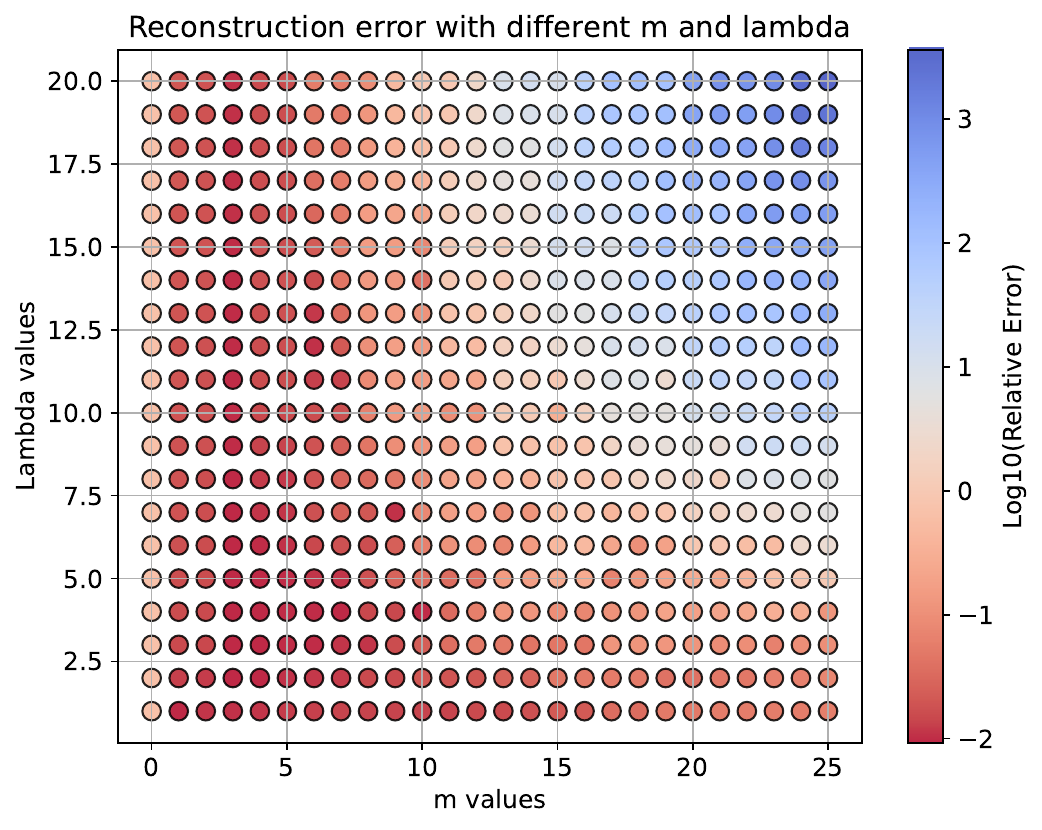}}
    \subfloat[$\sigma=2$]{\includegraphics[width=0.33\textwidth]{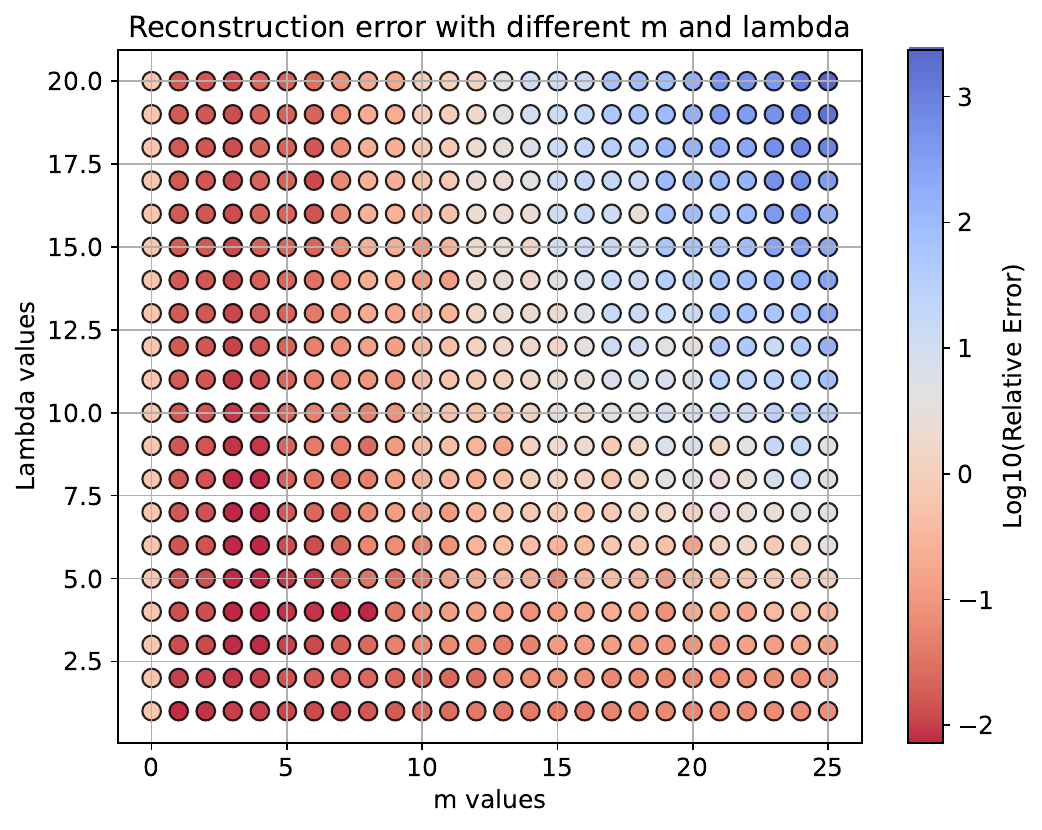}}
    \subfloat[$\sigma=10$]{\includegraphics[width=0.33\textwidth]{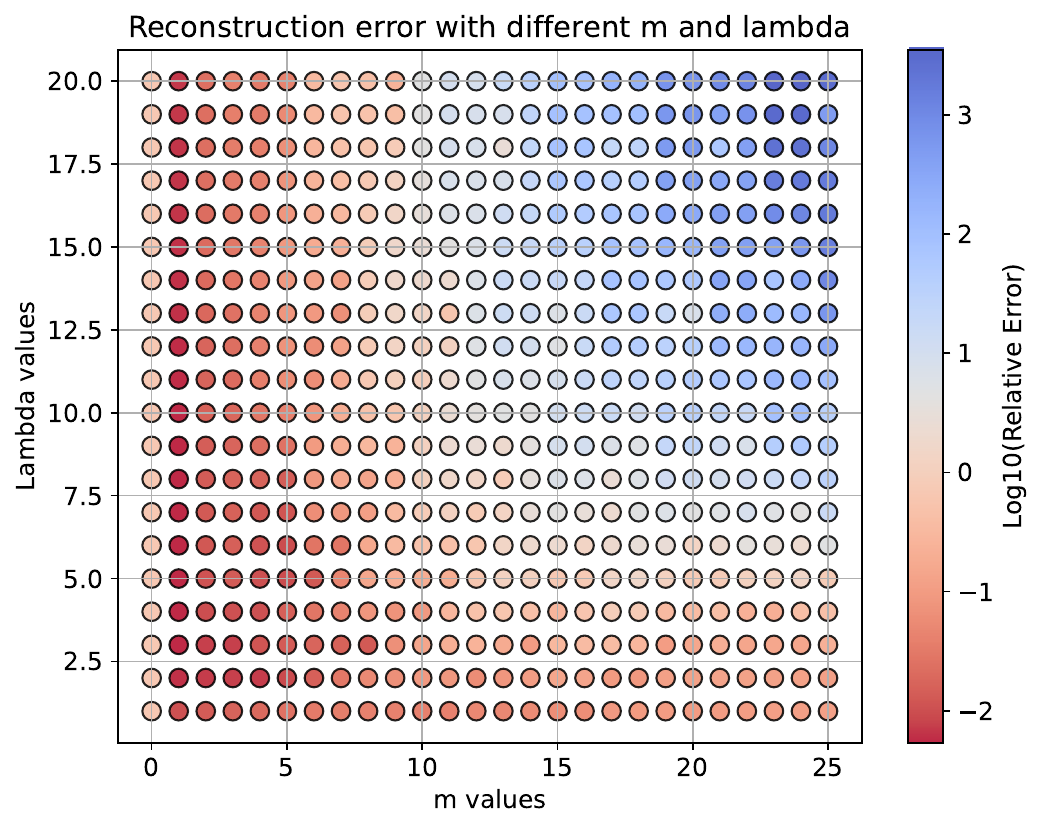}}
    \caption{Example 2. Errors of the Gegenbauer reconstruction of the G-ROM solution at $t=0.5$, with different parameters $\lambda$ and $m$.}
    \label{eg3_pod_sigma}
\end{figure}

We take the reduced dimension $r=25$ first, and
fix the reconstruction parameters to be $m=4$ and $\lambda=3.$
The results before and after post-processing of the G-ROM solution at $t=0.5$ are presented in Figure~\ref{eg3_t1} and Table~\ref{tab_combined_b_c2}. It is clear that in all cases (with $\sigma=0, 2, 10$), the post-processing is effective in removing oscillation and reducing the errors. We note that the results for $\sigma=10$ 
indicate the negative effects of using a large value of $\sigma$ on accuracy. We then vary the parameters $\lambda, m,$ and the results are reported in Figure ~\ref{eg3_pod_sigma}, which gives similar qualitative results as the previous example. 

\begin{table}[!htbp]
    \centering
    \caption{\label{tab_combined_b_c2}Example 2. Errors of the G-ROM solutions at $ t = 0.5 $ under different filter widths, with fixed reconstruction parameters $m=4$ and $\lambda=3$.}
    \begin{tabular}{c|l|c|c|c}
        \hline
        & & $\sigma = 0$ & $\sigma = 2$ & $\sigma = 10$ \\\hline
        \multirow{2}{*}{Relative Error} 
        & G-ROM& 7.5752e-02 & 8.9793e-02& 1.1856e-01 \\
        & G-ROM + Post  & 9.2824e-03 & 7.7443e-03 & 1.1306e-02\\\hline
        \multirow{2}{*}{Maximum Error} 
        & G-ROM & 2.8796e-01 & 3.5721e-01 &3.8437e-01 \\
        & G-ROM + Post  & 9.5881e-03 & 8.9730e-03 & 2.5133e-02 \\\hline
    \end{tabular}
\end{table}

\begin{figure}[!htbp]
    \centering
    \subfloat[$\sigma=0$]{\includegraphics[width=0.32\textwidth]{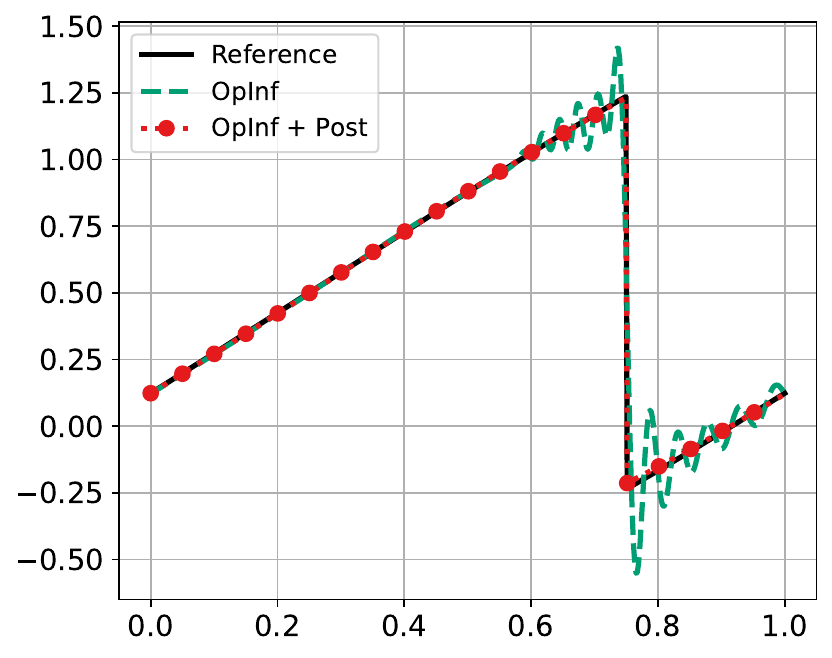}}
    \subfloat[$\sigma=2$]{\includegraphics[width=0.32\textwidth]{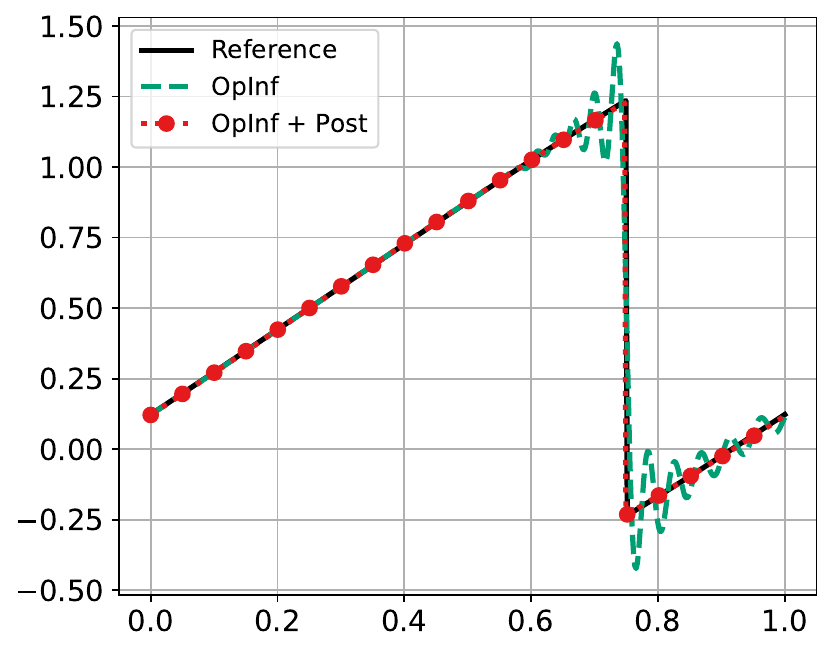}}
    \subfloat[$\sigma=10$]{\includegraphics[width=0.32\textwidth]{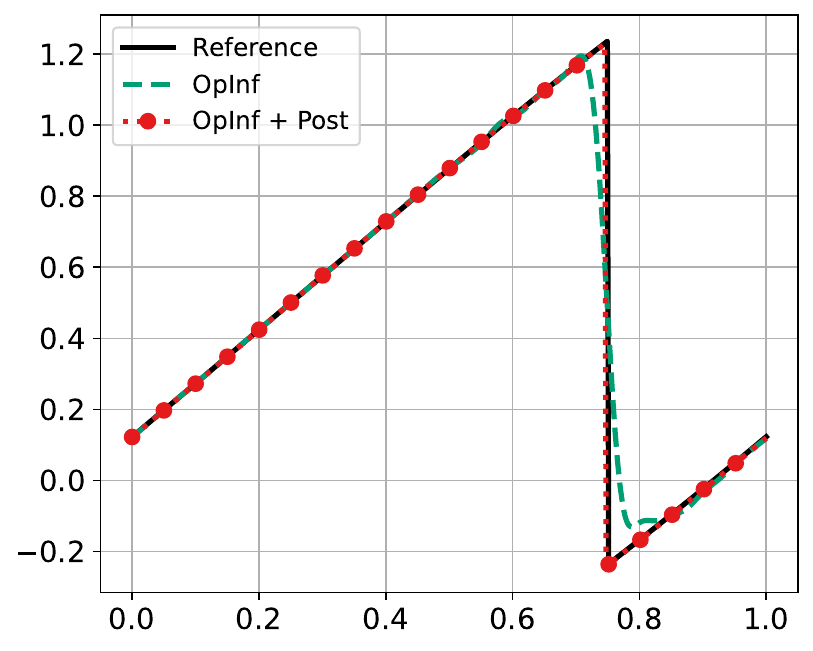}}
    \caption{Example 2. The OpInf solutions at $ t = 0.5 $ under different filter widths, with the same reconstruction parameters $m=3$ and $\lambda=3$.
    }
    \label{eg3_opinf_t1}
\end{figure}

\begin{figure}[!htbp]
    \centering
    \subfloat[$p_{Ridge}=1, \alpha_{Ridge}=0.001$]{\includegraphics[width=0.32\textwidth]{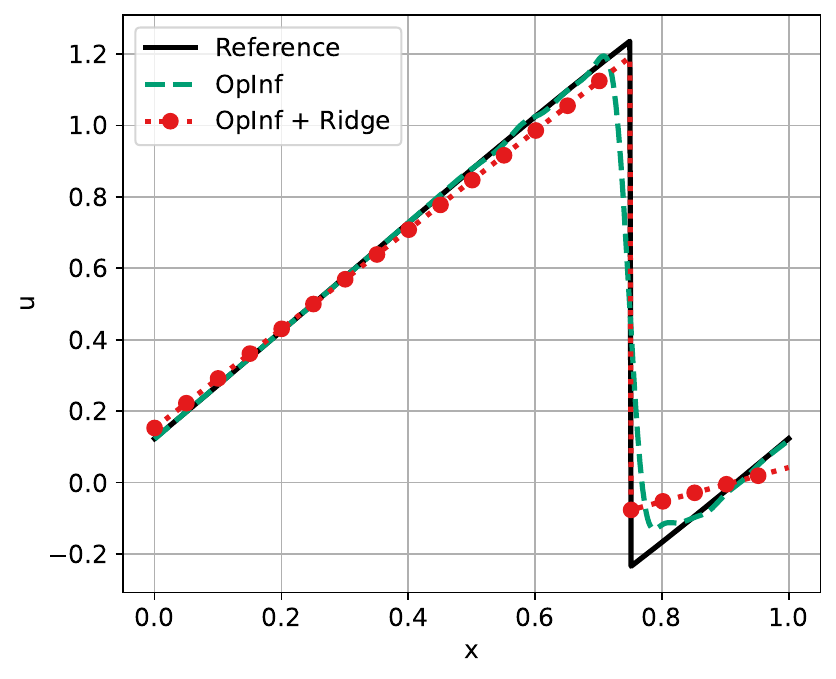}}
    \subfloat[$p_{Ridge}=3, \alpha_{Ridge}=0.001$]{\includegraphics[width=0.32\textwidth]{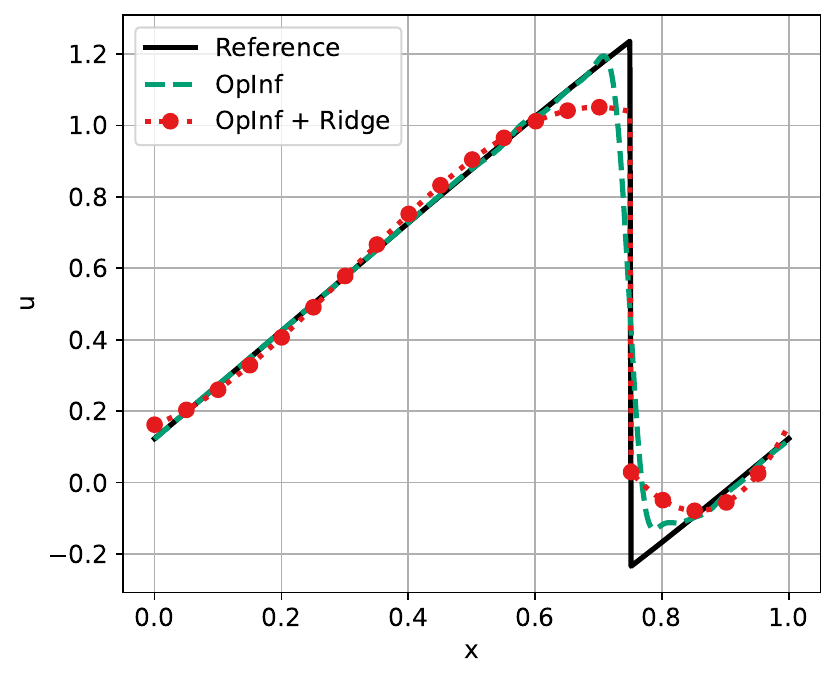}}
    \subfloat[$p_{Ridge}=3, \alpha_{Ridge}=0.1$]{\includegraphics[width=0.32\textwidth]{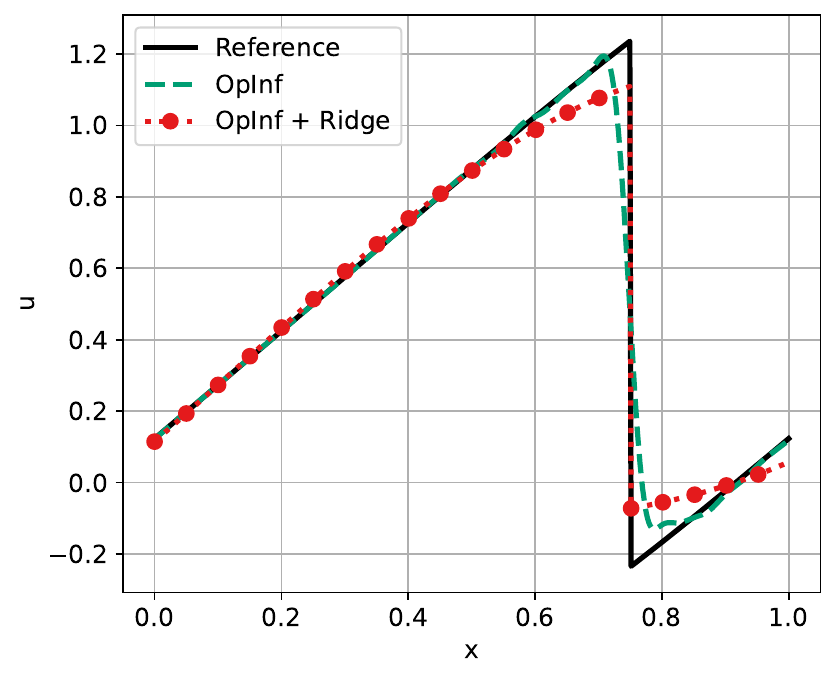}} 
    \caption{Example 2. Results of OpInf solution at $t=0.5$ with $\sigma=10$ corrected by Ridge regression under different polynomial degrees $p_{Ridge}$ and regularization parameters $\alpha_{Ridge}$.
    }
    \label{eg3_opinf_sigma10_Rideg}
\end{figure}

We then perform similar numerical studies for OpInf obtained with reduced dimension $r=25.$
Figure~\ref{eg3_opinf_t1} presents the reconstruction results of the OpInf solution at $t=0.5$ for three different filter widths ($\sigma=0,2,10$) under the same reconstruction parameters $m=3$ and $\lambda=3,$ and Table~\ref{tab_combined_opinf} reports the numerical errors. It is clear that post-processing is effective. In particular, what is interesting is the post-processing result for the case of $\sigma=10.$ As indicated in \cite{farcas2022filtering}, filtering can smear out the shock. 
By Gegenbauer post-processing, the shock transition is sharpened and the numerical error is greatly improved even for the smeared-out solution. In contrast, a standard method such as ridge regression \cite{hoerl1970ridge} is unable to recover such results even when the location of the discontinuity is known, see  Figure~\ref{eg3_opinf_sigma10_Rideg}.

\begin{table}[!htbp]
    \centering
    \caption{\label{tab_combined_opinf}Example 2. Errors of the OpInf solutions at $ t = 0.5 $ under different filter widths, with the same reconstruction parameters $m=3$ and $\lambda=3$.}
    \begin{tabular}{c|l|c|c|c}
        \hline
        & & $\sigma = 0$ & $\sigma = 2$ & $\sigma = 10$\\\hline
        \multirow{2}{*}{Relative Error} 
        & OpInf  & 8.6751e-02 &7.0919e-02 & 1.0733e-01 \\
        & OpInf + Post  & 7.3941e-03 & 3.0021e-03& 2.0312e-03\\\hline
        \multirow{2}{*}{Maximum Error} 
        & OpInf  & 3.3683e-01 & 2.2116e-01 & 4.9526e-01\\
        & OpInf + Post  & 1.9626e-02 & 4.7239e-03 &3.5026e-03\\\hline
    \end{tabular}
\end{table}

Finally, we consider the ROM obtained by CAE-LSTM. In the CAE architecture, the encoder begins by flattening the input data and passing it through a fully connected layer to reduce its dimensionality to 512, followed by reshaping it into a $512 \times 1$ structure. The network then applies a series of \texttt{Conv1D} layers with increasing filter sizes (4, 8, 16, 32, and 64), each followed by \texttt{MaxPooling1D} layers to progressively reduce the sequence length while extracting hierarchical characteristics. After the convolution and pooling operations, the feature maps are flattened, and the data is mapped to a 4-dimensional latent space  through a final dense layer. This compact representation serves as the output of the encoder. The decoder starts with a fully connected layer that transforms the latent space representation into a flat vector of size $16 \times 64$. The subsequent layers progressively upsample the data using multiple \texttt{Conv1DTranspose} layers with \textit{ELU} activations, doubling the sequence length at each step. Specifically, the decoder applies four \texttt{Conv1DTranspose} layers with decreasing filter sizes (64, 32, 16, and 8) while maintaining the same padding. After the upsampling process, the feature maps are flattened, passed through another dense layer to produce a vector of size 500, and finally reshaped into the desired $500 \times 1$ output. In the experiments, the optimizer is Adam with a learning rate of 0.0001, the batch size is set to 20, the maximum number of training epochs is 1000, and the loss function is the mean squared error.

The LSTM network consists of two LSTM layers and a fully connected layer. The first LSTM layer contains 40 units, takes inputs of shape $(\text{time\_window} = 10, \text{latent dimension}=4)$, and outputs a sequence. The second LSTM layer also consists of 40 units, processes the output of the previous layer, and retains only the result of the last time step. The fully connected layer generates the final prediction, with an output dimension which matches the input and no activation function. During training, the model utilizes the Adam optimizer (learning rate 0.001) and the mean squared error as the loss function. An early stopping mechanism monitors the training loss, halting training after 500 epochs of no significant improvement and saving the best weights. The model is trained over 5000 iterations with a batch size of 20. During network training, we observe that adding appropriate filtering in the preprocessing stage reduces both the CAE network and LSTM network training errors. 

\begin{figure}[!htbp]
    \centering
    \subfloat[$\sigma=0$]{\includegraphics[width=0.32\textwidth]{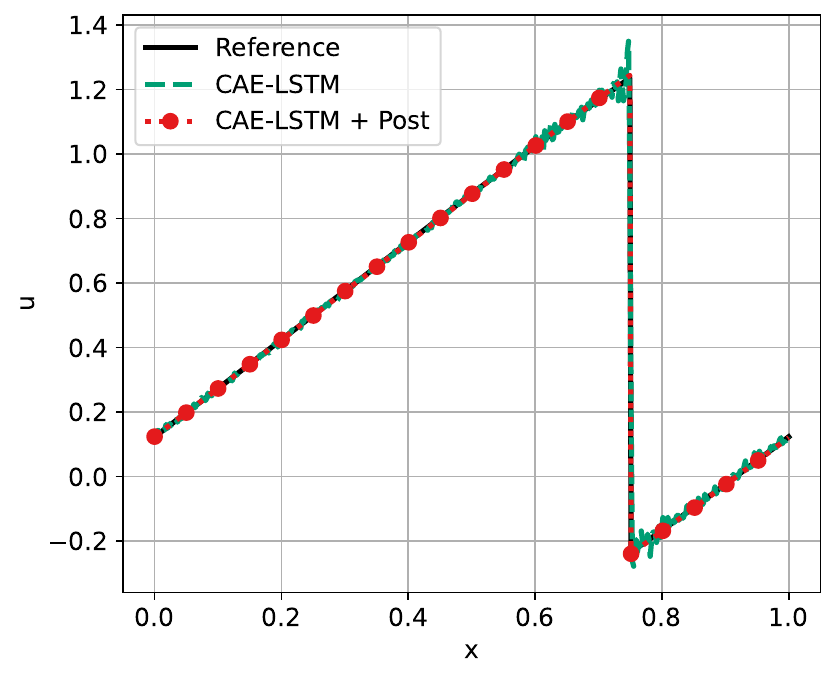}}
    \subfloat[$\sigma=2$]{\includegraphics[width=0.32\textwidth]{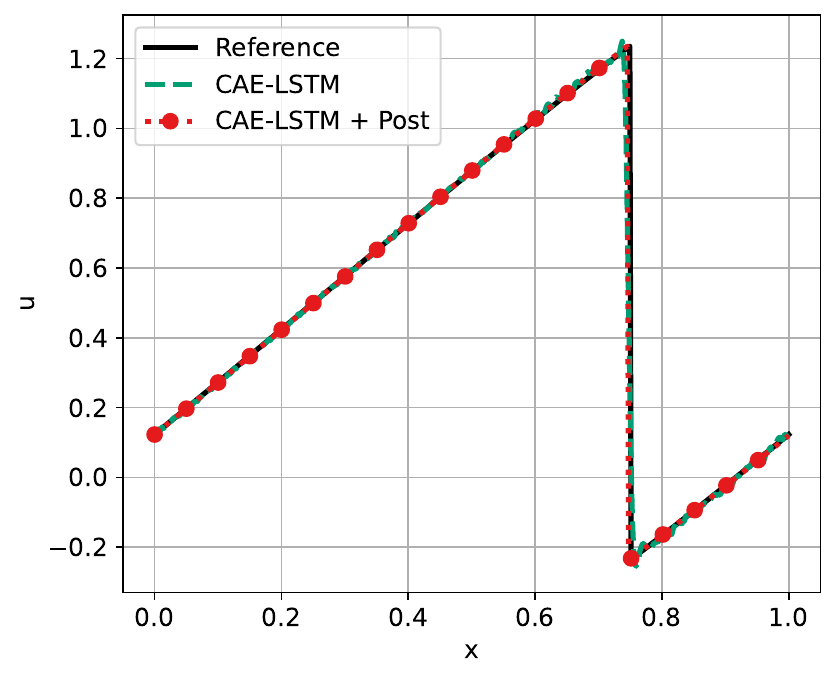}}
    \subfloat[$\sigma=10$]{\includegraphics[width=0.32\textwidth]{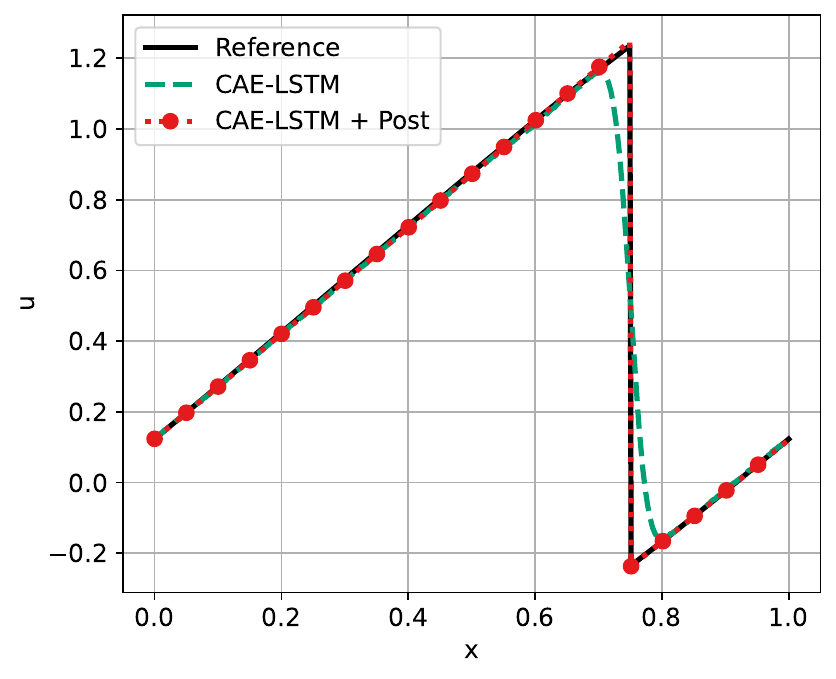}}
    \caption{Example 2. The CAE-LSTM solutions at t = 0.5 under different filter widths, with a fixed reconstruction parameters $m=3$ and $\lambda = 12$.
    }
    \label{eg3_caelstm_t1}
\end{figure}

\begin{figure}[!htbp]
    \centering
    \subfloat[$m=3\quad \lambda=3$]{\includegraphics[width=0.32\textwidth]{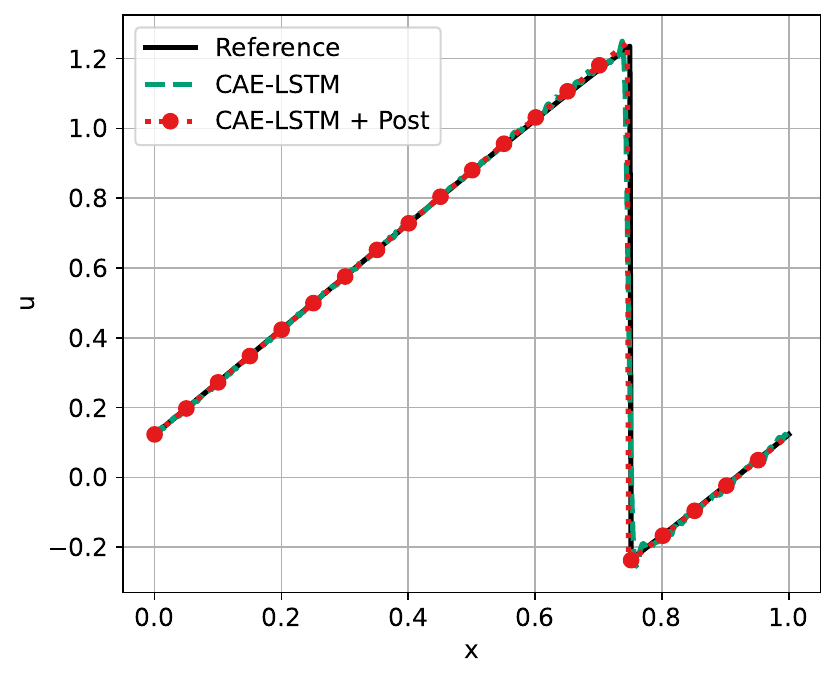}}
    \subfloat[$m=3\quad\lambda=12$]{\includegraphics[width=0.32\textwidth]{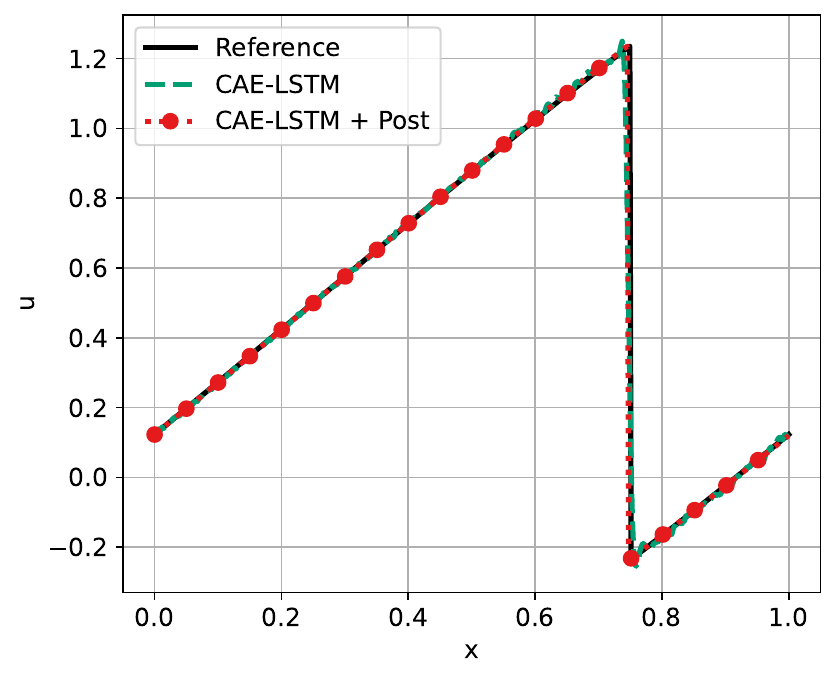}}
    \subfloat[$m=6\quad\lambda=12$]{\includegraphics[width=0.32\textwidth]{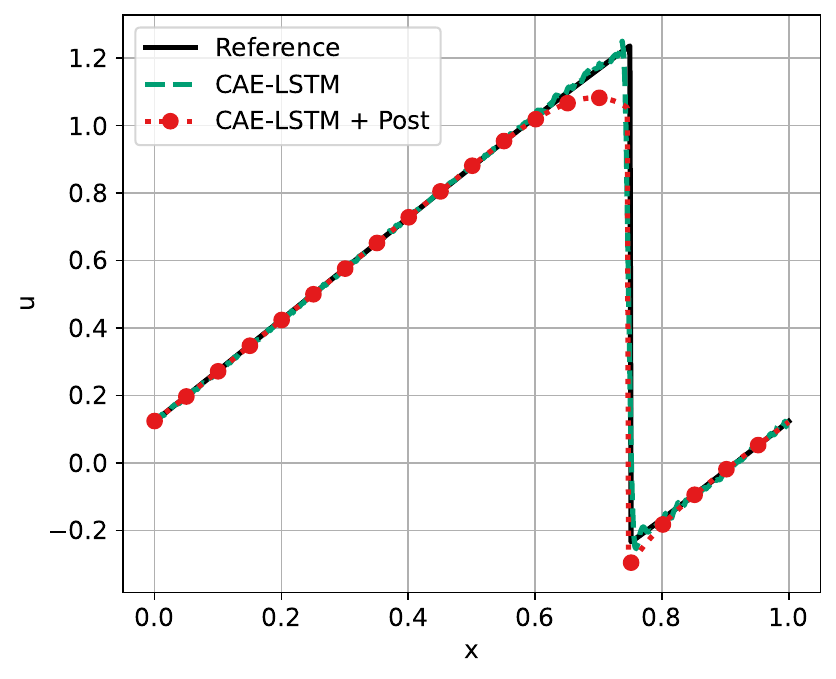}} 
    \caption{Example 2. The CAE-LSTM solutions at $ t = 0.5 $ with different reconstruction parameters, when fixing the filter width to 2 and the latent dimension to 4.
    }
    \label{eg3_caelstm_sigma2_lambda_m}
\end{figure}

The numerical solutions and corresponding errors are presented in Figure~\ref{eg3_caelstm_t1} and Table~\ref{tab_combined_caelstm_t1}, where the reconstruction parameters are set to $m=3$ and $\lambda=12$. The results of study on parameter $\lambda, m$ are reported in Figure~\ref{eg3_caelstm_sigma2_lambda_m}. Similar conclusions as the other ROMs can be made.

\begin{table}[!htbp]
    \centering
    \caption{\label{tab_combined_caelstm_t1}Example 2. Errors of the CAE-LSTM solutions at $ t = 0.5 $ under different filter widths, with a fixed reconstruction parameters $m=3$ and $\lambda = 12$.}
    \begin{tabular}{c|l|c|c|c}
        \hline
        && $\sigma = 0$ & $\sigma = 2$ & $\sigma = 10$\\\hline
        \multirow{2}{*}{Relative Error} 
        &CAE-LSTM  & 1.7165e-02 & 1.3458e-02 & 1.1671e-01 \\ 
        &CAE-LSTM + Post  & 3.4430e-03 & 2.7306e-03 & 6.2476e-03 \\ \hline
        \multirow{2}{*}{Maximum Error} 
        &CAE-LSTM &  8.4343e-02 & 1.1426e-01 & 5.9194e-01 \\ 
        &CAE-LSTM + Post & 8.0519e-03 & 5.0171e-03 & 1.1941e-02 \\ \hline
    \end{tabular}
\end{table}

To further validate the effectiveness of the proposed Gegenbauer post-processing framework for shock-dominated solutions, we additionally consider three classical full-order solvers for the inviscid Burgers' equation: a finite difference (FD) scheme, a finite volume (FV) scheme, and a physics-informed neural network (PINN). In all cases, the final time is fixed at $t=0.5$, at which the solution develops a moving discontinuity.

\begin{figure}[!htbp]
    \centering
    \subfloat[Finite difference]{\includegraphics[width=0.32\textwidth]{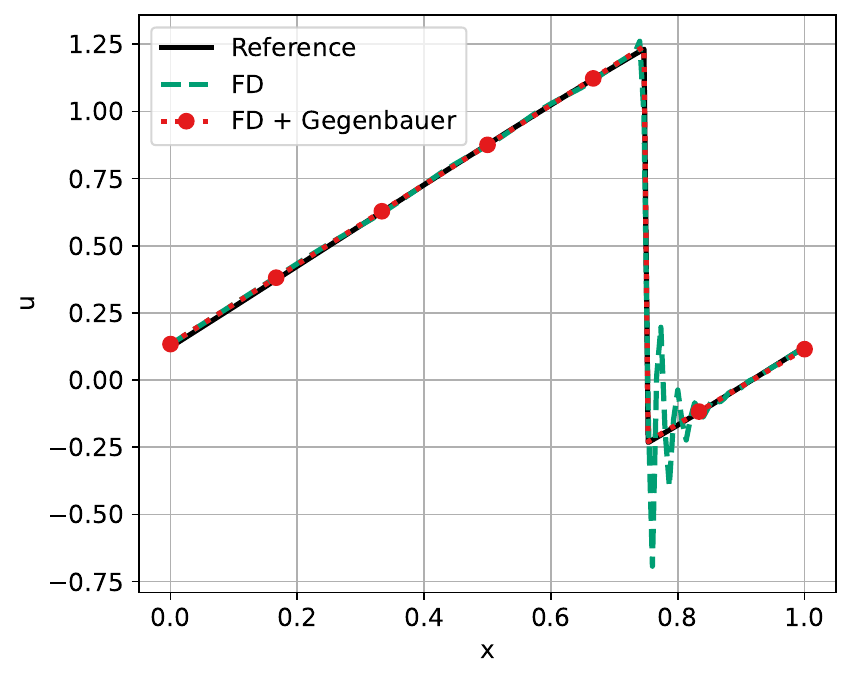}\label{eg2_fd}}
    \subfloat[Finite volume]{\includegraphics[width=0.32\textwidth]{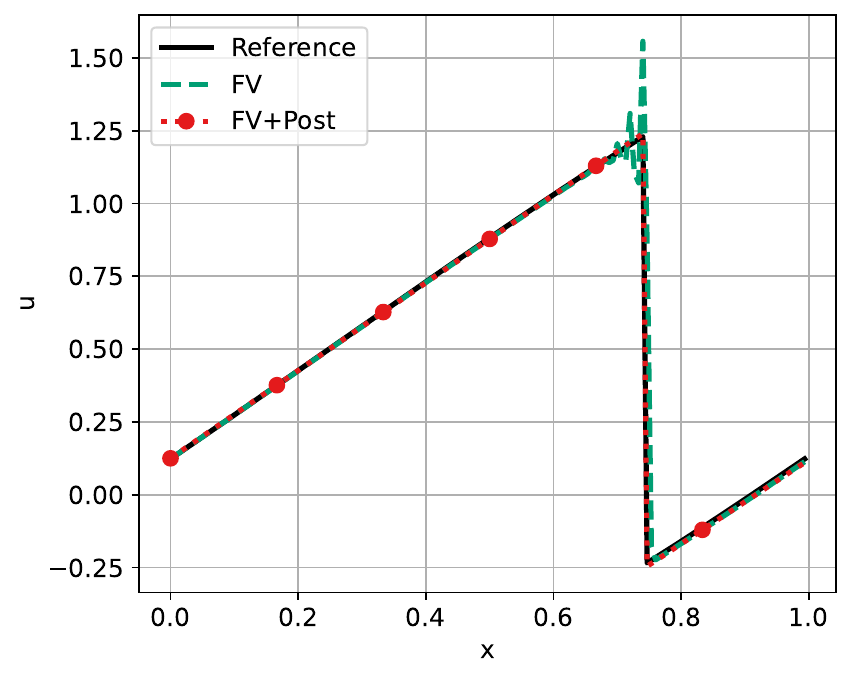}\label{eg2_fv}}
    \subfloat[Physics-informed neural network]{\includegraphics[width=0.32\textwidth]{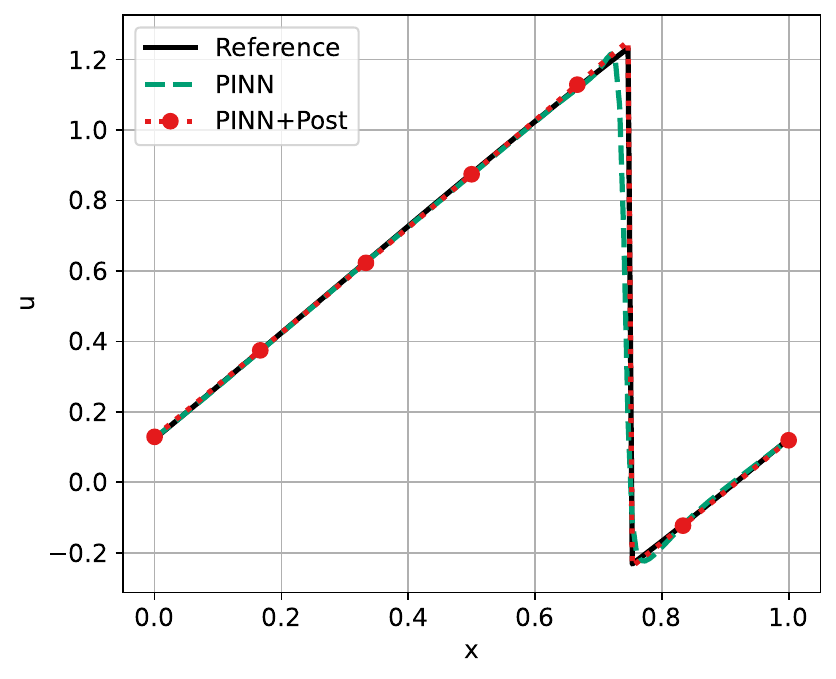}\label{eg2_pinn}} 
    \caption{Example 2. Numerical solutions of the inviscid Burgers' equation at $t=0.5$ obtained by the finite difference, finite volume, and physics-informed neural network methods. }
    \label{eg2_burgers_fd_fv_pinn}
\end{figure}

We first employ a conservative finite difference scheme with global Lax-Friedrichs flux splitting. The spatial domain $[0,1]$ is discretized using $N=151$ uniform grid points. A third-order upwind scheme is used to compute the numerical flux, and periodic boundary conditions are imposed through index wrapping. 
Time integration is carried out using a third-order Runge-Kutta method with an adaptive time step satisfying the CFL condition with $\mathrm{CFL}=0.2$.
The resulting numerical solution is then post-processed using the edge-detected Gegenbauer reconstruction with parameters $(\lambda,m)=(3,1)$, as shown in Figure~\ref{eg2_fd}.

We next employ a conservative finite volume scheme on a uniform mesh with $N=150$ cells. A second-order Lax-Wendroff flux is used for spatial discretization. Time integration is performed using a third-order Runge-Kutta method with a CFL condition $\mathrm{CFL}=0.3$. Periodic boundary conditions are imposed.
The resulting numerical solution is post-processed using the edge-detected Gegenbauer reconstruction with $(\lambda,m)=(3,1)$. As shown in Figure~\ref{eg2_fv}, the reconstruction effectively suppresses the spurious oscillations generated by the second-order scheme and produces a sharper transition near the shock.

Finally, we approximate the solution using a physics-informed neural network. The network takes $(x,t)$ as input and outputs $u_\theta(x,t)$. It consists of a fully connected feedforward architecture with three hidden layers of width 64 and $\tanh$ activations.
The training loss is defined as
\begin{equation*}
    \mathcal{L}
    = \mathcal{L}_{\mathrm{PDE}}
    + \mathcal{L}_{\mathrm{IC}}
    + \mathcal{L}_{\mathrm{BC}}
    + 10\,\mathcal{L}_{\mathrm{data}},
\end{equation*}
where $\mathcal{L}_{\mathrm{PDE}}=\|u_t+u u_x\|^2$ is evaluated at $N_r=15000$ collocation points, $\mathcal{L}_{\mathrm{IC}}$ enforces the initial condition with $N_{ic}=300$ samples, $\mathcal{L}_{\mathrm{BC}}$ imposes periodic boundary conditions using $N_{bc}=300$ samples, and $\mathcal{L}_{\mathrm{data}}$ matches reference data at $t=0.4$. The network is trained with the Adam optimizer for 4000 epochs.
The resulting PINN solution at $t=0.5$ exhibits a smeared shock (see Figure~\ref{eg2_pinn}). After applying the edge-detected Gegenbauer reconstruction with $(\lambda,m)=(3,1)$, the discontinuity is significantly sharpened and the local error is reduced.

Overall, these results demonstrate that the edge-detected Gegenbauer reconstruction is not restricted to reduced-order modeling frameworks. 
In addition to its effectiveness for ROMs, it also improves solutions produced by classical numerical discretizations such as finite difference and finite volume methods, as well as machine-learning approaches such as physics-informed neural networks. 
For FD and FV solvers, the reconstruction effectively suppresses spurious oscillations near discontinuities, while for PINN solutions it mitigates the overly smooth transition, resulting in a sharper profile near the shock. 
This highlights that the proposed reconstruction is a solver-independent and broadly applicable post-processing tool for shock-dominated nonlinear hyperbolic problems.

\subsection{Example 3: parametric viscous Burgers' equation}
We consider the viscous Burgers' equation with periodic boundary conditions which can be represented as
\begin{equation}
\begin{cases}
    \begin{aligned}
        &u_t + uu_x = \nu u_{xx}, \qquad x \in [0, 1],\\
        &u(x,0) = \frac{x}{1+\sqrt{\frac{1}{c}}\exp (Re\frac{x^2}{4})},\qquad c = \exp(\frac{Re}{8}),\quad Re = \frac{1}{\nu},
    \end{aligned}
\end{cases}
\end{equation}
where the parameter $Re = \frac{1}{\nu}$ is the Reynolds number and $\nu$ is kinematic viscosity coefficient. It is widely recognized that the aforementioned equation can produce large gradient when $\nu$ is sufficiently small, even when the initial conditions are smooth, due to the dominance of advection effects.
The analytical solution exists and is given by
\begin{equation}
    u(x,t) = \frac{\frac{x}{t+1}}{1+\sqrt{\frac{t+1}{c}}\exp (Re\frac{x^2}{4t+4})},
\end{equation}

In the following, we investigate G-ROM and CAE-Galerkin for this parametric problem. For both methods, the spatial domain is discretized into 256 uniform nodes in the interval $[0, 1]$. During snapshot construction, the parameter $ Re $ spans the range from 100 to 1900 with a step size of 100, resulting in a total of 19 parameter values. For each parameter, the final time is set to 1, and a snapshot is taken every 0.01 time interval.  We do not use pre-processing filters in this example. We then test the ROMs with $\text{Re} = 1530$ and $\text{Re} = 2500$ at $t=0.5$, where the former is within the parameter range used for training in the offline stage, while the latter exceeds this range.

\begin{figure}[!htbp]
    \centering
    \subfloat[$r=20$]{\includegraphics[width=0.33\textwidth]{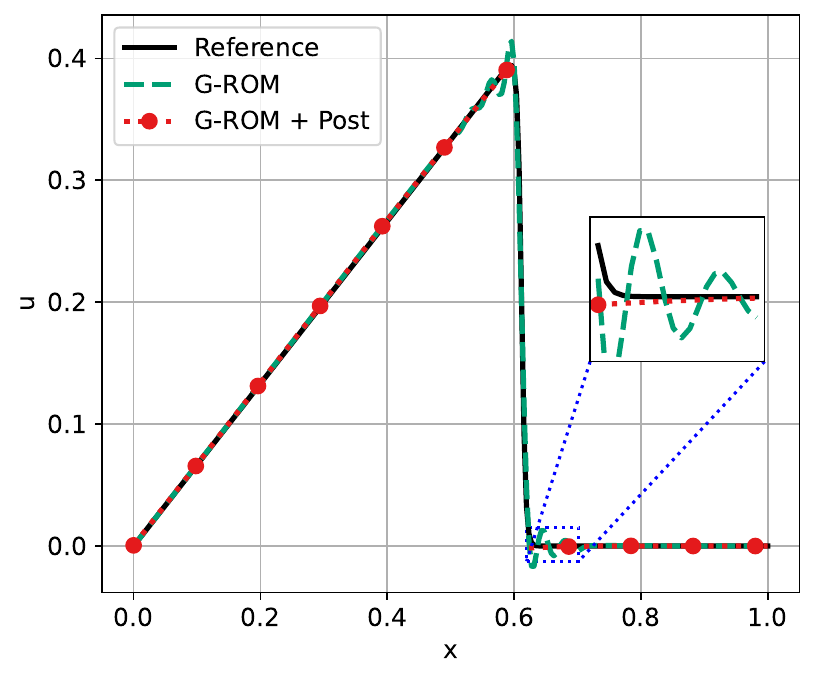} \label{inviscid_burgers_pod_t1_a}}
    \subfloat[$r=30$]{\includegraphics[width=0.33\textwidth]{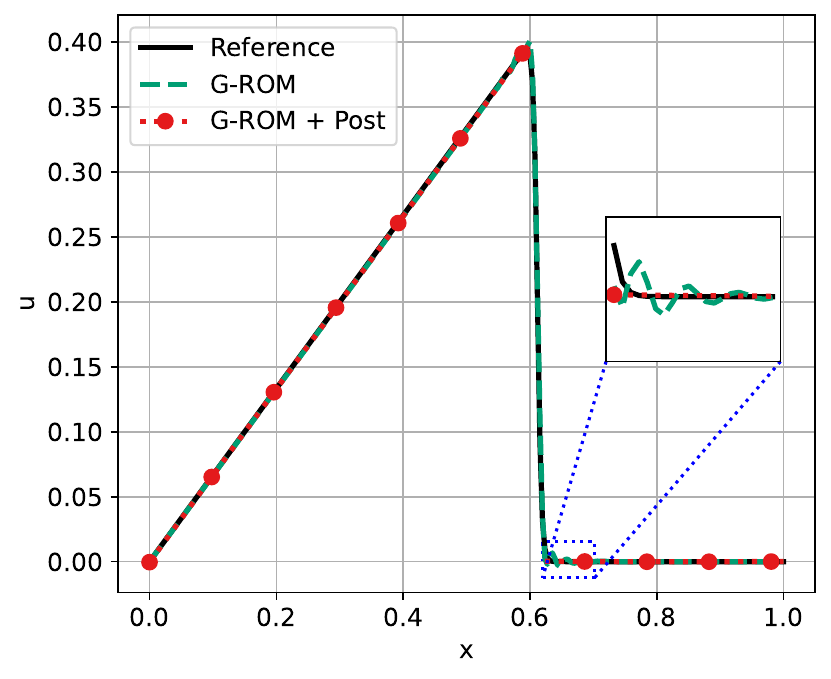}}
    \subfloat[$r=40$]{\includegraphics[width=0.33\textwidth]{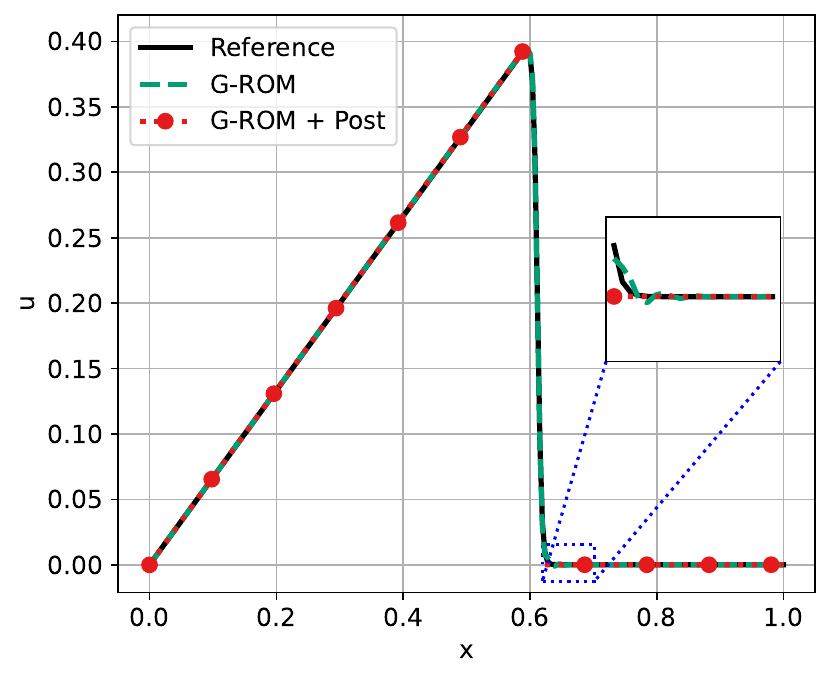}}
    \caption{Example 3. Numerical solutions of G-ROM at $t = 0.5$ under different reduced dimensions, with $Re = 1530$.    }
    \label{inviscid_burgers_pod_t1}
\end{figure}

\begin{figure}[!htbp]
    \centering
    \subfloat[$r=20$]{\includegraphics[width=0.33\textwidth]{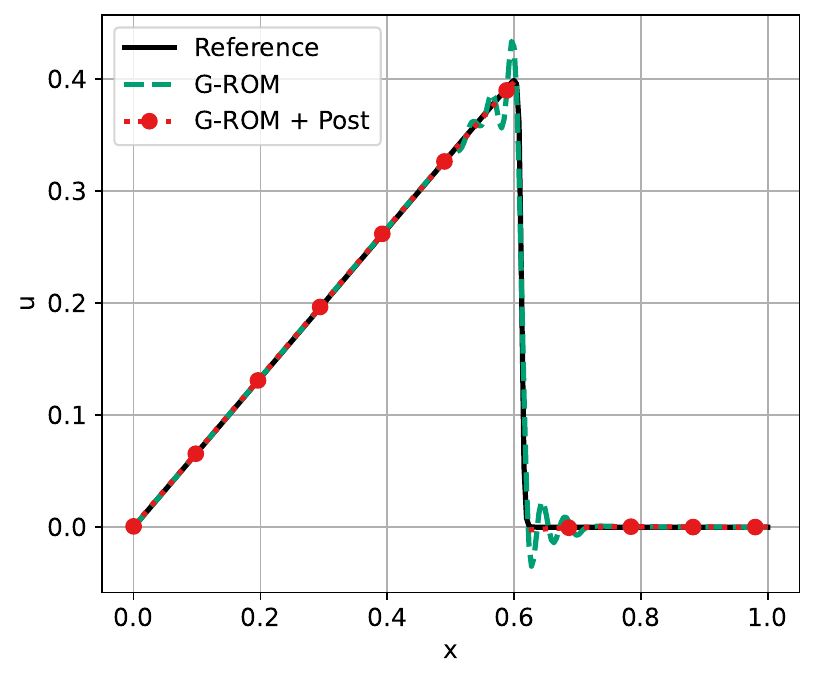}}
    \subfloat[$r=30$]{\includegraphics[width=0.33\textwidth]{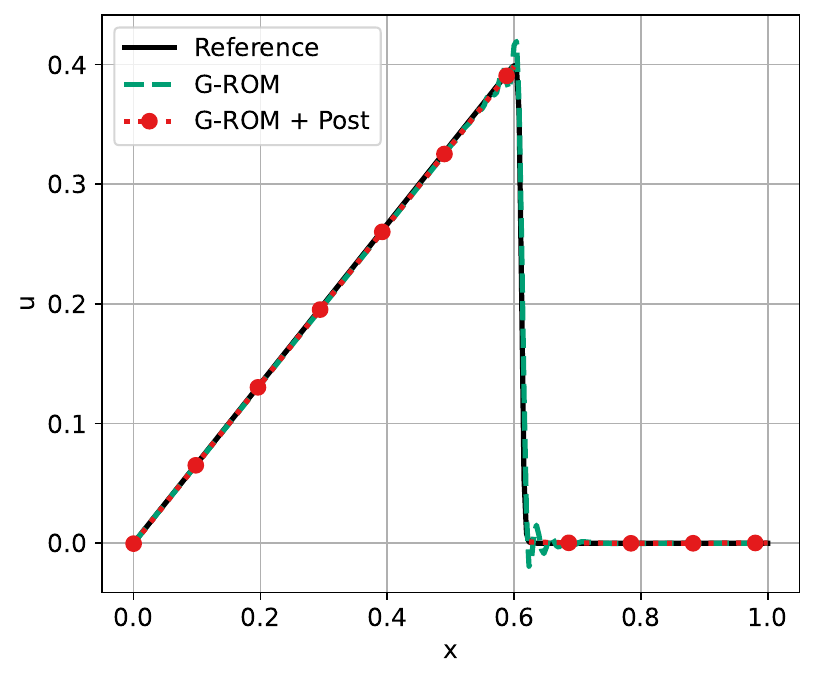}}
    \subfloat[$r=40$]{\includegraphics[width=0.33\textwidth]{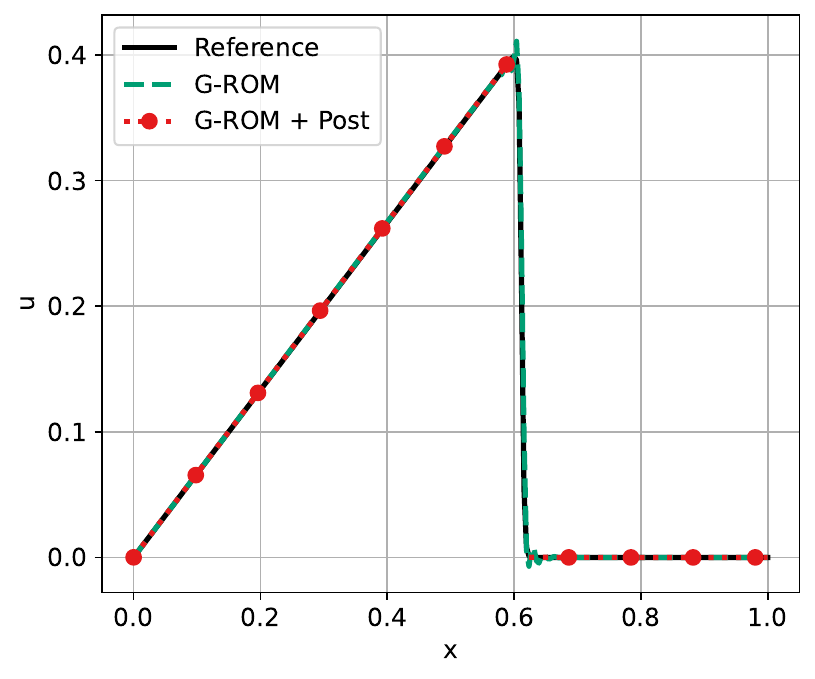}}
    \caption{Example 3. Numerical solutions of G-ROM at $t = 0.5$ under different reduced dimensions, with $Re = 2500$.    }
    \label{inviscid_burgers_pod_t1_2500}
\end{figure}

When solving with G-ROM online, we apply the exponential filter in the time evolution.  Since this problem does not exhibit discontinuity but contains points with large gradients, we identify these points as discontinuities. Two points on each side of the identified discontinuity are not used for reconstruction. For G-ROM solutions, we set the Gegenbauer reconstruction parameters to be $\lambda = m = 3$. Figure~\ref{inviscid_burgers_pod_t1} corresponds to the results for $\text{Re} = 1530$ with various reduced dimensions, while Figure~\ref{inviscid_burgers_pod_t1_2500} presents the results for $\text{Re} = 2500$ obtained by G-ROM.  Table~\ref{tab_inviscid_burgers_grom} contains the numerical errors before and after post-processing.

Analyzing the results, we observe that the Gegenbauer reconstruction is successful in suppressing the numerical oscillations.  
However, as the reduced dimension increases, the accuracy of G-ROM improves, and numerical oscillations are less apparent. From the data presented in Table~\ref{tab_inviscid_burgers_grom}, we notice that in this example, the reconstruction is not effective in significantly reducing numerical errors. 
This behavior can be attributed to the nature of the viscous Burgers’ equation at high Reynolds numbers: although the solution is mathematically smooth, it develops extremely steep, shock-like gradients rather than true discontinuities. Since the Gegenbauer reconstruction is theoretically designed for piecewise analytic functions with distinct discontinuities, treating these steep-gradient regions as “pseudo-discontinuities” introduces challenges. In particular, the limited number of grid points available within these narrow regions restricts the local reconstruction accuracy and may slightly degrade the overall error while still effectively removing oscillations near the sharp gradients.
Overall, these results indicate that the Gegenbauer reconstruction remains highly effective for suppressing oscillations in both discontinuous and shock-like regimes, while its quantitative error reduction is more pronounced for problems involving true discontinuities.
Actually, it may be possible to improve accuracy by choosing different values of $\lambda, m,$ but we did not experiment with this.

\begin{table}[!htbp]
    \centering
    \caption{\label{tab_inviscid_burgers_grom}Example 3. Errors of the G-ROM solutions at $ t = 0.5 $ under different reduced dimensions, with $Re = 1530$ and $Re = 2500$. }
    \begin{tabular}{c|l|c|c|c}
        \hline
        \multicolumn{5}{c}{Relative Error}\\ \hline
        & & $r = 20$ & $r = 30$ & $r = 40$\\\hline
        \multirow{2}{*}{$Re = 1530$} 
        & G-ROM  & 1.8800e-02 & 7.3482e-03 & 2.2165e-03 \\
        & G-ROM+ Post  & 6.1549e-03 & 5.6920e-03 & 5.3009e-03 \\\hline
        \multirow{2}{*}{$Re = 2500$} 
        & G-ROM  &  1.8081e-01 &  1.7411e-01 & 1.0687e-01 \\
        & G-ROM + Post  &1.8059e-01 & 1.7709e-01 & 1.1550e-01 \\\hline
        
        \multicolumn{5}{c}{Maximum Error}\\ \hline
        & & $r=20$  & $r=30$ & $r=40$\\\hline
        \multirow{2}{*}{$Re = 1530$} 
        & G-ROM & 1.9875e-02 & 1.0867e-02 & 1.6051e-03\\
        & G-ROM + Post  & 1.1398e-02 & 9.4667e-03 & 1.0580e-02\\\hline
        \multirow{2}{*}{$Re = 2500$} 
        & G-ROM  & 1.3894e-01 & 1.4308e-01 & 1.1494e-01\\
        & G-ROM + Post  &1.3960e-01 & 1.4226e-01 & 1.1982e-01 \\\hline
    \end{tabular}
\end{table} 

Next, we conduct a numerical study of the post-processing CAE-Galerkin solution.  
For CAE network, we use the Adam optimizer with a learning rate of $\eta = 10^{-4}$. The fraction of snapshots used for validation is $\omega = 0.1$. The number of minibatches is determined by a fixed batch size of $m_I(i) = 32$, and the maximum number of epochs is set to $n_{\text{epoch}} = 1000$.  The model employs early stopping and model checkpointing during training to monitor validation loss and ensure that the best model weights are saved. In the encoder, we apply a stride of length 1 $(s=1)$ for all convolutional layers, with pooling operations that downsample the spatial dimensions by a factor of 2 in each layer. In the decoder, we use transposed convolutional layers with a stride of 2 $(s=2)$ to upsample the data back to the original resolution. Each convolutional and transposed convolutional layer is followed by a nonlinear activation. 
For the architecture of CAE, the encoder network is composed of four convolutional layers with increasing filter sizes to extract hierarchical feature representations. Specifically, the first, second, third, and fourth convolutional layers contain 8, 16, 32, and 64 filters, respectively. Conversely, the decoder network reconstructs the input data using a series of four transposed-convolutional layers. The number of filters in these layers decreases progressively, with 64 filters in the first layer, 32 in the second, 16 in the third, and a single filter in the final layer to produce the reconstructed output. This symmetric structure enables the model to effectively learn compressed representations while maintaining essential features for reconstruction.

We train three separate CAE network with latent space dimensions $r$ of 2, 4, and 8, respectively. As the latent dimension $r$ increases, the total number of network parameters also increases, ranging from 33459 for $r=2$ to 45753 for $r=8$. Furthermore, the fitting capability of the network improves with increasing latent dimension, with the training error decreasing from $2.8985 \times 10^{-6}$ for $r=2$ to $9.7008 \times 10^{-7}$ for $r=8$, and the validation error decreasing from $2.5597 \times 10^{-6}$ to $8.7635 \times 10^{-7}$.  

\begin{figure}[!htbp]
    \centering
    \subfloat[$r=2$]{\includegraphics[width=0.33\textwidth]{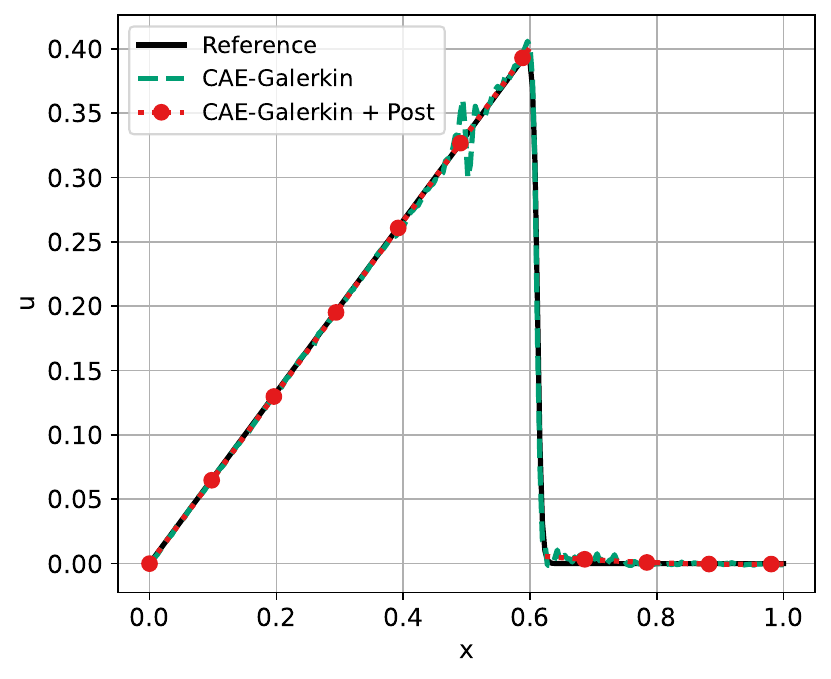}}
    \subfloat[$r=4$]{\includegraphics[width=0.33\textwidth]{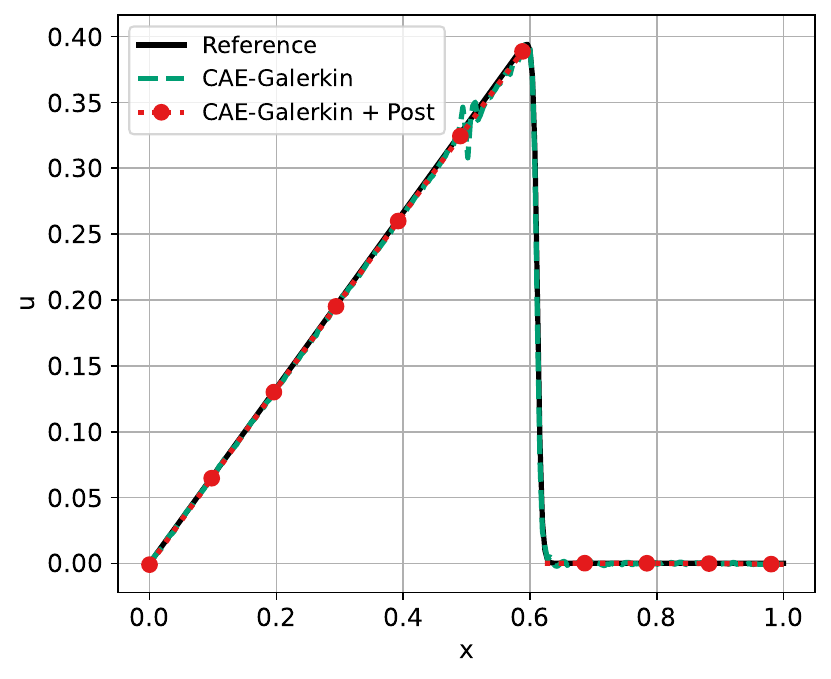}}
    \subfloat[$r=8$]{\includegraphics[width=0.33\textwidth]{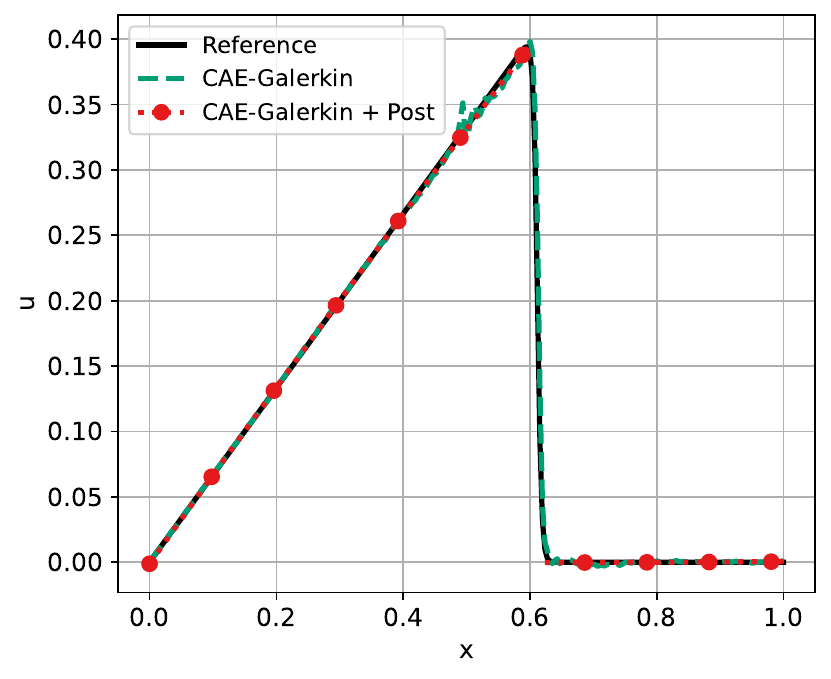}}
    \caption{Example 3. Numerical solutions of CAE-Galekin at $t=0.5$ under different latent space dimensions, with a fixed Reynolds number of 1530.
    }
    \label{CAE_1530_gegen}
\end{figure}

\begin{figure}[!htbp]
    \centering
    \subfloat[$r=2$]{\includegraphics[width=0.33\textwidth]{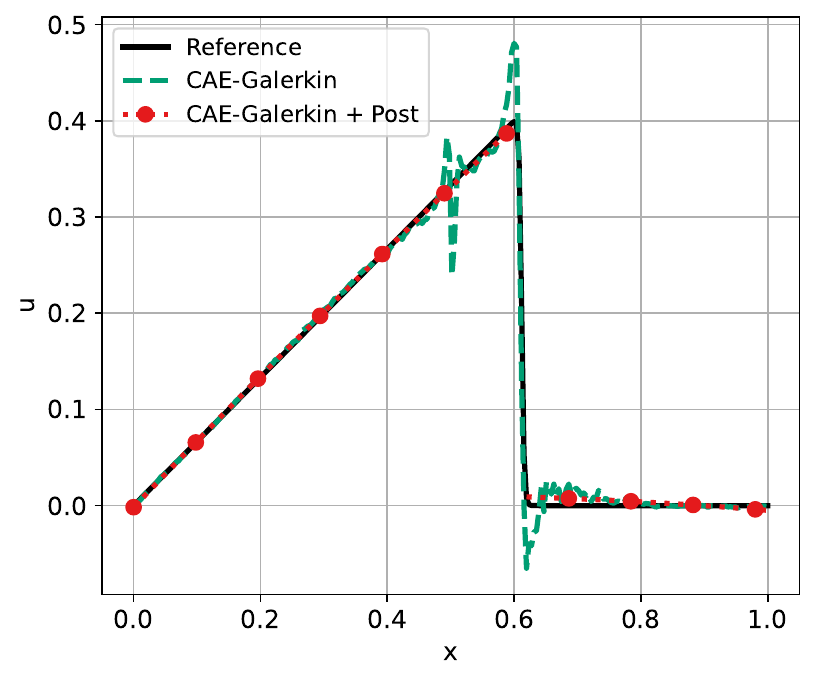}}
    \subfloat[$r=4$]{\includegraphics[width=0.33\textwidth]{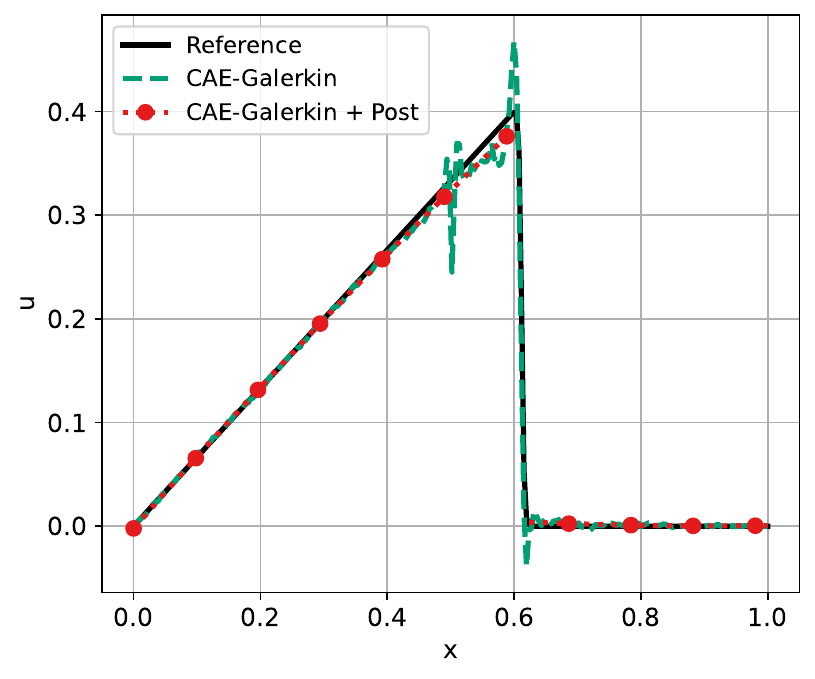}}
    \subfloat[$r=8$]{\includegraphics[width=0.33\textwidth]{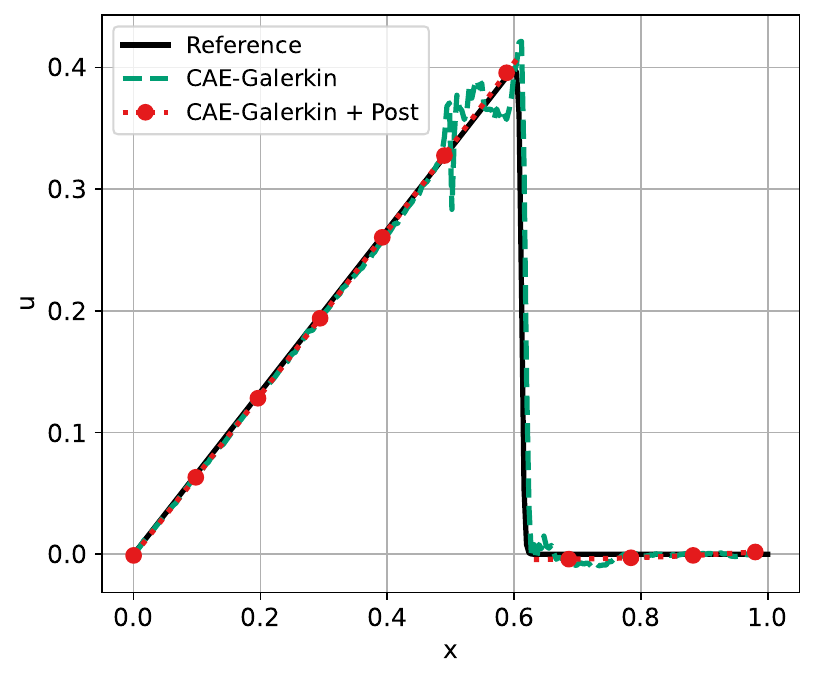}}
    \caption{Example 3. Numerical solutions of CAE-Galekin at $t=0.5$ under different latent space dimensions, with a fixed Reynolds number of 2500.
    }
    \label{CAE_2500_gegen}
\end{figure}

We evaluate the trained CAE-Galerkin model by testing it at $t=0.5$ with $\text{Re} = 1530$ and $\text{Re} = 2500$, and the results are presented in Figures ~\ref{CAE_1530_gegen} and ~\ref{CAE_2500_gegen}, respectively, which is obtained by Gegenbauer reconstruction parameters $\lambda = 1$ and $m = 2$. The conclusions are similar before. The Gegenbauer reconstructions can eliminate the numerical oscillations, but we do not expect strong performance in accuracy enhancement.

\subsection{Example 4: 2D linear equation}
\label{eg4}
We consider the two-dimensional linear problem that corresponds to the solid body rotation:
\begin{equation}
    u_t + yu_x - xu_y = 0 , \qquad (x,y)\in [-1,1]\times [-1,1],
\end{equation}
with periodic boundary condition in each direction and the initial condition 
\begin{equation}
    u(x,y,0) = u_0(x,y)=
    \begin{cases}
        1+ \sin(\pi x)\sin (\pi y),&\qquad \text{if} \ \frac{x^2}{0.7^2}+\frac{y^2}{0.5^2}\le 1,\\
        0, &\qquad \text{otherwise}.
    \end{cases}
\end{equation}
The true solution to the two-dimensional problem is given by $ u(x, y, t) = u_0(x\cos t - y\sin t, x\sin t + y \cos t) $. In the following, we will look at the solution obtained by POD-Galerkin and its post-processing.

\begin{figure}[!htbp]
    \centering
    \includegraphics[width=0.55\linewidth]{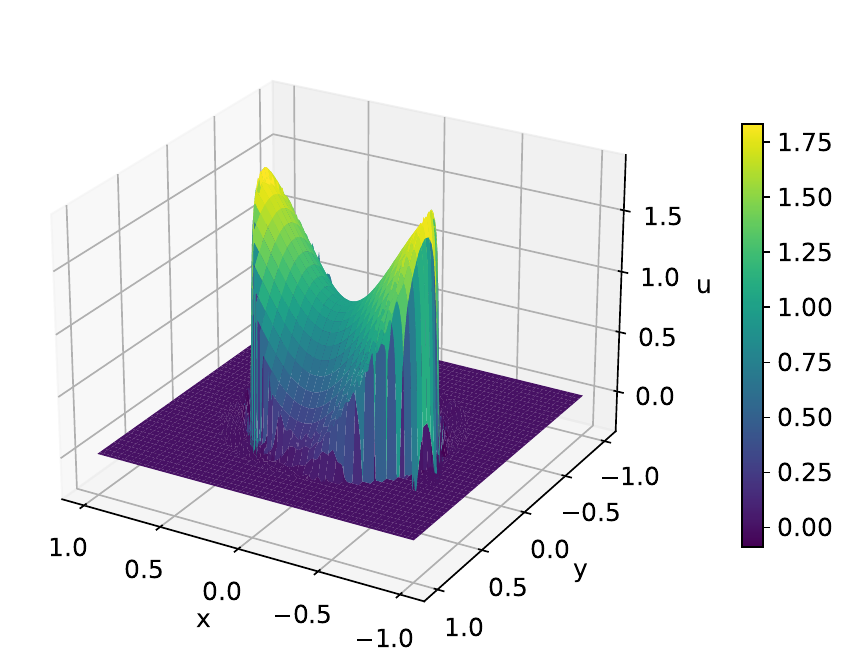}
    \caption{Example 4.  The G-ROM solution at $t=\pi / 4$, with the reduced dimension of 30.}
    \label{eg4_2d_grom}
\end{figure}

\begin{figure}[!htbp]
    \centering
    \subfloat[$y = -0.6094$]{\includegraphics[width=0.33\textwidth]{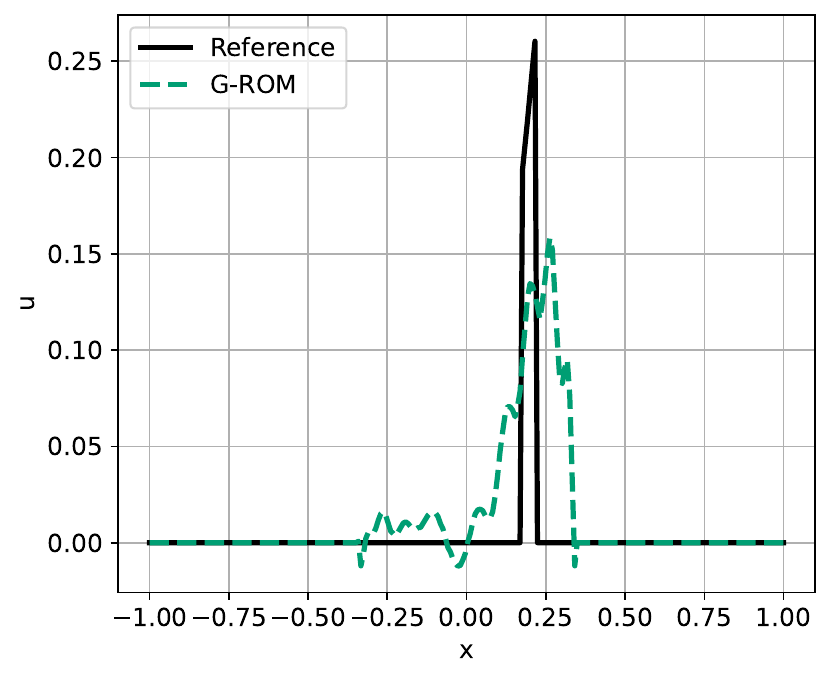}}
    \subfloat[$y = -0.4531$]{\includegraphics[width=0.33\textwidth]{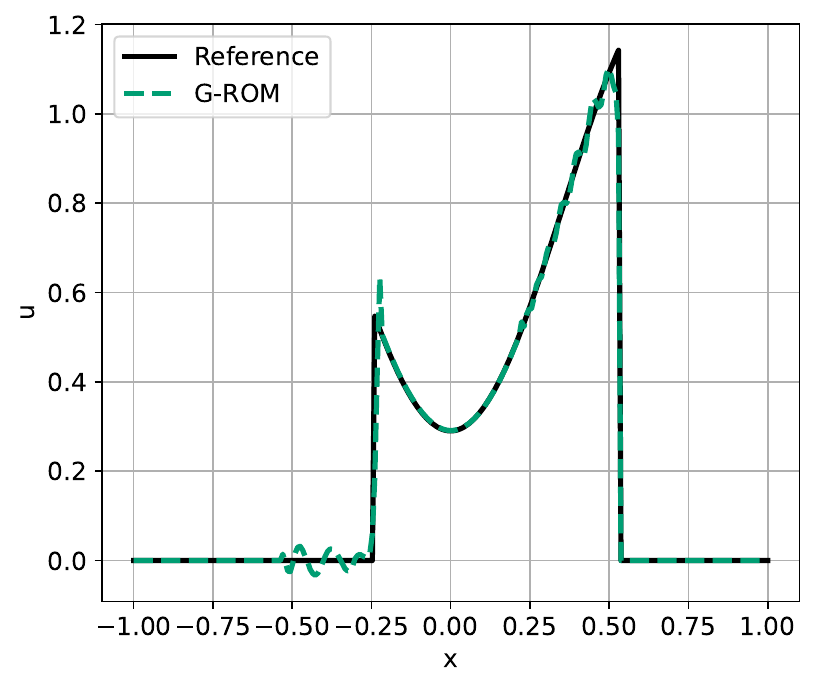}}
    \subfloat[$ y = 0.0 $]{\includegraphics[width=0.33\textwidth]{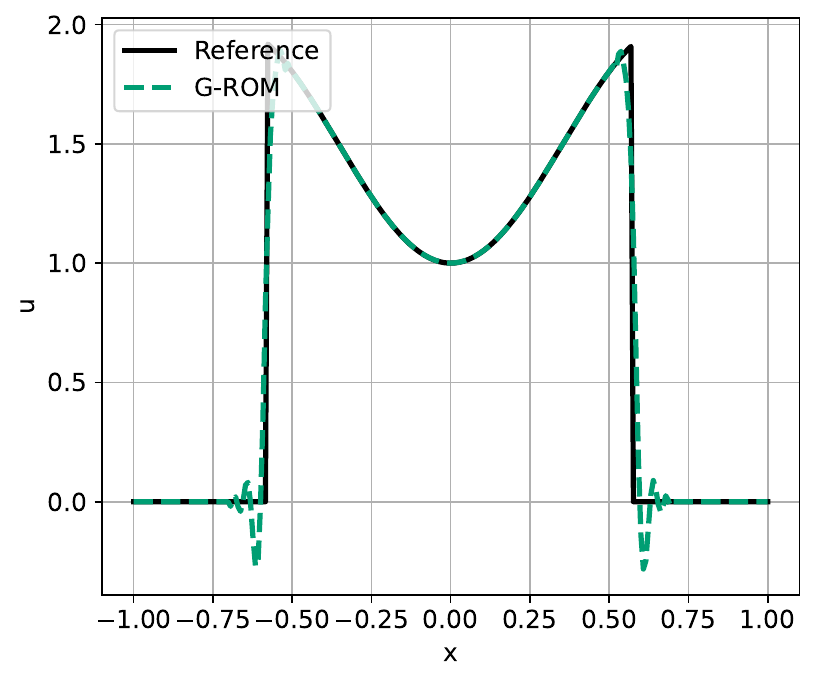}}
    \caption{Example 4. The G-ROM solution of $ t = \pi / 4 $ profiles along $ y = -0.6094 $, $ y = -0.4531$ and $ y =0.0 $.    }
    \label{eg4_2d_profiles_before}
\end{figure}

For the computation, the spatial domain is discretized uniformly with 256 points in both dimensions within the interval $[-1, 1]$. To construct the snapshot matrix, we sample 600 uniformly distributed time points in $[0, 2\pi]$.
For the snapshots, no filtering is applied. 
The POD method uses reduced space dimension of $30$ and computes the solution at $T = \pi / 4$ using G-ROM, which is plotted in Figure~\ref{eg4_2d_grom}. 
Additionally, several cross-sectional plots are presented in Figure~\ref{eg4_2d_profiles_before}, where we observe significant numerical oscillations near the discontinuity.  

\begin{figure}[!htbp]
    \centering
    \subfloat[$y=-0.6094$]{\includegraphics[width=0.33\textwidth]{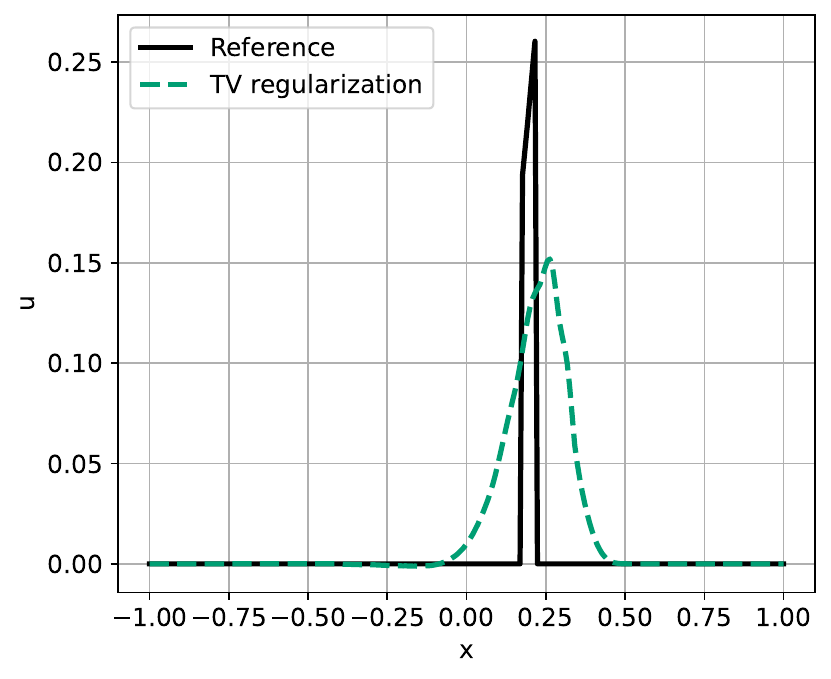}}
    \subfloat[$y=-0.4531$]{\includegraphics[width=0.33\textwidth]{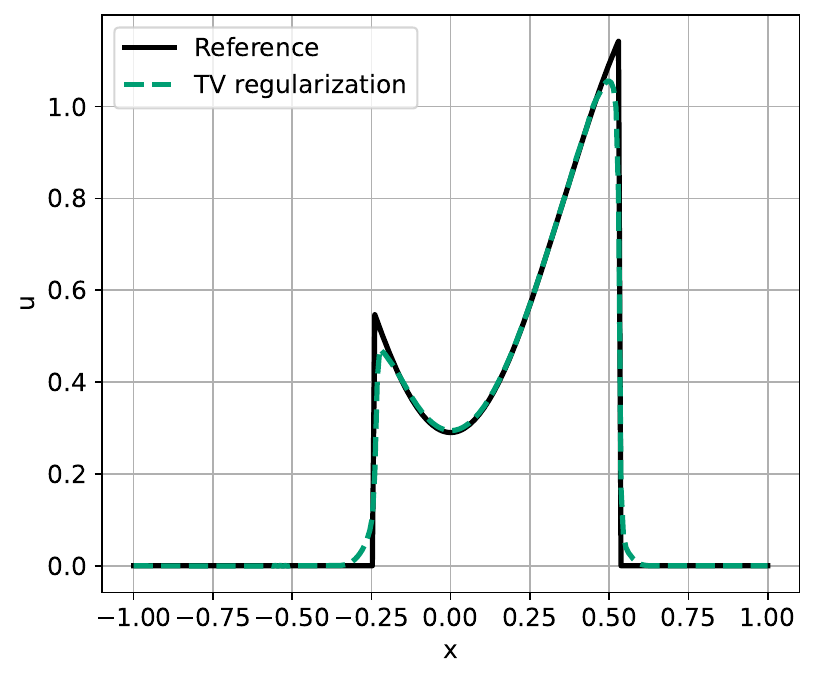}}
    \subfloat[$y= 0.0$]{\includegraphics[width=0.33\textwidth]{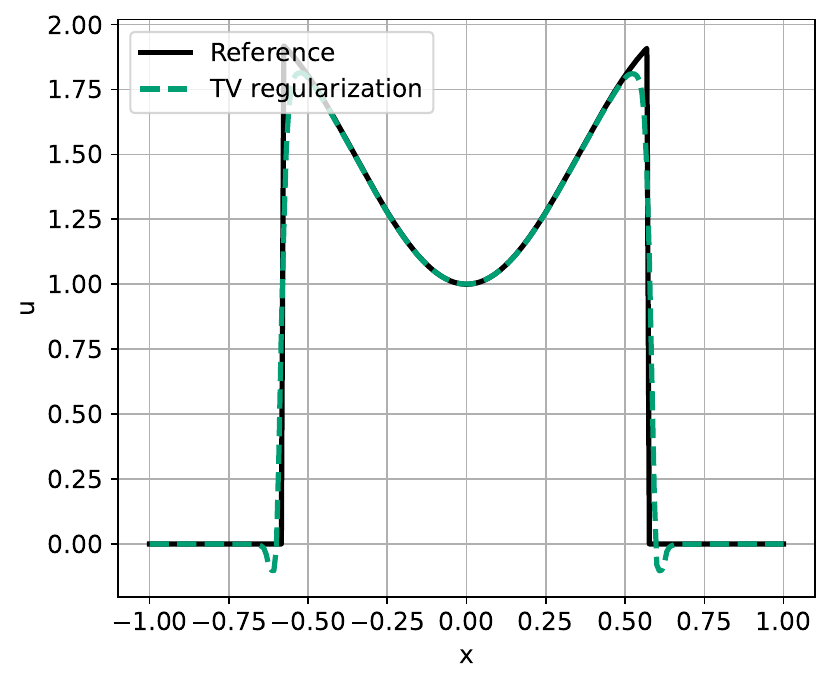}}
    \caption{Example 4. The G-ROM solution of $ t = \pi / 4 $, applied with TV regularization ( $\lambda_{TV}=0.2$ ),  profiles  along $ y = -0.6094 $, $ y = -0.4531$ and $ y =0.0 $.
    }
    \label{eg4_2d_profiles_tv}
\end{figure}

We then apply the standard TV regularization with $\lambda_{TV} = 0.2$ to the G-ROM solution. The results are presented in Figure   \ref{eg4_2d_profiles_tv}. While it is clear that the TV regularization is effective in suppressing spurious oscillations, the numerical results are smeared out, and the deviation from the reference solution is obvious at $y=-0.6094$.

\begin{figure}[!htbp]
    \centering
    \includegraphics[scale = 0.45]{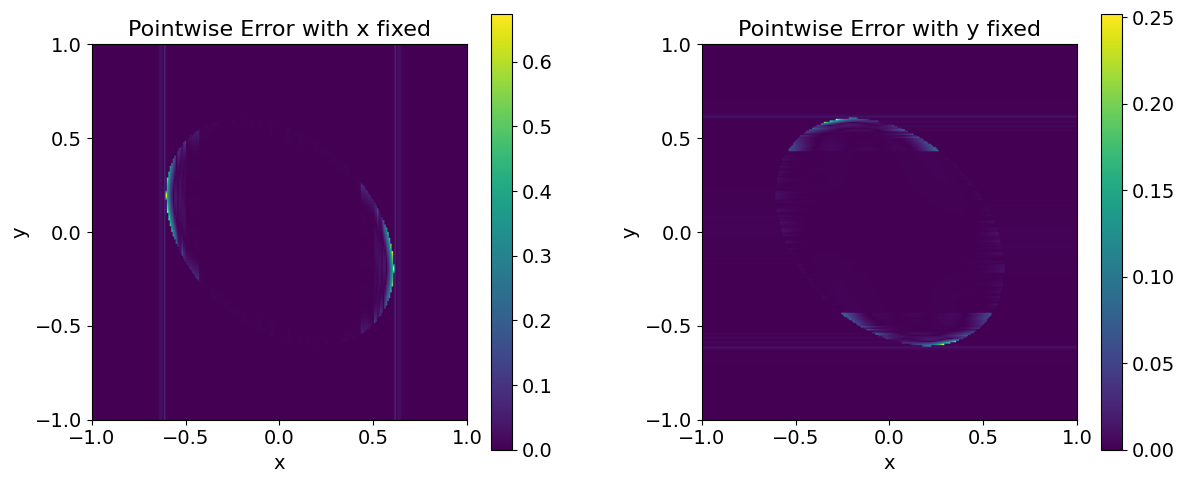}
    \caption{Example 4. The pointwise error of the Gegenbauer polynomial reconstruction at $t=\pi/4$ when $ x $ or $ y $ is fixed.}
    \label{eg4_x_y}
\end{figure}

We then post-process the G-ROM solution by Gegenbauer reconstruction. First, we follow Algorithm~\ref{alg_1} and Algorithm~\ref{alg_2} to perform a line-by-line reconstruction.  
In this example, we take the following heuristic rules for choosing the parameter value.
For the region outside the ellipse when the solution is flat, we set $\lambda = 2$ and $m = 0$ for the flat regions. Within the ellipse, when the number of grid points contained in the identified continuous interval is less than $0.4N$ (where $N$ denotes the total number of grid points along one direction), we choose $\lambda = 3$ and $m = 3$. Otherwise, we let $\lambda = 3$ and $m = 6$.  
We note that as mentioned in the 1D experiments, the choice of $\lambda, m$ is important and solution-dependent. Here, we take the empirical rule while commenting that a small change in parameters will not affect the reconstruction quality but a large deviation will deteriorate the solution results.
The numerical results presented in Figure \ref{eg4_x_y} indicate that for Gegenbauer polynomial reconstruction along fixed $x$ or fixed $y$, effective correction is achieved in most regions, except when the analyticity interval is too small. 

\begin{figure}[!htbp]
    \centering
    \includegraphics[width=0.95\linewidth]{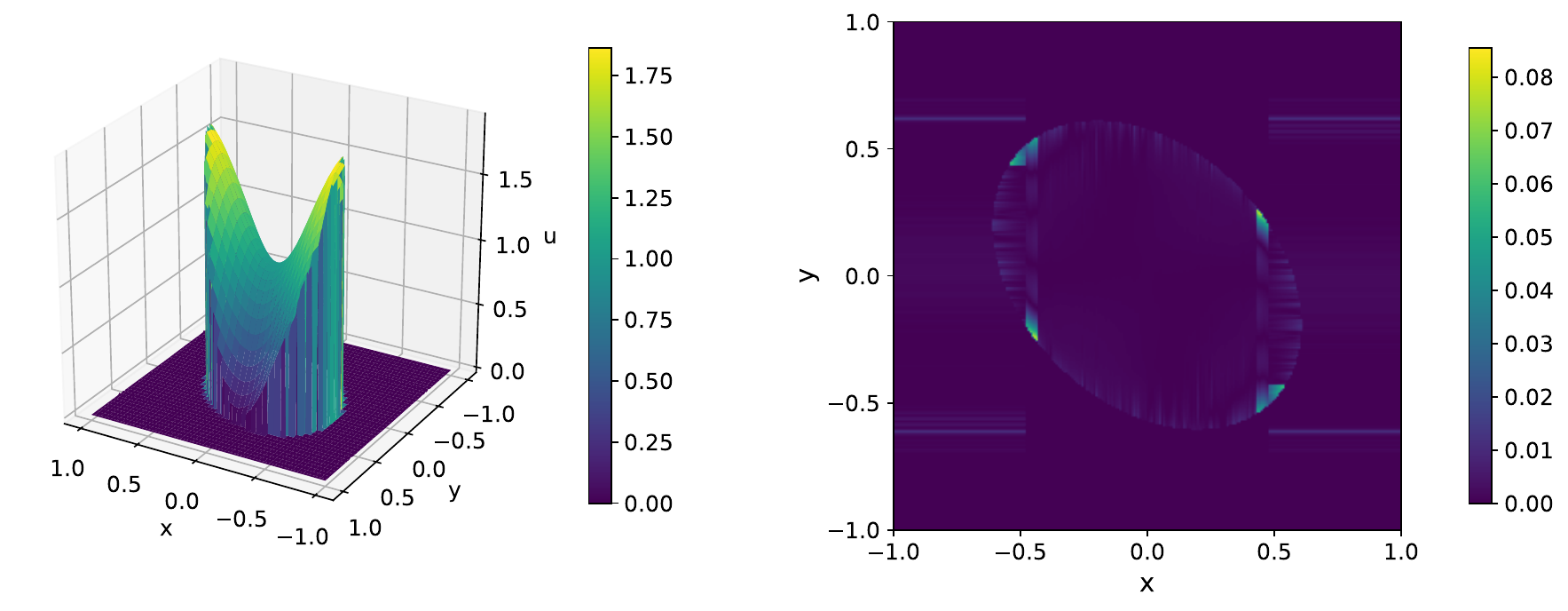}
    \caption{Example 4: The Gegenbauer reconstruction of the G-ROM solution at $t=\pi / 4$, and corresponding pointwise error.   
    }
    \label{eg4_2d_gegen}
\end{figure}

\begin{figure}[!htbp]
    \centering
    \subfloat[$y = -0.6094$]{\includegraphics[width=0.33\textwidth]{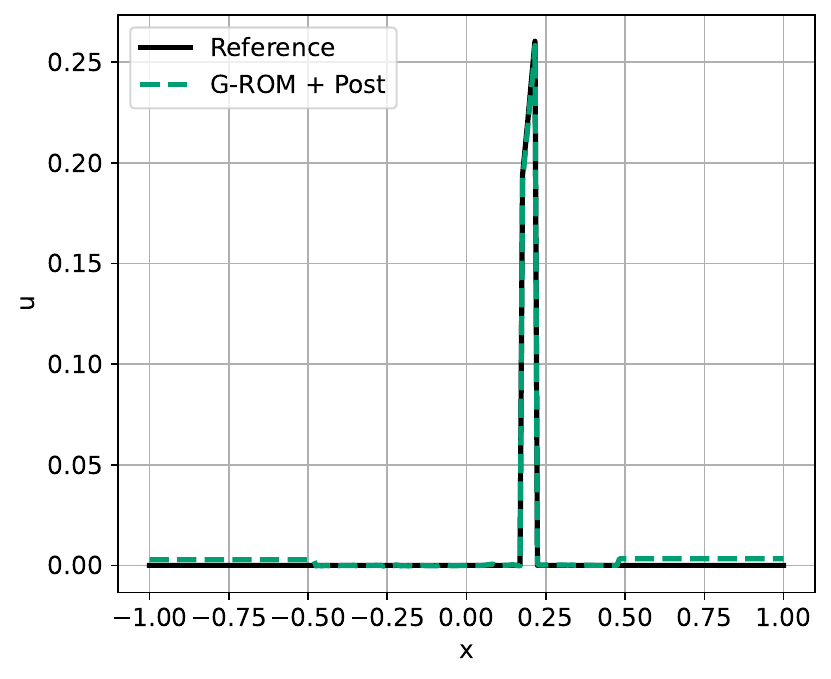}}
    \subfloat[$y = -0.4531$]{\includegraphics[width=0.33\textwidth]{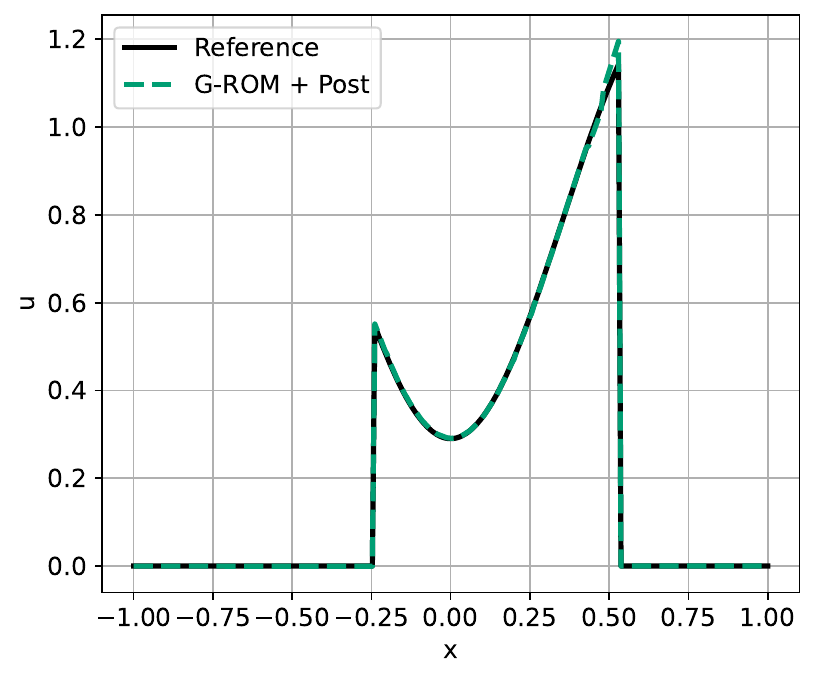}}
    \subfloat[$ y =0.0 $]{\includegraphics[width=0.33\textwidth]{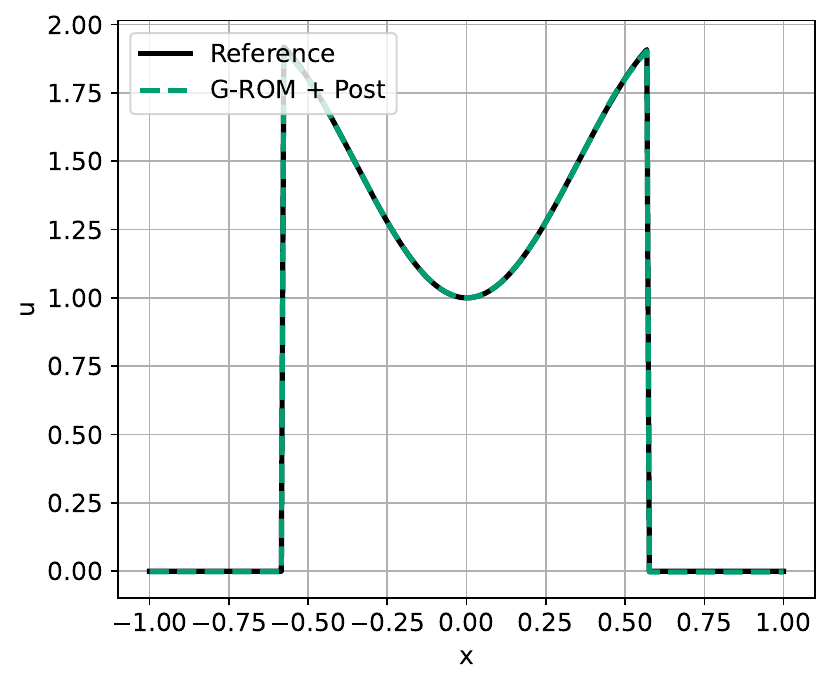}}
    \caption{Example 4. The Gegenbauer reconstruction of $ t = \pi / 4 $ profiles  along $ y = -0.6094 $, $ y = -0.4531$ and $ y =0.0 $.    }
    \label{eg4_2d_profiles_gegen}
\end{figure}

We then use Algorithm~\ref{alg_3} to improve the reconstruction. The corresponding results are presented in Figure~\ref{eg4_2d_gegen}. In this process, the threshold parameter in Algorithm~\ref{alg_3} is set to $\text{threshold}_1 = 0.35N$. It is clear that the numerical errors are reduced compared to the approach with fixed $x$ or $y$.
To better illustrate the effectiveness of Gegenbauer polynomial reconstruction in mitigating numerical oscillations, we also provide 1D cross-sectional plots in Figure ~\ref{eg4_2d_profiles_gegen}. Detailed numerical errors are presented in Table ~\ref{tab_eg4_2d}, showing that Gegenbauer polynomial reconstruction reduces both the relative and maximum errors by approximately an order of magnitude.

\begin{table}[!htbp]
    \centering
    \caption{\label{tab_eg4_2d}Example 4. Error comparison of the numerical solutions at $t=\pi / 4$.}
    \begin{tabular}{l|c|c}
        \hline
        & Maximum Error & Relative Error \\\hline
        G-ROM & 3.4717e-01 & 3.7396e-02 \\ \hline
        G-ROM + TV Regularization with $\lambda_{TV} = 0.2$ & 3.4498e-01 &2.2109e-02\\ \hline
        G-ROM + Reconstruction Based on Algorithm~\ref{alg_1} & 1.3202e-01 &1.4359e-02\\ \hline
        G-ROM + Reconstruction Based on Algorithm~\ref{alg_2} &   5.5341e-02 &  4.6207e-03\\ \hline
        G-ROM + Reconstruction Based on Algorithm~\ref{alg_3} &5.0860e-02 & 3.7432e-03\\ \hline
    \end{tabular}
\end{table}

\subsection{Example 5: 2D inviscid Burgers' equation}
\label{eg5}
We consider 2D Burgers' equation:
\begin{equation}
    u_t + uu_x+uu_y = 0 , \qquad (x,y)\in [0,4]\times[0,4],
\end{equation}
with periodic boundary condition in each direction and the initial condition 
\begin{equation}
    u(x,y,0) = u_0(x,y)= \frac{1}{2}+\sin(\frac{\pi}{2}(x+y)).
\end{equation}
We test the performance of Gegenbauer post-processing for the G-ROM solution in this example. The spatial domain is discretized into 201 uniformly distributed grid points in both $x$- and $y$-directions. To construct the snapshot matrix, the final time is set to $T = 2$, with 201 uniformly sampled snapshots. For each snapshot, the solution is computed using the method of characteristics. We apply the Gaussian filter with parameter $\sigma=1$ in this example. 

\begin{figure}[!htbp]
    \centering
    \subfloat[Before reconstruction]{\includegraphics[width=0.47\textwidth]{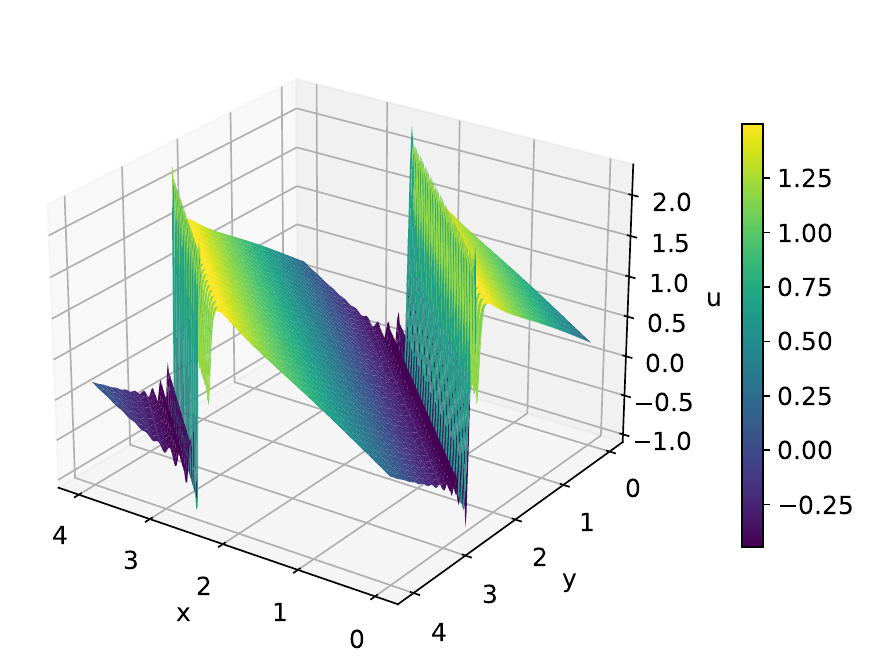}\label{eg5_2d_a}}
    \hspace{0.04\textwidth}
    \subfloat[After reconstruction]{\includegraphics[width=0.47\textwidth]{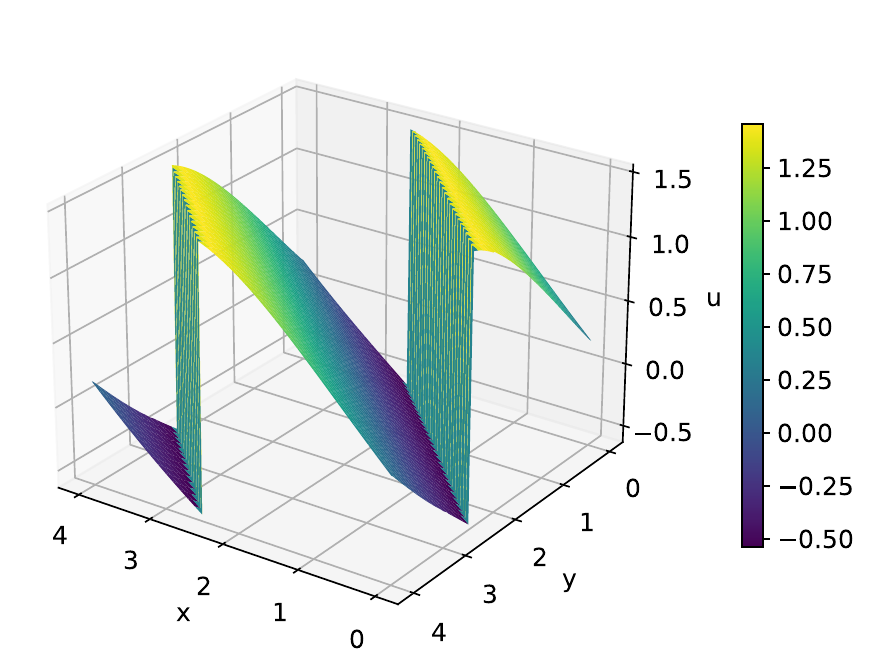}\label{eg5_2d_b}}
    \caption{Example 5. The G-ROM solution at $T=0.5$ before and after Gegenbauer reconstruction, with the reduced dimension of 30.}
    \label{eg5_2d_before_after}
\end{figure}

Fixing the reduced space dimension at 30, we employ the G-ROM method to compute the solution at $T=0.5$. 
Similar to the one-dimensional case, an exponential filter with $\alpha=2$ and $p=12$ is applied online for numerical stability.
The G-ROM solution is plotted in Figure ~\ref{eg5_2d_a}. It can be observed that significant numerical oscillations are present near the discontinuities. 

\begin{figure}[!htbp]
    \centering
    \subfloat[$y = 0.3980$]{\includegraphics[width=0.33\textwidth]{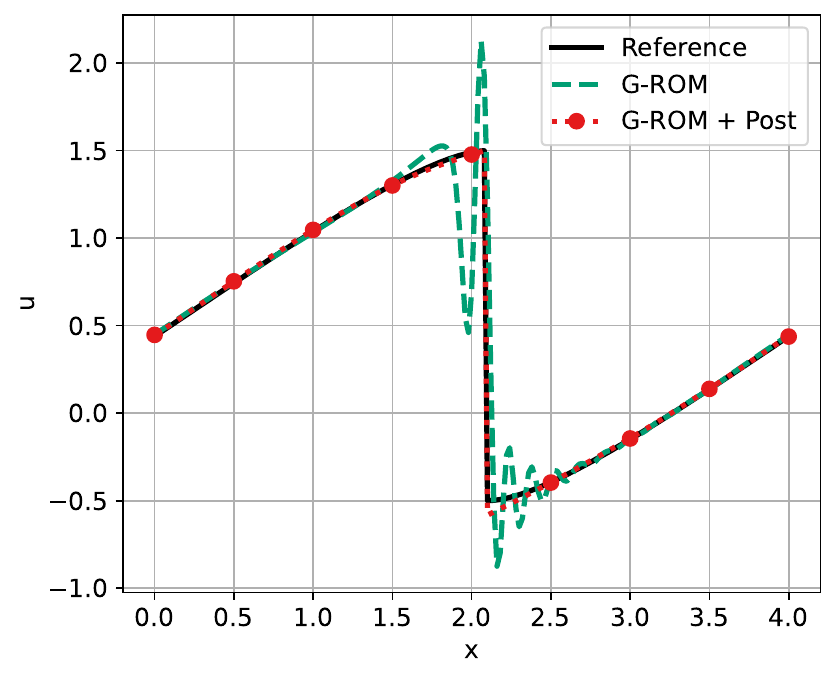}}
    \subfloat[$y = 1.3930$]{\includegraphics[width=0.33\textwidth]{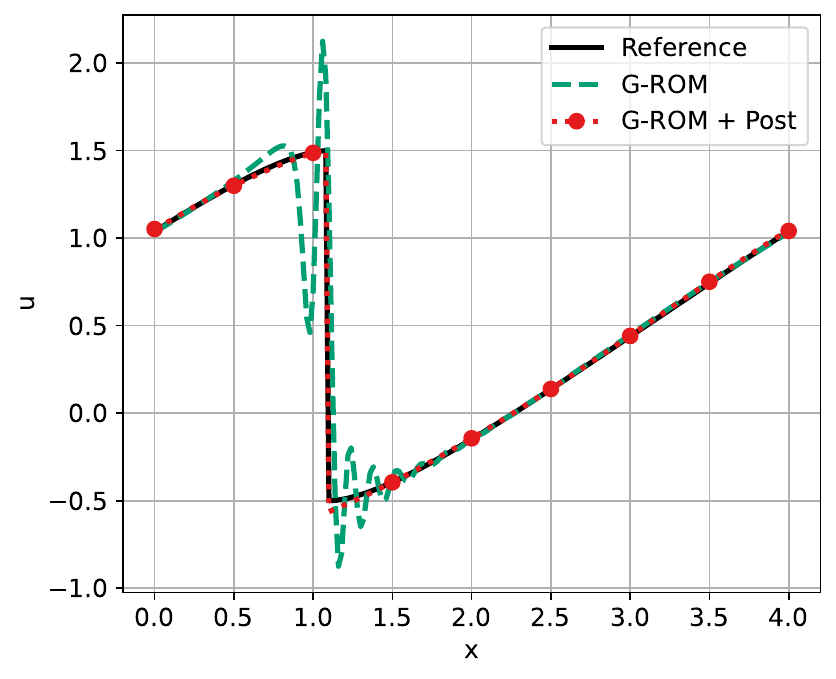}}
    \subfloat[$ y =2.3881$]{\includegraphics[width=0.33\textwidth]{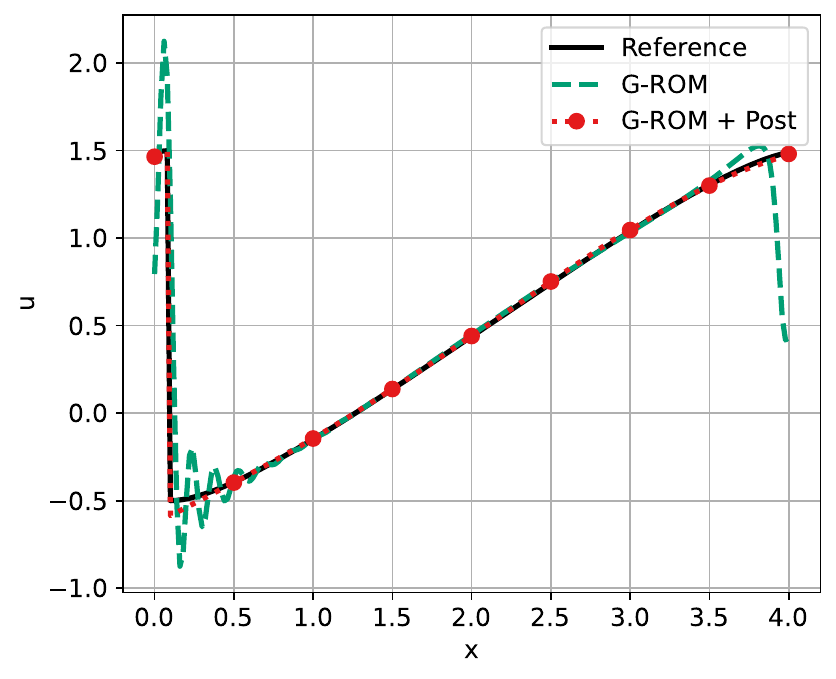}}
    \caption{Example 5. The Gegenbauer reconstruction of $ t =0.5$ profiles along $ y = 0.3980$, $ y = 1.3930$ and $ y =2.3881$.    }
    \label{eg5_2d_profiles_gegen}
\end{figure}

We then apply Gegenbauer reconstruction with parameter $\lambda = 2$ and $m = 4$
as in Algorithm~\ref{alg_3} with the threshold parameter set as $\textit{threshold}_1 = 0.1N$, where $N$ represents the number of grid points in one direction. The final reconstructed results are illustrated in Figure~\ref{eg5_2d_b}, demonstrating the overall suppression of numerical oscillations. We plot several 1D cross-sectional profiles in Figure~\ref{eg5_2d_profiles_gegen}, which shows the excellent performance of post-processing. The results are further validated in Table~\ref{tab_eg5_2d}, which shows the reduction in numerical error.

\begin{table}[!htbp]
    \centering
    \caption{\label{tab_eg5_2d}Example 5. Error comparison of the numerical solutions at $t=0.5$. }
    \begin{tabular}{l|c|c}
        \hline
        & Maximum Error & Relative Error \\\hline
        G-ROM & 2.1798e+00 &  2.4756e-01 \\ \hline
        G-ROM + TV Regularization with $\lambda_{TV} = 0.2$ & 2.2103e+00
        &2.2014e-01\\ \hline
        G-ROM + Reconstruction Based on Algorithm~\ref{alg_1} &2.0647e+00 &2.9213e-02\\ \hline
        G-ROM + Reconstruction Based on Algorithm~\ref{alg_2} &   2.0647e+00 & 2.9213e-02\\ \hline
        G-ROM + Reconstruction Based on Algorithm~\ref{alg_3} &9.0075e-02 &1.6067e-02\\ \hline
    \end{tabular}
\end{table}

\section{Conclusions and Future Work}
\label{Conclusions and Future Work}

This paper proposes a new approach to post-process the numerical solution from ROMs for mitigating spurious numerical oscillations for transport problems. We use Gegenbauer polynomial reconstruction method to re-project the numerical solution in each smooth segment of the solution, thereby suppressing the numerical oscillation and enhancing the accuracy at the same time.
We test our approach on both 1D and 2D transport-dominated problems and observe that the proposed procedure is effective for ROMs based on POD-Galerkin, OpInf and nonlinear reduction methods based on CAE. For inviscid problems, we observe accuracy enhancement of between one to two orders of magnitude. The Gegenbauer post-processing show clear advantages compared to TV regularization and ridge regression in obtaining a sharp resolution at the discontinuity. The experiments validate Gegenbauer reconstruction as a viable option to post-process ROM data.

Importantly, the method is not restricted to ROMs. 
Numerical experiments show that it also improves solutions from classical discretizations and machine learning solvers. In all cases, the post-processing either suppresses oscillations or sharpens diffused shocks, demonstrating that Gegenbauer reconstruction is a solver-independent and broadly applicable tool for nonlinear hyperbolic problems.

The main question that remains is the choice of reconstruction parameters $\lambda, m.$ While we made some heuristic observations about the optimal choice (e.g. $\lambda, m$ should not be chosen too large), we did not reach conclusive statements about the optimal choice of parameters which is clearly solution-dependent. 
In future work, we will investigate strategies to obtain the optimal parameter choice for the Gegenbauer reconstruction. One possible approach is to introduce additional regularization for selecting the parameter $m$, and it is also possible to obtain the optimal choice based on a supervised learning approach.

\section*{Acknowledgement}
\noindent 
Authors wish to acknowledge the significant contributions of Professor Yingda Cheng from Virginia Tech, who provided crucial input in the conceptualization, writing, and supervision of this work.


\begin{thebibliography}{99}

\bibitem{al2024rnn}
S. M. Al-Selwi, M. F. Hassan, S. J. Abdulkadir, A. Muneer, E. H. Sumiea, A. Alqushaibi, and M. G. Ragab. RNN-LSTM: From applications to modeling techniques and beyond—Systematic review. \textit{Journal of King Saud University-Computer and Information Sciences}, page 102068, 2024.

\bibitem{archibald2002method}
R. Archibald and A. Gelb. A method to reduce the Gibbs ringing artifact in MRI scans while keeping tissue boundary integrity. \textit{IEEE Transactions on Medical Imaging}, 21(4):305--319, 2002.

\bibitem{bergmann2009enablers}
M. Bergmann, C.-H. Bruneau, and A. Iollo. Enablers for robust POD models. \textit{Journal of Computational Physics}, 228(2):516--538, 2009.

\bibitem{berkooz1993proper}
G. Berkooz, P. Holmes, and J. L. Lumley. The proper orthogonal decomposition in the analysis of turbulent flows. \textit{Annual Review of Fluid Mechanics}, 25(1):539--575, 1993.

\bibitem{brunton2022data}
S. L. Brunton and J. N. Kutz. \textit{Data-driven science and engineering: Machine learning, dynamical systems, and control}. Cambridge University Press, 2022.

\bibitem{cagniart2017model}
N. Cagniart, R. Crisovan, Y. Maday, and R. Abgrall. Model order reduction for hyperbolic problems: a new framework. Working paper or preprint, 2017.

\bibitem{chaturantabut2010nonlinear}
S. Chaturantabut and D. C. Sorensen. Nonlinear model reduction via discrete empirical interpolation. \textit{SIAM Journal on Scientific Computing}, 32(5):2737--2764, 2010.

\bibitem{demers1992non}
D. DeMers and G. Cottrell. Non-linear dimensionality reduction. \textit{Advances in Neural Information Processing Systems}, 5, 1992.

\bibitem{driscoll2001pade}
T. A. Driscoll and B. Fornberg. A Pad\'e-based algorithm for overcoming the Gibbs phenomenon. \textit{Numerical Algorithms}, 26:77--92, 2001.

\bibitem{farcas2022filtering}
I. Farcas, R. Munipalli, and K. E. Willcox. On filtering in non-intrusive data-driven reduced-order modeling. In \textit{AIAA AVIATION 2022 Forum}, page 3487, 2022.

\bibitem{gelb2004parameter}
A. Gelb. Parameter optimization and reduction of round off error for the Gegenbauer reconstruction method. \textit{Journal of Scientific Computing}, 20:433--459, 2004.

\bibitem{gottlieb1997gibbs}
D. Gottlieb and C.-W. Shu. On the Gibbs phenomenon and its resolution. \textit{SIAM Review}, 39(4):644--668, 1997.

\bibitem{hesthaven2017numerical}
J. S. Hesthaven. \textit{Numerical methods for conservation laws: From analysis to algorithms}. SIAM, 2017.

\bibitem{hesthaven2007spectral}
J. S. Hesthaven, S. Gottlieb, and D. Gottlieb. \textit{Spectral methods for time-dependent problems}, volume 21. Cambridge University Press, 2007.

\bibitem{hochreiter1997long}
S. Hochreiter. Long short-term memory. \textit{Neural Computation MIT-Press}, 1997.

\bibitem{hoerl1970ridge}
A. E. Hoerl and R. W. Kennard. Ridge regression: Biased estimation for nonorthogonal problems. \textit{Technometrics}, 12(1):55--67, 1970.

\bibitem{issan2023predicting}
O. Issan and B. Kramer. Predicting solar wind streams from the inner-heliosphere to Earth via shifted operator inference. \textit{Journal of Computational Physics}, 473:111689, 2023.

\bibitem{jackiewicz2004determination}
Z. Jackiewicz. Determination of optimal parameters for the Chebyshev–Gegenbauer reconstruction method. \textit{SIAM Journal on Scientific Computing}, 25(4):1187--1198, 2004.

\bibitem{kashima2016nonlinear}
K. Kashima. Nonlinear model reduction by deep autoencoder of noise response data. In \textit{2016 IEEE 55th Conference on Decision and Control (CDC)}, pages 5750--5755, 2016.

\bibitem{kawai2022gegenbauer}
S. Kawai, W. Yamazaki, and A. Oyama. Gegenbauer reconstruction method with edge detection for multi-dimensional uncertainty propagation. \textit{Journal of Computational Physics}, 468:111505, 2022.

\bibitem{lee2020model}
K. Lee and K. T. Carlberg. Model reduction of dynamical systems on nonlinear manifolds using deep convolutional autoencoders. \textit{Journal of Computational Physics}, 404:108973, 2020.

\bibitem{li2025bayesian}
T. Li and A. Gelb. A Bayesian framework for spectral reprojection. \textit{Journal of Scientific Computing}, 102(3):78, 2025.

\bibitem{lu2020lagrangian}
H. Lu and D. M. Tartakovsky. Lagrangian dynamic mode decomposition for construction of reduced-order models of advection-dominated phenomena. \textit{Journal of Computational Physics}, 407:109229, 2020.

\bibitem{lu2021dynamic}
H. Lu and D. M. Tartakovsky. Dynamic mode decomposition for construction of reduced-order models of hyperbolic problems with shocks. \textit{Journal of Machine Learning for Modeling and Computing}, 2(1), 2021.

\bibitem{maulik2021reduced}
R. Maulik, B. Lusch, and P. Balaprakash. Reduced-order modeling of advection-dominated systems with recurrent neural networks and convolutional autoencoders. \textit{Physics of Fluids}, 33(3), 2021.

\bibitem{mcquarrie2021data}
S. A. McQuarrie, C. Huang, and K. E. Willcox. Data-driven reduced-order models via regularised operator inference for a single-injector combustion process. \textit{Journal of the Royal Society of New Zealand}, 51(2):194--211, 2021.

\bibitem{min2007fourier}
M. Min, S. Kaber, and W. Don. Fourier–Padé approximations and filtering for spectral simulations of an incompressible Boussinesq convection problem. \textit{Mathematics of Computation}, 76(259):1275--1290, 2007.

\bibitem{ohlberger2015reduced}
M. Ohlberger and S. Rave. Reduced basis methods: Success, limitations and future challenges. In \textit{Proceedings of ALGORITMY}, pages 1--12, 2016.

\bibitem{peherstorfer2020model}
B. Peherstorfer. Model reduction for transport-dominated problems via online adaptive bases and adaptive sampling. \textit{SIAM Journal on Scientific Computing}, 42(5):A2803--A2836, 2020.

\bibitem{peherstorfer2022breaking}
B. Peherstorfer. Breaking the Kolmogorov barrier with nonlinear model reduction. \textit{Notices of the American Mathematical Society}, 69(5):725--733, 2022.

\bibitem{peherstorfer2016data}
B. Peherstorfer and K. Willcox. Data-driven operator inference for nonintrusive projection-based model reduction. \textit{Computer Methods in Applied Mechanics and Engineering}, 306:196--215, 2016.

\bibitem{qian2022reduced}
E. Qian, I.-G. Farcas, and K. Willcox. Reduced operator inference for nonlinear partial differential equations. \textit{SIAM Journal on Scientific Computing}, 44(4):A1934--A1959, 2022.

\bibitem{reiss2018shifted}
J. Reiss, P. Schulze, J. Sesterhenn, and V. Mehrmann. The shifted proper orthogonal decomposition: A mode decomposition for multiple transport phenomena. \textit{SIAM Journal on Scientific Computing}, 40(3):A1322--A1344, 2018.

\bibitem{rudin1992nonlinear}
L. I. Rudin, S. Osher, and E. Fatemi. Nonlinear total variation based noise removal algorithms. \textit{Physica D: Nonlinear Phenomena}, 60(1-4):259--268, 1992.

\bibitem{shizgal2003towards}
B. D. Shizgal and J.-H. Jung. Towards the resolution of the Gibbs phenomena. \textit{Journal of Computational and Applied Mathematics}, 161(1):41--65, 2003.

\bibitem{sirisup2004spectral}
S. Sirisup and G. E. Karniadakis. A spectral viscosity method for correcting the long-term behavior of POD models. \textit{Journal of Computational Physics}, 194(1):92--116, 2004.

\bibitem{sobel19683x3}
I. Sobel and G. Feldman. A 3 $\times$ 3 isotropic gradient operator for image processing, a talk at the Stanford Artificial Project in 1968.

\bibitem{tadmor2007filters}
E. Tadmor. Filters, mollifiers and the computation of the Gibbs phenomenon. \textit{Acta Numerica}, 16:305--378, 2007.

\bibitem{wells2017evolve}
D. Wells, Z. Wang, X. Xie, and T. Iliescu. An evolve-then-filter regularized reduced order model for convection-dominated flows. \textit{International Journal for Numerical Methods in Fluids}, 84(10):598--615, 2017.

\bibitem{yu2022model}
J. Yu and J. S. Hesthaven. Model order reduction for compressible flows solved using the discontinuous Galerkin methods. \textit{Journal of Computational Physics}, 468:111452, 2022.


\end{thebibliography}
\end{document}